\documentclass[11pt]{article}
\usepackage[english]{babel}
\usepackage[cp850]{inputenc}
\usepackage[dvips]{graphicx}
\usepackage{amsmath,amsfonts,amsthm,amssymb}
\usepackage[usenames, dvipsnames]{color}
\usepackage{fancyhdr}
\usepackage{stmaryrd}
\usepackage[colorlinks=true,citecolor=red,linkcolor=blue,urlcolor=blue,pdfstartview=FitH]{hyperref}
\usepackage{dsfont}
\usepackage{xcolor}
\usepackage{epsfig}

\font\tenmath=msbm10 scaled 1200

\font\sevenmath=msbm7 scaled 1200
\font\fivemath=msbm5 scaled 1200

\newfam\mathfam \textfont\mathfam=\tenmath
\scriptfont\mathfam=\sevenmath \scriptscriptfont\mathfam=\fivemath

\def\R{{\mathbb R}}
\def\N{{\mathbb N}}
\def\E{{\mathbb E}}
\def\L{{\cal L}}

\def\P{{\mathbb P}}

\def\Q{{\mathbb Q}}

\def\D{{\mathbb D}}

\newtheorem{Thm}{Theorem}[section]
\newtheorem{Lem}{Lemma}[section]
\newtheorem{Pro}{Proposition}[section]
\newtheorem{Cor}{Corollary}[section]
\newtheorem{Dfn}{Definition}[section]

\newfam\mathfam \textfont\mathfam=\tenmath
\scriptfont\mathfam=\sevenmath \scriptscriptfont\mathfam=\fivemath

\def \^#1{\if#1i{\accent"5E\i}\else{\accent"5E#1}\fi}

\def \D {I\!\!D}

\def \cqfd{\quad\Box}
\def \ms{\medskip}

\def \bs{\bigskip}
\def \ni{\noindent}

\def\ni{\noindent}

\def \supess{{\rm esssup}}

\title{\bf Convex order for path-dependent derivatives: a dynamic programming approach}

\author{\textsc{Gilles Pag\`es} \thanks{Laboratoire de Probabilit\'es et Mod\`eles al\'eatoires, UMR~7599, UPMC, case 188, 4, pl. Jussieu, F-75252 Paris Cedex 5, France. E-mail: \texttt{gilles.pages@upmc.fr}}}

\begin{document}

\maketitle




\begin{abstract} We investigate the (functional) convex order of for various continuous martingale processes, either with respect to their diffusions coefficients for L\'evy-driven SDEs  or their integrands for stochastic integrals. Main  results are bordered by counterexamples.   Various upper and lower bounds can be derived for   pathwise European option prices in local volatility models. In view of numerical applications, we adopt a systematic (and symmetric) methodology: (a) propagate the convexity in a {\em simulatable} dominating/dominated discrete time model through a backward induction (or linear dynamical principle); (b) Apply  functional weak convergence results to numerical schemes/time discretizations of the continuous time martingale satisfying (a) in order  to transfer the   convex order properties. Various bounds are derived for European options written on convex pathwise dependent payoffs.  We retrieve and  extend former results  obtains by several authors (\cite{NEKJESH, BRusch1, HOB,YORetal}) since the seminal  paper~\cite{HAJ} by Hajek. In a second part, we extend this approach to Optimal Stopping problems using a that the Snell envelope satisfies  (a') a  Backward Dynamical Programming Principle to propagate convexity in discrete time; (b') satisfies  abstract convergence results under non-degeneracy assumption on filtrations. Applications to the comparison of American option prices  on convex pathwise payoff processes are given obtained by a purely probabilistic arguments.
     \end{abstract}


\bs 
\ni {\em 2010 AMS Classification.} Primary: 62 P05, 60E15, 91B28, secondary : 60J75, 65C30

%
%

\section{Introduction}\label{Intro}

The first aim of this paper is to propose a systematic and unified approach to establish {\em functional   convex order} results for discrete and continuous time martingale stochastic processes using the propagation of convexity through some kind of backward dynamic programing principles  (in discrete time) and weak functional limit theorems (to switch to continuous time. The term ``functional" mainly refers to the ``parameter" we deal with: thus, for diffusions processes (possibly with jumps) this parameter is the  diffusion coefficient or, for stochastic integrals, their  integrand. Doing so we will retrieve, extend and sometimes establish new results on functional convex order. As a second step, we will tackle the same type of question  in the framework of Optimal Stopping Theory for the  Snell envelopes and their means,  the so-called {\em r\' eduites} (which maybe provides a better justification for the terminology ``dynamic programming approach" used in the title).

Let us first briefly recall  that if $X$ and $Y$ are two real-valued random variables, $X$ is dominated by $Y$ for the convex order~--~ denoted $X\preceq_c Y$~--~if, for every convex functions $f:\R\to \R$ such that $f(X),\, f(Y)\!\in L^1(\P)$,  
\[
\E\, f(X)\le \E\, f(Y).
\]
 Thus, if $(M_{\lambda})_{\lambda>0}$ denotes  a martingale indexed by a parameter $\lambda$, then $\lambda\mapsto M_\lambda$ is non-decreasing for the convex order as a straightforward consequence of Jensen's Inequality. The converse is clearly not true but, as first established by Kellerer in~\cite{KEL}, whenever  $\lambda\mapsto X_\lambda$ is non-decreasing for the convex order, there exists a martingale $(\widetilde X_\lambda)_{\lambda \ge 0}$ such that $X_\lambda\stackrel{d}{=}\widetilde X_\lambda$ for every $\lambda \ge 0$ (we will say that $(X_{\lambda})$ and $(\widetilde X_{\lambda})$ coincide en $1$-marginal distributions.

The connection with Finance and, to be more precise with the pricing and hedging of derivative products is straightforward~: let $(X_t^{(\theta)})_{t\in [0,T]}$ be a family of non-negative $\P$-martingales on a probability space $(\Omega,{\cal A}, \P)$ indexed by a parameter $\theta$. Such a family can be seen as  possible models for  the discounted  price dynamics of a risky asset under its/a  risk-neutral probability where $\theta$ (temporarily) is  a real parameter ($e.g.$ representative of the volatility). If $\theta \mapsto X^{(\theta)}_{_T}$   is non-decreasing for the convex order, then for every convex {\em vanilla payoff} function $f: \R_+\to \R_+$, the function $\theta \mapsto \E\, f(X^{(\theta)}_{_T})$ is non-decreasing or equivalently its {\em greek} $\frac{\partial }{\partial \theta} \E f(X^{(\theta)}_{_T})$  with respect to $\theta$  is non-negative. Typically, in a discounted Black-Scholes model
\[
X^{\sigma,x}_t = x e^{\sigma W_t-\frac{\sigma^2}{2}t}, \,x,\,\sigma>0,
\]
the function $\sigma\longmapsto \E f\big( x \,e^{\sigma W_{_T}-\frac{\sigma^2T}{2}}\big)$ since 
\[
 \forall\, t\!\in [0,T], \quad x \,e^{\sigma W_{_T}-\frac{\sigma^2T}{2}}\stackrel{\cal L}{\sim}\Big[e^{B_u-\frac u2}\Big]_{|u= \sigma^2T}
\]
where $u\mapsto e^{B_u-\frac u2}$ is a martingale as as well as its composition with $\sigma\mapsto \sigma^2T$. So $(X^{\sigma,x}_{_T})_{\sigma \ge 0}$coincides in $1$-marginal distributions with a martingale.  The same result holds true for the premium of convex {\em Asian payoff} functions of the form
\[
 \E\, f\left(\frac 1T\int_0^T  x e^{\sigma W_t-\frac{\sigma^2t}{2}}dt\right)
\]
but, by contrast, its proof is significantly more involved (see~\cite{Carretal} or, more recently, the proof in~\cite{YORetal} where an explicit martingale based on the Brownian sheet coinciding in $1$-dimensional martingale is exhibited). Both results turn out to be  examples  of a general result dealing with convex pathwise dependent functionals (see $e.g.$~\cite{YORetal} or~\cite{PAG} where a functional co-monotony argument is used).

A natural question  at this stage is to try establishing  a {\em functional version} of these results in terms of $\theta$-parameter $i.e.$ when $\theta$ is longer a real number or a vector  but lives  in a subset of a functional space or even of space of stochastic processes. A typical example where $\theta$ is a function is the case where $X^{(\theta)}$ is a diffusion process, (weak) solution to a Stochastic Differential equation $(SDE)$ of the form
\[
dX^{(\theta)}_t = \theta(t,X^{(\theta)}_{t-}) dZ_t, \;X^{(\theta)}_0= x, \; t\!\in [0,T], 
\]
with  $Z=(Z_t)_{t\in [0,T]}$   a martingale L\'evy process (having moments of order at least $1$). The parameter  $\theta$ can also be a (predictable) stochastic process when
\[
X^{(\theta)}_t = \int_0^t \theta_sdZ_s, \; t\!\in [0,T]. 
\]

When dealing with optimal stopping problems, $i.e.$ with the {\em r\'eduite} of a target process $Y_t = F(t,X^{(\theta),t})$, $t\!\in [0,T]$,  (where $X^{(\theta),t}_s=X^{(\theta)}_{s\wedge t}$ is the stopped process $X^{(\theta)}$ at $t$) and all  the functionals $F(t,.)$ are (continuous) convex functionals defined on the path space of the process $X$, the {\em functional convex order} as defined above amounts to determine the sign of the sensitivity with respect to the functional parameter $\theta$ of an American option with payoff functional $F(t,.)$ at time $t\!\in [0,T]$, ``written" on $X^{(\theta)}$: if the holder of the American option contract exercises the option at time $t$, she receives $F(t,X^{(\theta),t})$.

More generally, various notions of convex order in Finance are closely related to risk modeling and come out in many other frameworks than the pricing and hedging of derivatives.

\smallskip
Many  of these questions have already been investigated  for a long time: thus,  the first result known to us goes back to Hajek in~\cite{HAJ} where convex order is established for Brownian martingale diffusions ``parametrized" by their (convex) diffusion coefficients (with  an extension to drifted diffusions  with non-decreasing convex  drifts but with a restriction to  non-decreasing convex functions   $f$ of $X_{_T}$).  The first application to the sensitivity of (vanilla) options of  both European and American style, is due to~\cite{NEKJESH}.  It is shown that the options with convex payoffs in a $[\sigma_{\min}, \sigma_{\max}]$-valued local volatility model with bounded volatility can be upper- and lower-bounded by the  premium of the the same option contracts in a Black-Scholes model with volatilities $\sigma_{\min}$ and $\sigma_{\max}$ respectively (note however that a PDE approach relying on a maximal principle provides an alternative easier proof). See also~\cite{HOB} for a result on lookback options. More recently, in a series of papers (see~\cite{BRusch0, BRusch1, BRusch2})  Bergenthum and R\" uschendorf extensively investigated the above mentioned problems   (for both fixed maturity  and for optimal stopping problemss) for various classes of continuous and jump processes, including general semi-martingales in~\cite{BRusch0} (where the  comparison is carried out  in terms of their  predictable local characteristics, assuming  one of them  propagates convexity, then proving this last fact). In several of these papers, the  convexity is  --~but not always~(see~\cite{BRusch1})~--   propagated directly  in continuous time  which is clearly an elegant way to proceed but also  seems to more heavily rely on specific features of the investigated class of processes (see~\cite{YORetal}). In this paper, we  propose a generic and systematic systematic --~but maybe also more  ``symmetric"~--  two-fold approach which turns out to be efficient for many classes of  stochastic dynamics  and processes which is based on a swathing from discrete times to continuous time using functional weak limit theorems ``\`a la Jacod-Shyriaev" (see~\cite{JASH2}). To be more precise:

\smallskip
-- As first step, we study the propagation of convexity ``through" a discrete time dynamics -- typically a $GARCH$ model~-- in a very elementary way for path-dependent convex functionals  relying on repeated elementary backward inductions and conditional Jensen's inequality. These inductions take advantage of the ``linear" backward dynamical programming principle resulting from from a discrete time  martingale property written in a  step-by-step manner. This terminology borrowed from stochastic control can be viewed as a bit excessive but refers to a second aspect of the paper devoted to optimal stopping theory  (see further on).

\smallskip
-- As a second step, we use  that these discrete time $GARCH$ model are   discretization schemes for the ``target" continuous time dynamics (typically the ``Euler schemes as concerns diffusion processes) and we transfer to this target the searched functional convex order property by calling upon functional limit theorems for the convergence of stochastic integrals (typically borrowed form~\cite{JaMePa} and/or~\cite{KUPR}.

\smallskip Our typical result for jump diffusions reads as follows (for a more complete and  rigorous statements see~$e.g.$ Theorems~\ref{FCOEuro} and \ref{FCOEuroJump} in Section~\ref{sec:FCO}). If $0\le \kappa_1\le \kappa\le \kappa_2$ are continuous functions with linear growth defined on $\R$ and $\kappa$ is convex then the existing weak solutions $X^{(\kappa_i)}$, $i=1,2$, to the $SDE$s 
\[
X^{(\kappa_i)}_t= x+\int_{(0,t]} \kappa_i(X^{(\kappa_i)}_{s-})dZ_s
\]  
where $Z=(Z_t)_{t\in [0,T]}$ is a martingale L\'evy process with L\'evy measure $\nu$ satisfying $\nu(z^2)<+\infty$, then $X^{(\kappa_1)}\preceq_{fc} X^{(\kappa_2)}$ for the convex order defined on (continuous for the Skorokhod topology) convex functionals (with polynomial growth). Note that when $Z$ is a Brownian motion, the continuity of the functional appears as a consequence of its convexity (under the polynomial growth assumption, see the remark in Section~\ref{sec:FCOW}). Equivalently, we have $X^{(\kappa_1)}\preceq_{fc} X^{(\kappa)} \preceq_{fc} X^{(\kappa_2)}$ as soon as both functions $\kappa_i$ are convex. Results in the same spirit are obtained for stochastic integrals, Dol\'eans exponentials  (which unfortunately requires   one of the  two integrands $H_1$ and $H_2$  to be  deterministic). Counter-examples to put the main results in perspective are exhibited to prove the consistency of these assumptions in both settings.

\smallskip We also deal with non-linear problems, typically optimal stopping problems, framework where we use the same approach from discrete to continuous time, taking advantage of the Backward Dynamic Programming Principle in the first framework and using various convergence results for the Snell envelope~(see~\cite{LAPA}).  In fact, a   similar approach in discrete time has already been been developed to solve the propagation of convexity in a stochastic control problem ``through"  the dynamic programming principle  in a pioneering work by Hern\'andez-Lerma and Runggaldier~~\cite{HERU}.
%

\smallskip
 The main reason for developing in a systematic manner this approach is related with Numerical Probability: our discrete time models appear as {\em simulatable} discretization schemes of the continuous time dynamics of interest. It is important for applications, especially in Finance, to have at hand discretization schemes which {\em both} preserve the (functional) convex order and can be simulated at a reasonable cost. So is the case of the Euler scheme for L\'evy driven diffusions (as soon as the underling L\'evy processes is  itself simulatable). This is not always the case: think $e.g.$ to the (second order) Milstein scheme for Brownian diffusions, in spite of its better performances in term of strong convergence rate.
 
\medskip The paper is organized as follows. Section~\ref{sec:FCO} is devoted to functional convex order for path-dependent functionals of Brownian and L\'evy driven martingale diffusion processes. Section~\ref{sec:ItoProc} is devoted to comparison results for It\^o processes based on comparison of their integrand.  Section~\ref{sec:Amer} deals with  {\em r\'eduites}, Snell envelopes of path dependent obstacle processes  (American options)  in both Brownian and L\'evy driven martingale diffusions. In the two-fold appendix, we provide short proofs of functional weak convergence of the Euler scheme toward a weak solution of $SDE$s   in both Brownian and L\'evy frameworks under natural continuity-linear growth  assumptions on the diffusion coefficient.

\medskip
 \noindent {\sc Notation:} $\bullet$ For every $T>0$ and  every integer $n\ge 1$, one denotes the uniform mesh of $[0,T]$ by $t^n_k= \frac{kT}{n}$, 
$k=0,\ldots,n$. Then for every $t\!\in [\frac{kT}{n}, \frac{(k+1)T}{n})$, we set $\underline{t}_n =\frac{kT}{n}$ and $\overline t^n=\frac{(k+1)T}{n}$ with the convention $\underline T_n=T$. We also set $\underline t_{n-}= \lim_{s\to t} \underline s_n = \frac{kT}{n}$ if $t\!\in\big(\frac{kT}{n}, \frac{(k+1)T}{n}\big]$. 

\noindent $\bullet$ For every $u=(u_1,\ldots,u_d)$, $v=(v_1,\ldots,v_d)\!\in\R^d$, $(u|v)= \sum_{i=1}^d u_iv_i$, $|u|= \sqrt{(u|u)}$ and $x_{m:n}=(x_m,\ldots,x_n)$ (where $m\le n$, $m$, $n\!\in \N\setminus\{0\}$).

\noindent $\bullet$ ${\cal F}([0,T], \R)$ denotes the $\R$-vector space of $\R$-valued functions $f:[0,T]\to\R$ and ${\cal C}([0,T], \R)$ denotes the sub\-space of $\R$-valued continuous functions defined over $[0,T]$.

\noindent $\bullet$ For every $\alpha\!\in {\cal F}([0,T], \R)$, we define ${\rm Cont}(\alpha)= \big\{t\!\in [0,T]\, :\, \alpha \mbox{ is continuous at } t\big\}$ with the usual left- and right- continuity conventions at $0$ and $T$ respectively. We also define the {\em uniform continuity modulus} of $\alpha$ by where $\displaystyle w(\alpha,\delta)= \sup\big\{|\alpha(u)-\alpha(v)|,\,u,v\in [0,T], |u-v|\le \delta\big\}$ ($\delta\!\in [0,T]$).

\noindent $\bullet$ $L^p_{_T}= L^p([0,T],dt)$, $1\le p\le +\infty$, $|f|_{L^p_{_T}}= \big(\int_0^T |f(t)|^pdt\big)^{\frac 1p}\le +\infty$, $1\le p<+\infty$ and $|f|_{L^{\infty}_{_T}}= dt\mbox{-}{\supess}|f|$ where $dt$ stands for the Lebesgue measure on $[0,T]$ equipped with its Borel $\sigma$-field.

\noindent $\bullet$ For a function  $f:[0,T]\to \R$, we denote $\|f\|_{\sup}= \sup_{u\in [0,T]}|f(u)|$.

\noindent $\bullet$ Let $(\Omega,{\cal A}, \P)$ be a probability space and let $p\!\in (0,+\infty)$. For every random vector $X:(\Omega,{\cal A})\to \R^d$ we set $\|X\|_p = \big(\E|X|^p\big)^{\frac 1p}$. $L^p_{\R^d}(\Omega,{\cal A}, \P) $ denotes the vector space of (classes) of $\R^d$-valued random vectors $X$ such that $\|X\|_p<+\infty$. $\|\,.\,\|_p$ is a norm on $L_{\R^d}^p(\Omega,{\cal A}, \P)$  for $p\!\in[1,+\infty)$ (the mention of $\Omega$, ${\cal A}$ and the subscript $_{\R^d}$ will be dropped when there is no ambiguity).

\noindent $\bullet$ If ${\cal F}=({\cal F}_t)_{t\in [0,T]}$ denotes a filtration on  $(\Omega,{\cal A}, \P)$,  let ${\cal T}^{\cal F}_{[0,T]}=\{\tau:\Omega\to [0,T], {\cal F}\mbox{-stopping time}\}$.

\noindent $\bullet$  ${\cal F}^Y=({\cal F}^Y_t)_{t\in [0,T]}$ is the smallest right continuous  filtration $({\cal G}_t)_{t\in [0,T]}$ that makes the process $Y=(Y_t)_{t\in [0,T]}$ a $({\cal G}_t)_{t\in [0,T]}$-adapted process.

\noindent $\bullet$  $\D([0,T],\R^d)$ denotes the set of $\R^d$-valued right continuous left limited (or c\`adl\`ag following the French acronym) function defined on the interval $[0,T]$, $T>0$.
It is usually endowed with the  Skorokhod topology denoted  $Sk$ (see~\cite{JASH}, chapter VI or~\cite{BIL}, chapter~3, for an introduction to Skorokhod topology). 

\noindent $\bullet$ If two random vectors $U$ and $V$ have the same distribution, we write $U\stackrel{d}{\sim} V$. If an $(S,d_S)$-valued sequence of random variable ($(S,d)$ Polish space  equipped with its Borel $\sigma$-field ${\cal B}or(S)$) {\em weakly converges} toward an $(S,d)$-valued random variable $Y_{\infty}$ (we will also say {\em converge in distribution or in law}), we will denote $\displaystyle Y_n \stackrel{{\cal L}(S,d_S)}{\longrightarrow} Y_{\infty}$ or, if no ambiguity, $\displaystyle Y_n \stackrel{{\cal L}(d_{S})}{\longrightarrow} Y_{\infty}$ . 

\smallskip
We will extensively make use the following classical result:

\smallskip
Let  $(Y_n)_{n\ge 1}$be  a sequence of tight random variables taking values in a Polish space $(S, d_{_S})$ (see~\cite{BIL}, Chapter~1). If $Y_n $ weakly converges toward $Y_{\infty}$ and $(\Phi(Y_n))_{n\ge 1}$ is uniformly integrable where $\Phi:S\to \R$ is a  Borel function 
then, for every $\P_{Y_{\infty}}$-$a.s.$ continuous Borel  functional $F:S\to \R$   such that $|F(u)|\le C(1+\Phi(u))$ for every $u\!\in S$, one has $\E \,F(Y_n)\to \E \,F(Y_{\infty})$. 

\section{Functional convex order}\label{sec:FCO}
\subsection{Brownian martingale diffusion}\label{sec:FCOW}

The main result of this section is the proposition below.
\begin{Thm}\label{FCOEuro}Let $\sigma$, $\theta:[0,T]\times \R\to \R$ be two continuous functions with linear growth in $x$ uniformly in $t\!\in [0,T]$. Let $X^{(\sigma)}$ and $X^{(\theta)}$ be two Brownian martingale diffusions, supposed to be the unique {\em weak} solutions starting from $x$ at time $0$, to the stochastic differential equations (with $0$ drift)
\begin{equation}\label{eq:martdiff}
dX^{(\sigma)}_t = \sigma(t,X^{(\sigma)}_{t}) dW^{(\sigma)}_t, \; X^{(\sigma)}_0=x\quad \mbox{ and }\quad dX^{(\theta)}_t = \theta(t,X^{(\theta)}_{t}) dW^{(\theta)}_t, \; X^{(\theta)}_0=x
\end{equation}
respectively, where $W^{(\sigma)}$ and $W^{(\theta)}$ are standard one dimensional Brownian motions.

\smallskip
\noindent $(a)$ {\em Partitioning assumption}: Let  $\kappa:[0,T]\times \R\to \R_+$ be  a continuous function with (at most)  linear growth in $x$ uniformly in $t\!\in[0,T]$, satisfying
$$
\kappa(t,.)\; \mbox{ is convex for every $t\!\in [0,T]$ and }\; 0\le \sigma\le \kappa\le \theta.
$$
 Then, for every convex  functional $F : {\cal C}([0,T],\R)\to\R$
 with  $(r,\|\,.\,\|_{\sup})$-polynomial growth, $r\ge 1$, in the following sense
\[
\forall\, \alpha\!\in {\cal C}([0,T],\R),\quad |F(\alpha)|\le C(1+\|\alpha\|^r_{\sup}) 
\]
 one has
\[
\E\,F(X^{(\sigma)})\le \E\,F(X^{(\theta)}).
\]
From now on, the function $\kappa$ is called a {\em partitioning function}. 

\smallskip
\noindent $(a')$ Claim $(a)$ can be reformulated equivalently as follows: if   either $\sigma(t,.)$ is convex for every $t \!\in [0,T]$ or  $\theta(t,.)$ is convex for every $t \!\in [0,T]$ and $0\le \sigma\le \theta$, then  the conclusion of $(a)$ still holds true.

\smallskip
\noindent $(b)$ {\em Domination assumption}:  If $|\sigma| \le \theta$ and   $\theta$ is  convex, then 
\[
\E\,F(X^{(\sigma)}) \le \E\,F(X^{(\theta)}).
\]
\end{Thm}

\noindent {\bf Remarks.}  $\bullet$ The linear growth assumption on the convex functional $F$ implies its everywhere local boundedness  on the  Banach space $\big( {\cal C}([0,T],\R),\|\,.\,\|_{\rm sup}\big)$, hence  its $\|\,.\,\|_{\rm sup}$-continuity (see $e.g.$ Lemma~2.1.1 in~\cite{LUC}, p.22).

\smallskip
\noindent $\bullet$ The introduction of two standard Brownian motions $W^{(\sigma)}$ and $W^{(\theta)}$ in the above claim~$(a)$ is just a way to recall that the two diffusions processes can be defined on different probability spaces, although it  may be considered as an abuse of notation. By ``unique weak solutions", we mean classically  that two such solutions (with respect to possibly different Brownian motions) share  the same distribution on the Wiener space.

\smallskip
\noindent $\bullet$ Weak uniqueness holds true as soon as strong uniqueness holds $e.g.$ as soon as $\sigma$ and $\theta$ are Lipschitz continuous in $x$, uniformly in $t\in [0,T]$, (as it can easily be derived from Theorem A.3.3, p.271, in~\cite{BOLE}).
%

\medskip
The proof of  this theorem can be decomposed in two main steps: the first one is   a dynamic programming approach in discrete time detailed in Proposition~\ref{FCOEuroDiscTime} below which relies itself on a revisited version of Jensen's Inequality. The second one remiss on a functional  weak approximation argument. 
%

\medskip
The first ingredient is a simple reinterpretation of the celebrated Jensen Lemma.

%
%
%
%
%
%

\begin{Lem}\label{Jen}(Revisited Jensen's Lemma) \label{Lem:Jensen}Let $Z:(\Omega,{\cal A}, \P)\to \R$ be an integrable centered $\R$-valued random vector.

\smallskip 
\noindent $(a)$ Assume that $Z\!\in L^r(\P)$ for an  $r\ge 1$. For every Borel function $\varphi:\R\to \R$ such that $|\varphi(x)|\le C(1+|x|^{r})$, $x\!\in \R$, we define
\begin{equation}\label{eq:Q}
\forall\, u\!\in \R,\; Q\varphi(u)= \E\, \varphi\big(uZ\big).
\end{equation}
If $\varphi$ is convex, then, $Q\varphi$  is convex and  $u\mapsto Q\varphi(u)$ is non-decreasing  on $\R_+$, non-increasing on  $\R_-$. 

\smallskip 
\noindent $(b)$ If $Z$ has exponential moments in the sense that
\[
\forall\, u\!\in \R,\; \E (e^{uZ})<+\infty
\]
(or equivalently  $\E( e^{a|Z|}) <+\infty$ for every $a\ge 0$), then claim~$(a)$ holds true for any convex function $\varphi:\R\to\R$ satisfying an  exponential growth condition of the form $|\varphi(x)|\le Ce^{C|x|}$, $x\!\in \R$, for a real constant $C\ge 0$.

\smallskip 
\noindent $(c)$ If $Z$ has a symmetric distribution  (i.e. $Z$ and $-Z$ have the same distribution) and $\varphi:\R\to \R$ is convex, then $Q\varphi$ is an even function, hence satisfying the following {\em  maximum principle}:  
\[
\forall\, a\!\in \R_+,\quad  \sup_{|u|\le a}Q\varphi(u)= Q\varphi(a).
\] 
\end{Lem}

\noindent {\bf Proof.} $(a)$-$(b)$ Existence and convexity of $Q\varphi$ are obvious. 
 The function $Q\varphi$ is clearly finite on $\R$ and convex. Furthermore, Jensen's Inequality implies that  
\[
Q\varphi(u) = \E\,\varphi (u \,Z) \ge \varphi (\E\,u \,Z)= \varphi(0) = Q\varphi(0)
\]
since $Z$ is centered. Hence $Q\, \varphi$ is convex and minimum at $u=0$ which implies that it is non-increasing on $\R_-$ and non-decreasing on $\R_+$.

\smallskip
\noindent $(c)$ is obvious.$\cqfd$

\begin{Pro}\label{FCOEuroDiscTime}Let $(Z_k)_{1\le k\le n}$ be a sequence of independent, centered, $\R$-valued random vectors lying in  $L^r(\Omega, {\cal A}, \P)$, $r\ge 1$, and let  $({\cal F}^Z_k)_{0\le k\le n}$ denote its natural filtration. Let $(X_k)_{0\le k\le n}$ and $(Y_k)_{0\le k\le n}$ be two  sequences of random vectors recursively defined by
\begin{equation}\label{modele}
X_{k+1}= X_k + \sigma_k(X_k) Z_{k+1}, \; Y_{k+1}= Y_k + \theta_k(Y_k) Z_{k+1}, \; 0\le k\le n-1,\; X_0=Y_0=x
\end{equation}
where $\sigma_k$, $\theta_k:\R \to \R$, $k=0,\ldots,n-1$,  are Borel functions   with linear growth $i.e.$ $|\sigma_k(x)|+|\theta_k(x)|\le C(1+|x|)$, $x\!\in \R$, for a real constant $C\ge 0$.

\smallskip
\noindent $(a)$  Assume that, either $\sigma_k$ is convex for every $ k\!\in \{0, \ldots,n-1\}$, or $\theta_k$ is convex for every $ k\!\in \{0, \ldots,n-1\}$, and that
\[
\forall\, k\!\in \{0, \ldots,n-1\},\,\; 0\le \sigma_k \le \theta_k.
\]
Then, for every convex function $\Phi:\R^{n+1}\to \R$ with $r$-polynomial growth, $r\ge1$,    $i.e.$ satisfying  $|\Phi(x)|\le C(1+|x|^{r})$, $x\!\in \R$, for a real constant $C\ge 0$, 
\[
\E\, \Phi(X_{0:n})\le \E\, \Phi(Y_{0:n}).
\] 

\noindent$(b)$ If the random variables $Z_k$ have symmetric distributions, if the functions $\theta_k$ are all convex and if
\[
\forall\, k\!\in \{0, \ldots,n-1\},\; | \sigma_k |\le \theta_k,
\]
then the conclusion of claim $(a)$ remains valid.
\end{Pro}


\noindent {\bf Proof.} $(a)$ First  one shows by an easy induction that the random variables $X_k$ and $Y_k$ all lie in $L^r$. Let $Q_k$, $k=1,\ldots,n$, denote the operator attached to $Z_k$ by \eqref{eq:Q} in Lemma~\ref{Jen}.

Then, one defines the following martingales
\[
M_k=\E \big(\Phi(X_{0:n})\,|\, {\cal F}^Z_k\big)\quad \mbox{ and }\quad N_k=\E\big( \Phi(Y_{0:n})\,|\, {\cal F}^Z_k\big), \; 0\le k\le n.
\]
Their existence follows from the growth assumptions on $\Phi$, $\sigma_k$ and $\theta_k$, $k=1,\ldots,n$. Now we define recursively in a backward way two sequences of functions $\Phi_k$ and $\Psi_k:\R^{k+1}\to \R$, $k=0,\ldots,n$, 
\[
\Phi_n =\Phi \; \mbox{ and }\;\Phi_k (x_{0:k})=\big(Q_{k+1}\Phi_{k+1}(x_{0:k}, x_k+.)\big)(\sigma_k(x_k)),\;  x_{0:k}\!\in \R^{k+1},\; k=0,\ldots,n-1, 
\]
on the one hand and, on the other hand, 
\[
\Psi_n =\Phi \; \mbox{ and }\;\Psi_k (x_{0:k})=\big(Q_{k+1}\Psi_{k+1}(x_{0:k}, x_k+.)\big)(\theta_k(x_k)),\; x_{0:k}\!\in \R^{k+1},\; k=0,\ldots,n-1.
\]
This can be seen as a {\em linear Backward Dynamical  Programming Principle}. It is clear by a (first) backward induction and the definition of the operators $Q_k$ that,  for every $k\!\in \{0,\ldots,n\}$, 
\[
M_k=\Phi_k(X_{0:k})\quad \mbox{ and }\quad N_k=\Psi(Y_{0:k}).
\]
Let $k\!\in \{0,\ldots,n-1\}$. One derives from  the properties of the operator $Q_{k+1}$ (and the  representation below as an expectation)  that, for any convex function $G:\R^{k+2}\to \R$ with $r$-polynomial growth, $r\ge 0$, the function 
\begin{equation}\label{eq:Gtilde}
\widetilde G: (x_{0:k},u)\longmapsto \big(Q_{k+1}G(x_{0:k},x_k+\,.\,)\big)(u)= \E \,G(x_{0:k},x_k+uZ_{k+1})
\end{equation}
 is convex. Moreover, owing to Lemma~\ref{Lem:Jensen}$(a)$, for fixed $x_{0:k}$,  $\widetilde G$ is non-increasing  on $(-\infty,0)$, non-decreasing  on $(0,+\infty)$ as a function of $u$.  In turn, this implies that, if $\gamma:\R\to \R_+$ is convex (and non-negative), then $\xi\mapsto \widetilde G\circ \gamma(\xi) = Q_{k+1}G(x_{0:k},x_k+.)\big(\gamma(\xi)\big)$ is convex in $\xi$.

\smallskip
$\rhd$ Assume all  the functions $\sigma_k$, $k\!=\!0,\ldots, n-1$, are non-negative and convex. One shows by a (second) backward induction that  the functions $\Phi_k$ are all convex. 

 Finally, we prove that $\Phi_k\le \Psi_k$ for every $k=0,\ldots,n-1$, using  again a (third) backward induction on $k$. First note that $\Phi_n=\Psi_n=\Phi$. If $\Phi_{k+1}\le \Psi_{k+1}$, then 
\begin{eqnarray*}
\Phi_k (x_{0:k})=\big(Q_{k+1}\Phi_{k+1}(x_{0:k},x_k+.)\big)(\sigma_k(x_k))&\le& \big(Q_{k+1}\Phi_{k+1}(x_{0:k},x_k+.)\big)(\theta_k(x_k))\\
&\le&  \big(Q_{k+1}\Psi_{k+1}(x_{0:k},x_k+.)\big)(\theta_k(x_k))= \Psi_k(x_{0:k}).
\end{eqnarray*}
In particular, when $k=0$, we get $\Phi_0(x)\le \Psi_0(x)$ or, equivalently, taking advantage of the martingale property, $\E\, \Phi(X_{0:n})\le \E\, \Phi(Y_{0:n})$. 

\smallskip $\rhd$ If all the functions $\theta_k$, $k=0,\ldots,n-1$  are convex, then all  functions $\Psi_k$, $k=0,\ldots,n$, are convex and one shows like wise that $\Phi_k\le \Psi_k$ for every $k=0,\ldots,n-1$.

\smallskip
\noindent $(b)$ The proof follows the same lines as $(a)$ calling upon Claim~$(c)$ of Lemma~\ref{Jen}. In particular, the functions $u\mapsto \widetilde G(x_{0:k},u)$ is also even so that $\sup_{u\in [-a,a]} \widetilde G(x_{0:k},u)= \widetilde G(x_{0:k},a)$ for any $a\ge 0$. 
$\cqfd$

\bigskip
To prove Theorem~\ref{FCOEuro} we need to transfer the above result into a continuous time setting by a functional weak approximation result. To this end, we introduce the notion of {\em piecewise affine interpolator} and recall an elementary weak convergence lemma. 

 \begin{Dfn} $(a)$  For every integer $n\ge 1$,  let $i_n : \R^{n+1}\to {\cal C}([0,T], \R)$ denote the  {\em piecewise affine  interpolator} defined by 
\[
\forall\, x_{0:n}\!\in \R^{n+1},\; \forall\,k=0,\ldots,n-1,\; \forall\, t\!\in [t^n_k,t^n_{k+1}],\quad i_n(x_{0:n})(t) =\frac nT\big((t^n_{k+1}-t)   x_{k}+(t-t^n_k)  x_{k+1}\big).
\]

\noindent $(b)$  For every integer $n\ge 1$,  let     $I_n : {\cal F}([0,T], \R)\to {\cal C}([0,T], \R)$ denote  the  {\em functional interpolator}  defined by 
\[
\forall\, \alpha\!\in {\cal F}([0,T], \R),\quad I_n(\alpha)= i_n\big(\alpha(t^n_0),\ldots, \alpha(t^n_n)\big).
\]
\end{Dfn}

We will use extensively the following obvious fact
\[
 \sup_{t\in [0,T]} |I_n(\alpha)_t|\le \sup_{t\in [0,T]}|\alpha(t)|
\]
in particular for uniform integrability purpose.

 \begin{Lem} \label{Interpol}  Let $X^n$, $n\ge 1$, be a sequence of continuous processes weakly converging towards $X$ for the $\|\,.\,\|_{\rm \sup}$-norm. 
 Then  the sequence of continuously interpolated processes  $\widetilde X^n=I_n(X^n) $  of $X^n$, $n\ge 1$,   is weakly converging toward $X$  for the $\|\,.\,\|_{\rm \sup}$-norm topology.
 \end{Lem}

 \noindent {\bf Proof.} For every integer $n\ge 1$ and every $\alpha\!\in {\cal F}([0,T], \R^d)$, the interpolation operators $I_n(\alpha)$ reads 
 \[
 I_n(\alpha)= \frac nT\big((t^n_{k+1}-t)   \alpha(t^n_k)+(t-t^n_k)   \alpha(t^n_{k+1})\big), \;t\!\in [t^n_k,t^n_{k+1}], \; k=0,\ldots,n-1.
 \]
Note that $I_n$ maps ${\cal C}([0,T], \R^d)$ into itself.  One easily checks that $\|I_n(\alpha)-\alpha\|_{\rm \sup}\le w(\alpha, T/n)$  (keep in mind that  $w$ denotes the uniform continuity modulus of $\alpha$) and $\|I_n(\alpha)-I_n(\beta)\|_{\rm \sup}\le \|\alpha-\beta\|_{\rm \sup}$. We use the standard distance $d_{wk}$ for weak convergence on Polish metric spaces defined by 
 \[
d_{wk}\big({\cal L}(X),{\cal L}(Y)\big) = \sup\big\{|\E\, F(X)-\E\, F(Y)|,\;  [F]_{\rm Lip}\le 1,\; \|F\|_{\rm \sup}\le 1 \big\}.
 \]
 Then 
 \begin{eqnarray*}
d_{wk}\big({\cal L}(I_n(X^n)),{\cal L}(X)\big)&\le&d_{wk}\big({\cal L}(I_n(X^n)),{\cal L}(I_n(X))\big) +  d_{wk}\big({\cal L}(I_n(X)),{\cal L}(X)\big)\\
&\le&  d_{wk}\big({\cal L}(X^n),{\cal L}(X)\big) +\E\,\big(w(X,T/n)\wedge 2\big)
 \end{eqnarray*}
 which goes to $0$ since $X$ has continuous paths.$\cqfd$

\bigskip
\noindent{\bf Proof of Theorem~\ref{FCOEuro}.} 
%
%
We consider now for both SDEs (related to  $X^{(\sigma)}$ and $X^{(\theta)}$) their continuous (also known as  ``genuine") Euler schemes with step $\frac Tn$, starting at $x$ with respect to a given standard Brownian motion $W$ defined on an appropriate probability space. $E.g.$, to be more precise, the Euler scheme related to $X^{(\sigma)}$ is defined by  
 \begin{eqnarray*}
\bar X^{(\sigma),n}_{t^n_{k+1}}&=& \bar X^{(\sigma),n}_{t^n_k} +
\sigma(t^n_k,\bar X^{(\sigma),n}_{t^n_k})\big(W_{t^n_{k+1}}-W_{t^n_k}\big), \; k=0,\ldots,n-1,\; \bar X^{(\sigma),n}_{0}= x\\
\bar X^{(\sigma),n}_{t}&=& \bar X^{(\sigma),n}_{t^n_k}+ 
\sigma(t^n_k,\bar X^{(\sigma),n}_{t^n_k})\big(W_t-W_{t^n_k}\big),\quad t\!\in [t^n_k,t^n_{k+1}).
\end{eqnarray*}
It is clear that both sequences $( \bar X^{(\sigma),n}_{t^n_k} )_{k=0:n}$ and $( \bar X^{(\theta),n}_{t^n_k} )_{k=0:n}$  are of the form~(\ref{modele}) with the Gaussian  white noise sequence $Z_k= W_{t^n_{k}}-W_{t^n_{k-1}}$, $k=1,\ldots,n$. Furthermore, owing to the linear growth assumption made on  $\sigma$ and $\theta$, the sup-norm of these  Euler schemes of Brownian diffusions  lie in $L^p(\P)$ for any $p\!\in (0,+\infty)$, uniformly in $n$,  (see $e.g.$ Lemma B.1.2 p.275 in~\cite{BOLE} or Proposition~\ref{pro:EulerWeakCv} in Appendix~A)
\[
\sup_{n\ge 1}\big\|\sup_{t\in [0,T]}|\bar X^{(\sigma),n}_t|   \big\|_{p}+\sup_{n\ge 1}\big\|\sup_{t\in [0,T]}|\bar X^{(\theta),n}_t|   \big\|_{p}<+\infty.
\]

Furthermore, $I_n(\bar X^{(\sigma),n})= i_n\big((\bar X^{(\sigma),n})_{t^n_{0:n}}\big)$ is but the piecewise affine interpolated Euler scheme (which coincide with $\bar X^{(\sigma),n}$ at times $t^n_k$). Note that the sup-norm of $I_n(\bar X^{(\sigma),n})$ also has finite polynomial moments uniformly in $n$ like the genuine Euler scheme.

\smallskip
Let $F:{\cal C}([0,T], \R)\to \R$ be a convex functional with $(r,\|\,.\,\|_{\sup})$-polynomial growth. For every   integer $n\ge 1$, we define on $\R^{n+1}$ the function $F_n$   by
\begin{equation}\label{eq:Fn}
F_n(x_{0:n})= F\big(i_n(x_{0:n})\big), \; x_{0:n}\!\in \R^{n+1}.
\end{equation}
It is clear that the convexity of  $F$ on ${\cal C}([0,T], \R)$ is transferred to the functions $F_n$, $n\ge 1$. So does  the polynomial  growth property. Moreover, $F$ is $\|\,.\,\|_{\sup}$-continuous since it is convex with $\|\,.\,\|_{\sup}$-polynomial growth (see~Lemma~2.1.1 in~\cite{LUC}). It follows from Proposition~\ref{FCOEuroDiscTime} applied with $\Phi=F_n$, $(Z_k)_{1\le k\le n}= (W_{t^n_k}-W_{t^n_{k-1}})_{1\le k\le n}$, $\sigma_k= \sigma(t^n_k,.)$ and $\theta_k= \theta(t^n_k,.)$, $k=0,\ldots,n$ which obviously satisfy the required linear growth and integrability assumptions,  that, for every $n\ge 1$, 
\begin{equation}\label{eq:Comparsigthet}
\E\, F\big(I_n(\bar X^{(\sigma),n}) \big)= \E\,F_n\big((\bar X^{(\sigma),n}_{t^n_k})_{k=0:n}\big)\le \E\,F_n\big((\bar X^{(\theta),n}_{t^n_k}))_{k=0:n}\big) = \E\, F\big(I_n(\bar X^{(\theta),n})\big).
\end{equation}
On the other hand, it is classical background 
that the {\em genuine (continuous) Euler schemes} $\bar X^{(\sigma),n}$ 
 weakly converges for the $\|\,.\,\|_{\sup}$-norm topology toward $X^{(\sigma)}$, unique weak solution to the  $SDE\equiv dX_t =\sigma(X_t)dW_t$, $X_0=x$, as  $n\to +\infty$. For a proof we refer $e.g.$  to exercise~23 in~\cite{PRO}, p.359 when $\sigma$ is homogeneous in $t$, see also~\cite{JaMePa, KUPR}; we  also provide a short  self-contained proof in Proposition~\ref{pro:EulerWeakCv} in Appendix~\ref{App:A}). The key in all them being the weak convergence theorem for stochastic integrals first established in~\cite{JaMePa}.
%
%

It follows from Lemma~\ref{Interpol} and the $L^p(\P)$-boundedness of the sup-norm of the sequence $(I_n (\bar X^{(\sigma),n}))_{n\ge 1}$ for $p>r$  that
\[
\E \,F(X^{(\sigma)} )= \lim_n \E\, F\big(I_n (\bar X^{(\sigma),n})\big) = \lim_n\E\, F_n\big((\bar X^{(\sigma),n}_{t^n_k})_{0\le k\le n}\big).
\]
The same holds true for the diffusion $X^{(\theta)}$ and its Euler scheme. The conclusion follows.

\smallskip
\noindent $(a)$ Applying successively what precedes  to the couples $(\sigma, \kappa)$ and $(\kappa, \theta)$ until Equation~\eqref{eq:Comparsigthet} respectively, we derive that for every $n\ge 1$, 
\[
\E\, F\big(I_n(\bar X^{(\sigma),n}) \big) \le  \E\, F\big(I_n(\bar X^{(\kappa),n})\big)  \le  \E\, F\big(I_n(\bar X^{(\theta),n})\big)
\]
and one concludes likewise by letting $n$ go to infinity in the resulting inequality
\[
\E\, F\big(I_n(\bar X^{(\sigma),n}) \big) \le     \E\, F\big(I_n(\bar X^{(\theta),n})\big).
\]

\noindent $(b)$  The proof follows the same lines by calling upon item $(c)$ of the above Lemma~\ref{Jen}, having in mind that the  distribution of a standard Brownian increment  is symmetric   with polynomial moments at any order as a Gaussian random vector.~$\cqfd$

\bigskip
\noindent {\bf Remarks.} $\bullet$ Note that no ``weak uniqueness" assumption is requested for the  function $\kappa$.

\smallskip
\noindent   $\bullet$ The Euler scheme has already been successfully used to establish convex order in~\cite{BRusch1}.

\medskip The following corollaries can be obtain with obvious adaptations of the above proof. 

\begin{Cor} Under the above assumption of Claim $(a)$, if, furthermore,  the $SDE$
\[
dX^{(\kappa)}_t = \kappa(t,X^{(\kappa)}_t) dW_t, \; X^{(\kappa)}_0=x
\]
has a unique weak solution, then, for every convex
 functional $F:{\cal C}([0,T], \R)\to \R$  with  $(r,\|\,.\,\|_{\sup})$-polynomial growth,
 \[
\E\,F(X^{(\sigma)})\le \E\,F(X^{(\kappa)})\le  \E\,F(X^{(\theta)}).
\]  
\end{Cor}

\begin{Cor}  If $\sigma,\theta:[0,T]\times I\to \R$, where $I$ is a nontrivial interval of $\R$, are continuous with polynomial growth  and if the related Brownian  $SDE$s satisfy a weak uniqueness assumption for every  $I$-valued  weak solution starting
 from $x\!\in I$, at time $t=0$. Then the above Proposition remains true (the extension of the functional weak convergence of the Euler scheme established in Appendix~\ref{App:A} (Proposition~\ref{pro:EulerWeakCv}) under the assumption made on the drift $b$  is left to the reader).
\end{Cor}

%
%
This approach based on the combination of a (linear) dynamic programming principle and a functional weak approximation argument also allows us to retrieve  Hajek's result  for drifted diffusions.
\begin{Pro}[Extension to drifted diffusions, see~\cite{HAJ}] Let $\sigma$ and $\theta$ be two functions on $[0,T]\times \R$ satisfying     the partitioning or the dominating assumptions $(a)$ or $(b)$ from Theorem~\ref{FCOEuro} respectively. Let $b:[0,T]\times \R\to \R$ be a continuous function with linear growth in $x$ uniformly in $t$ and such that $b(t,.)$ is convex for every $t\!\in [0,T]$. Let $Y^{(\sigma)}$ and $Y^{(\theta)}$  be the   {\em weak} solutions, supposed to be unique, starting from $x$ at time $0$ to the SDEs
$dY^{(\sigma)}_t = b(t,Y^{(\sigma)}_t ) dt +\sigma(t,Y^{(\sigma)}_t )dW^{(\sigma)}_t$ and $dY^{(\theta)}_t = b(t,Y^{(\theta)}_t ) dt +\theta(t,Y^{(\theta)}_t )dW^{(\theta)}_t$.
Then, for every   {\em non-decreasing} convex function $f:\R\to \R$,
\[
\E\, f(X^{(\sigma)})\le \E\, f(X^{(\theta)}).
\]
\end{Pro}

\noindent {\bf Proof.} We have to introduce the operators $Q_{b,\gamma, t}$, $\gamma>0$,  $t\!\in [0,T]$, defined for every Borel function $f: \R\to \R$ (satisfying the appropriate polynomial growth assumption in accordance with the existing moments of $Z$)  by
\[
Q_{b,\gamma,t}(f)(x,u)= \E \,f\big(x+\gamma b(t,x)+uZ\big).
\]
One shows like  in Lemma~\ref{Jen} above that, if the function $f$ is convex, $Q_{b,\gamma,t }f$ is convex in $(x,u)$,  non-decreasing in   $u$ on $\R_+$, non-increasing in $u\!\in \R_-$.  $\cqfd$

\subsection{Applications to (Brownian) functional peacocks} 
We consider a local volatility model on the discounted risky asset dynamics given by 
\begin{equation}\label{LocalVol}
dS^{(\sigma)}_t =S^{(\sigma)}_t\sigma(t,S^{(\sigma)}_t)dW^{(\sigma)}_t, \; S^{(\sigma)}_0=s_0>0,
\end{equation}
where $\sigma:[0,T]\times \R\to \R$ is a bounded continuous function so that the above equation has at least a weak solution $(S^{(\sigma)}_t)_{t\in [0,T]}$ with distribution on a probability space $(\Omega,{\cal A}, \P)$ on which lives  a Brownian motion $(W^{(\sigma)}_t)_{\in [0,T] }$ (with augmented filtration $({\cal F}^{W^{(\sigma)}}_t)_{t\in [0,T]}$). This  follows from the proof of Proposition~\ref{pro:EulerWeakCv} in Appendix~\ref{App:A} (see also~\cite{REYO}, p. ??). Then, $(S^{(\sigma)}_t)_{t\in [0,T]}$ is a true $({\cal F}^{W^{(\sigma)}}_t)_{t\in [0,T]}$-martingale satisfying 
\[
S^{(\sigma)}_t = s_0 \exp{\Big(\int_0^t \sigma(s,S^{(\sigma)}_s) dW^{(\sigma)}_s-\frac 12\int_0^t \sigma^2(s,S^{(\sigma)}_s)ds\Big)}
\]
so that  $S^{(\sigma)}_t >0$ for every  $t\!\in [0,T]$. One  introduces likewise the local volatility model $(S^{(\theta)}_t)_{t\in [0,T]}$ related to the bounded volatility function $\theta:[0,T]\times\R_+\to \R$, still starting from $s_0>0$. Then, the following proposition holds which appears as a functional or non-parametric  extension of the fact that $\Big(\int_0^Te^{\sigma B_t-\frac{\sigma^2t}{2}} dt\Big)_{\sigma\ge 0}$ is a peacock (see~$e.g.$~\cite{Carretal, YORetal}).

\begin{Pro}[Functional peacocks] \label{pro:FuncPeac}Let $\sigma$ and $\theta$ be two real valued bounded continuous functions defined on $[0,T]\times \R$.  Assume that $S^{(\sigma)}$ is the unique weak solution to~\eqref{LocalVol} as well as   $S^{(\theta)}$ for its {\em mutatis mutandis} counterpart  involving $\theta$. If one of the following additional conditions holds
\begin{eqnarray*}
(i)\!&\hskip -0.75 cm \mbox{{\em Partitioning function}: there exists a function } \kappa:[0,T]\times \R_+ \to \R_+ \mbox{ such that, for every $t\!\in [0,T]$, }\hskip 8 cm \\ 
\!&\!\! x\mapsto x\,\kappa(t,x)    \mbox{ is convex  on $\R_+$}  \mbox{ and }\; 0\le \sigma(t,.)\le \kappa(t,.)\le \theta(t,.) \mbox{ on $\R_+$,} \hskip 13 cm  & \\
\mbox{or }&&\\
(ii)
\!&\!\!  \mbox{\em Domination property: for every $t\!\in [0,T]$ the function } x\mapsto x\,\theta(t,x)   \mbox{ is convex on $\R_+$} \mbox{ and }\hskip 10 cm &\\
\!&\!\!\ \;  |\sigma(t,.)|\le \theta(t,.),\hskip 6 cm\hskip 6 cm&
\end{eqnarray*}
then, for every  convex (hence continuous) function $f: \R\to \R$ with  polynomial growth 
\[
\E\, f\left(\int_0^T S^{(\sigma)}_s\mu(ds)\right)\le \E\, f\left(\int_0^T S^{(\theta)}_s\mu(ds)\right)
\]
where $\mu$ is a signed (finite) measure on $([0,T], {\cal B}or([0,T]))$. More generally, for every  
convex functional $F: {\cal C}([0,T],\R_+)\to \R$ with  $(r,\|\,.\,\|_{\sup})$-polynomial growth  polynomial growth,  one has,
\begin{equation}\label{ComparLocVol}
  \E\, F\big( S^{(\sigma)}\big)\le \E\, F\big( S^{(\theta)}\big).
\end{equation}
\end{Pro}

\noindent {\bf Proof.} We focus  on the first {\em partitioning} setting. The second one can be treated likewise. First note that $\kappa$ is bounded since $\theta$ is. As a consequence, the function $x\mapsto x\,\kappa(t,x)$ is zero at $x=0$ and can be extended into  a convex function on the whole real  line  by setting $x\,\kappa(t,x)= 0$ if $x\le 0$. One  extends $x\,\sigma(t,x)$ and $x\,\theta(t,x)$ by zero on $\R_-$ likewise. Once this has been done, this claim appears as a straightforward consequence of Theorem~\ref{FCOEuro} for the (martingale) diffusion processes  whose  diffusion coefficients are given by   (the extension) of $x\,\sigma(t,x)$ and $x\,\theta(t,x)$ on the whole real line. As above, the sup-norm continuity follows from the convexity and polynomial growth.
In the end we take advantage of the {\em a posteriori} positivity of $S^{(\theta)}$ and $S^{(\sigma)}$  when starting from  $s_0>0$ to conclude.~$\cqfd$

\bigskip
%
%
\noindent{\sc Applications to volatility comparison results.} The corollary below shows that comparison results  for vanilla European options  established  in~\cite{NEKJESH} appear as  a special case of Proposition~\ref{pro:FuncPeac}.

\begin{Cor} Let $\sigma:[0,T]\times \R \to \R_+$ be a bounded continuous function 
\[
0\le \sigma_{\min}(t)\le \sigma(t,.)\le  \sigma_{\max}(t), \; t\!\in [0,T],
\]
then for every  convex functional $F: {\cal C}([0,T], \R_+)\to \R$ with $(r,\|\,.\,\|_{\sup})$-polynomial growth ($r\ge 1$),
\begin{equation}\label{ComparLocVol2}
\E\, F\left( S^{(\sigma_{min})}_s\right)\le \E\, F\left( S^{(\sigma)}_s\right)\le \E\, F\left( S^{(\sigma_{\max})}_s\right).
\end{equation}
\end{Cor}

\noindent {\bf Proof.} We successively apply the former Proposition~\ref{pro:FuncPeac} to the couple $(\sigma_{\min}, \sigma)$ and the partitioning function $\kappa(t,x)= \sigma_{\min}(t)$ to get the left inequality and to the couple  $(\sigma, \sigma_{\max})$ with $\kappa= \sigma_{\max}$ to get the right inequality. $\cqfd$

\medskip
Note that  the left and right hand side of the above inequality are usually considered as quasi-closed forms since they correspond to Hull-White model (or even to  the regular Black-Scholes model  if $\sigma_{\min}$, $\sigma_{\max}$ are constant). Moreover, it has to be emphasized that no convexity assumption on $\sigma$ is requested.

\subsection{Counter-example (discrete time setting)} The above comparison results for the convex order can fail when the assumptions of  Theorem~\ref{FCOEuro} are not satisfied by the diffusion coefficient. In fact, for simplicity, the counter-example below is developed in a  discrete time framework corresponding to Proposition~\ref{FCOEuroDiscTime}. 
We consider the $2$-period dynamics $X=X^{\sigma,x}= (X^{\sigma,x}_{0:2})$ satisfying 
\[
X_1= x+ \sigma Z_1\quad\mbox{ and }\quad
X_2 = X_1 + \sqrt{2v(X_1)}Z_2
\]
where $Z_{1:2}\stackrel{\cal L}{\sim} {\cal N}(0;I_2)$, $\sigma\ge 0$,  and $v: \R\to \R_+$ is a bounded ${\cal C}^2$-function such that  $v$ has a strict local maximum at $x_0$ satisfying  $v'(x_0)\!=\!0$ and $v''(x_0)<-1$ (so is the case if $v(x)\! =\! v(x_0) -\rho (x-x_0)^2 + o((x-x_0)^2)$, $0<\rho<\frac 12$, in the neighbourhood of $x_0$). Of course this implies that $\sqrt{v} $ cannot be convex.

Let $f(x) =e^x$. It is clear that 
\[
\varphi(x,\sigma) := \E f(X_2)= e^x \E \big(e^{\sigma Z_1+v(x+\sigma Z_1)}\big).
\]
Elementary computations show that 
\begin{eqnarray*}
\varphi'_{\sigma}(x,\sigma)&=& e^x \E\Big(e^{\sigma Z_1+v(x+\sigma Z_1)}\big(1+v'(x+\sigma Z_1)\big)Z_1\Big)\\
\varphi''_{\sigma^2}(x,\sigma)&=&e^x\left(\E\Big(e^{\sigma Z_1+v(x+\sigma Z_1)}\big(1+v'(x+\sigma Z_1)\big)^2Z^2_1\big)+ \E\Big(e^{\sigma Z_1+v(x+\sigma Z_1)}v''(x+\sigma Z_1)Z^2_1\Big)\right).
\end{eqnarray*}
In particular
\[
\varphi'_{\sigma}(x,0)= e^{x+v(x)}(1+v'(x)) \E\, Z_1 = 0\quad\mbox{ and }\quad
\varphi''_{\sigma^2}(x,0)= e^{x+v(x)}\Big((1+v'(x))^2+v''(x)\Big)
\]
so that $\varphi''_{\sigma^2}(x_0,0)<0$ which implies that there exists a small enough  $\sigma_0>0$ such that 
$$
\varphi'_{\sigma}(x_0,\sigma)<0\quad \mbox{ for every  }\quad\sigma\!\in (0, \sigma_0],
$$
This clearly exhibits a counter-example to Proposition~\ref{FCOEuroDiscTime}
when the convexity assumption is   fulfilled neither by   the functions $(\sigma_k)_{k=0:n}$ nor the functions $(\kappa_k)_{k=0:n}$ (here with $n=1$).

\subsection{L\'evy driven diffusions} \label{Levy}
 Let $Z=(Z_t)_{t\in [0,T]}$ be a L\'evy process with L\'evy measure $\nu $ satisfying $\displaystyle \int_{|z|\ge 1} |z|^p\nu(dz) < +\infty$,  $p\!\in [1,+\infty)$. Then $Z_t\!\in L^1(\P)$ for every $t\!\in [0,T]$. Assume furthermore that $\E\,Z_1=0$:  then    $(Z_t)_{t\in [0,T]}$ is an  ${\cal F}^Z
 $-martingale.
 

 \begin{Thm}\label{FCOEuroJump} Let $Z=(Z_t)_{t\in [0,T]}$ be a martingale L\'evy process with L\'evy measure $\nu $ satisfying 
 $\nu(|z|^p)<+\infty$ for a $p\!\in (1,+\infty)$ if $Z$ has no Brownian component and $\nu(z^2)<+\infty$ if $Z$ does have a Brownian component. Let $\kappa_i\!:\![0,T]\times \R\to \R$, $i\!=\!1,2$, be  continuous functions with linear growth in $x$ uniformly in $t\!\in [0,T]$. For $i=1,2$, let $X^{(\kappa_i)}\!=\!(X^{(\kappa_i)}_t)_{t\in [0,T]}$ be the weak   solution, assumed to be unique,  to 
 \begin{equation}\label{LevyEDS}
 dX^{(\kappa_i)}_t = \kappa_i(t,X^{(\kappa_i)}_{t-})dZ^{(\kappa_i)}_t, \quad X^{(\kappa_i)}_0=x\!\in \R, 
 \end{equation}
where $Z^{(\kappa_i)}$, $i=1,2$ have the same distribution as $Z$. Let $F :I\!\!D([0,T],\R)\to\R$ be a Borel convex functional, $\P_{\!X^{(\kappa_i)}}$-$a.s.$ $Sk$-continuous, $i\!=\!1,2$, with $(r,\|.\|_{\sup})$-polynomial growth for some $r\!\in [1,p)$~$i.e.$ 
\[
\forall\, \alpha\!\in I\!\!D([0,T],\R),\quad |F(\alpha)|\le C(1+\|\alpha\|^r_{\sup}).
\]

\noindent $(a)$  {\em Partitioning function}: If there exists a function $\kappa:[0,T]\times \R\to \R_+$ such that $\kappa(t,.)$ is convex for every $t\!\in [0,T]$ and $0\le \kappa_1\le \kappa\le \kappa_2$,   then 
\[
\E\,F(X^{(\kappa_1)})\le \E\,F(X^{(\kappa_2)}).
\]

\noindent $(a')$ An equivalent form for claim $(a)$ is:  if $0\le \kappa_1\le \kappa_2$ and, either $\kappa_1(t,.)$ is convex for every $t \!\in [0,T]$, or  $\kappa_2(t,.)$ is convex for every $t \!\in [0,T]$, then the conclusion of $(a)$ still holds true.

\smallskip
\noindent $(b)$ {\em Domination property}: If $Z$ has a symmetric distribution, $|\kappa_1| \le \kappa_2$ and  $\kappa_2$ is  convex, then 
\[
\E\,F(X^{(\kappa_1)}) \le \E\,F(X^{(\kappa_2)}).
\]
\end{Thm}

\noindent {\bf Remarks.} $\bullet$ The $\P_{X^{(\kappa_i)}}$-$a.s.$ $Sk$-continuity of the functional $F$, $i=1,2$, is now requested since$Sk$-continuity no longer follows form the convexity ($\big(I\!\!D([0,T],\R), Sk\big)$ is a Polish space but not even a topological vector space). Thus the function  $\alpha \mapsto |\alpha(t_0)|$ for a fixed $t_0\!\in (0,T)$ is  continuous at a given $\beta \!\in I\!\!D([0,T], \R)$ if and only if $\beta $ is continuous at $t_0$ (see~\cite{BIL}, Chapter 3).

\smallskip
\noindent $\bullet$ The result remains true under the less stringent moment assumption on the L\'evy measure $\nu$: $\nu(|z|^p\mbox{\bf 1}_{\{|z|\ge 1\}}<+\infty$ but would require much more technicalities since one has to carry out the reasoning of the proof below between two large jumps of $Z$ and ``paste" these inter-jump results.  

\smallskip
The following technical lemma is the key that solves the approximation part of the proof in this c\`adl\`ag setting. 
\begin{Lem}\label{Lem:6.3}\label{lem:Jacod} Let $\alpha\!\in I\!\!D([0,T], \R)$. The sequence of stepwise constant  approximations defined by
\[
\alpha_n(t)= \alpha(\underline t_n),\; t\!\in [0,T],
\]
converges toward $\alpha$ for the Skorokhod topology.
\end{Lem}


 \noindent {\bf Proof.} See~\cite{JASH2}Proposition~ VI.6.37~p.387 (second edition).~$\cqfd$

 \bigskip
 \noindent {\bf Proof of Theorem~\ref{FCOEuroJump}.}    {\sc Step~1.} Let $(\bar X^n_t)_{t\in [0,T]}$ be  the genuine  Euler scheme defined by
 \[
 \bar X^n_t= x + \int_{(0,t]}\  \kappa(\underline s_n, \bar X^n_{{\underline s_n}_-})dZ_s
 \]
where $  \kappa= \kappa_1$ or $\kappa_2$. Then, owing to the linear growth of $  \kappa$, we derive (see$e.g.$~Proposition~\ref{pro:highmoment}  in Appendix~\ref{app:B}) that 
\[
\Big\| \sup_{t\in[0,T]}|  X_t |  \Big \|_p+ \sup_{n\ge 1}\Big \|  \sup_{t\in [0,T]}|\bar X^n_t |   \Big\|_p<+\infty.
\]
 We know, $e.g.$ from  form Proposition~\ref{pro:WeakEulerLevy} in Appendix~\ref{app:B},  that $(\bar X^n)_{n\ge 1}$ functionally weakly converges for the Skorokhod topology toward the unique weak solution $X$ of the   $SDE$ $dX_k=   \kappa(t,X_{t_-})dZ_t$, $X_0=x$.  In turn, Lemma~\ref{Lem:6.3} implies that $(\bar X^n_{\underline t_n})_{t\in [0,T] }$  $Sk$-weakly converges  toward   $X$.

\medskip
 \ni  {\sc Step~2.} Let  $F:\D([0,T], \R)\to \R$ be a $\P_{_X}$-$Sk$-continuous convex functional. For every integer $n\ge 1$, we still define the sequence of convex functional $F_n :\R^{n+1}\to \R$  by \\ $\displaystyle F_n(x_{0:n}) = F \Big( \sum_{k=0}^{n-1}x_k\mbox{\bf 1}_{[t^n_k,t^n_{k+1})}+ x_n\mbox{\bf 1}_{\{T\}}\Big)$ so that $F_n\big((\bar X^n_{t^n_k})_{0:n}\big) =F\big((\bar X^n_{\underline t_n})_{t\in [0,T]}\big)$.

\smallskip Now, for every $n\ge 1$,  the discrete time Euler schemes $\bar X^{(\kappa_i),n}$, $i\!=\!1,2$, related to the jump diffusions with diffusion coefficients $\kappa_1$ and $\kappa_2$ are of the form~(\ref{modele}) and $|F_n(x_{0:n})|\le C(1+\|x_{0:n}\|^r)$, $r\!\in [1,p)$. 

\medskip \ni $(a)$ Assume $0\le \kappa_1\le \kappa_2$. Then, taking advantage of the partitioning function $\kappa$, it follows from Proposition~\ref{FCOEuroDiscTime}$(a)$ that, for every $n\ge 1$,  $\E\, F_n\big((\bar X^{(\kappa_1),n}_{t^n_k})_{0:n}\big)\le \E\,   F_n\big((\bar X^{(\kappa_2),n}_{t^n_k})_{0:n}\big)$ $i.e.$ $\E\, F\big((\bar X^{(\kappa_1),n}_{\underline t_n})_{t\in[0,T]}\big) \le \E\,   F\big((\bar X^{(\kappa_2),n}_{\underline t_n})_{t\in[0,T]}\big)$. Letting $n\to+ \infty$ completes the proof like in   that of  Theorem~\ref{FCOEuro} since $F$ is $\P_{_X}$-$a.s.$ $Sk$-continuous. $\cqfd$

\ms
\ni $(b)$ is an easy consequence of Proposition~\ref{FCOEuroDiscTime}$(b)$.~$\cqfd$


%

\section{Convex order   for non-Markovian It\^o and Dol\'eans martingales}\label{sec:ItoProc}
The results of this section illustrates  another aspects of our paradigm in order to establish functional convex order 
for various classes of continuous time  stochastic processes.
Here we deal with (couples of)  It\^o-intregrals 
with the restriction that   one of the two integrands needs to be  deterministic. 

\subsection{ It\^o martingales}
\begin{Pro}\label{NonMarkovIto} Let $(H_t)_{t\in [0,T]}$ be an $({\cal F}_t)$-progressively measurable process defined on a )  filtered probability space $(\Omega,{\cal A}, ({\cal F}_t)_{t\in [0,T]}, \P)$ satisfying the usual conditions and let $h=(h_t)_{t\in [0,T]}\!\in L^2_{_T}$. Let $F : {\cal C}([0,T],\R)\to\R$ be a convex 
functional with $(r,\|.\|_{\sup})$-polynomial growth,~$r\!\ge\!1$.
%

\medskip
\noindent $(a)$ If $ |H_t| \le   h_t$ $\P$-$a.s.$ for every $t\in [0,T]$, then 
\[
\E \,F\left(\int_0^.H_sdW_s\right)\le \E\, F\left(\int_0^.h_sdW_s\right).
\]

\noindent $(b)$ If  $H_t \ge h_t \ge 0$ $\P$-$a.s.$ for every $t\in [0,T]$ and $\|H\|_{L^2_{_T}}\!\in L^{r'}(\P)$ for $r'>r$, then 
\[
\E\, F\left(\int_0^.H_sdW_s\right)\ge \E\,F \left(\int_0^. h_sdW_s\right).
\]
\end{Pro}

\noindent {\bf Remarks.} $\bullet$
%
In the ``marginal" case where $F$ is of the from $F(\alpha)= f(\alpha(T))$,  it has been shown in~\cite{HIYO2} that the above assumptions on $H$ and $h$  in $(a)$ and $(b)$ are too stringent and can be relaxed into 
\[
 \int_0^T \E\, H_t^2 dt  \le   \int_0^T h^2_tdt \quad \mbox{ and }\quad  \int_0^T \E\, H_t^2 dt  \ge   \int_0^T h^2_tdt 
\]
respectively. The main ingredient of the proof is the Dambis-Dubins-Schwartz representation theorem for one-dimensional Brownian martingales (see $e.g.$~Theorem 1.6 in~\cite{REYO}, p.181).

\smallskip
\noindent$\bullet$ The first step of the proof below, can be compared to Proposition~\ref{FCOEuroDiscTime} in a Markov framework as an autonomous proposition devoted to discrete time setting.


\bigskip
\noindent {\bf Proof.} {\sc Step~1} (Discrete time). Let $(Z_k)_{1\le k\le n}$ be an $n$-tuple of independent symmetric (hence centered)  $\R$-valued random variables satisfying $Z_k\!\in L^r(\Omega, {\cal A}, \P)$, $r\ge 1$,  and let  ${\cal F}^Z_0=\{\emptyset, \Omega\}$, ${\cal F}^Z_k= \sigma\big(Z_, \ldots,Z_k\big)$, $k=1,\ldots,n$ be its natural filtration. Let $(H_k)_{0\le k\le n}$ be an $({\cal F}^Z_k)_{0\le k\le n}$-adapted sequence such that $H_k\!\in L^r(\P)$,  $k=1,\ldots,n$.

Let $X\!=\!(X_k)_{0\le k\le n}$ and $Y\!=\!(Y_k)_{0\le k\le n}$ be the two sequences of random variables recursively defined by
\[
X_{k+1}= X_k +H_k Z_{k+1},\quad Y_{k+1}= Y_k +   h_k Z_{k+1}, \quad 0\le k\le n-1,\quad X_0=Y_0=x_0.
\]
These are the discrete time stochastic integrals of $(H_h)$ and $(h_k)$ with respect to $(Z_k)_{1\le k\le n}$.  It is clear by induction that $X_k$, $Y_k\!\in L^r(\P)$ for every $k=0,\ldots,n$ since $H_k$ is ${\cal F}^Z_k$-measurable and $Z_{k+1}$ is independent of ${\cal F}^Z_k$. 

Let $\Phi:\R^{n+1}\to \R$ be a convex function such that $|\Phi(x)|\le C(1+|x|^{r})$ where $C\ge 0$ is a real constant. 
 Let us focus on the first inequality (discrete time counterpart of claim $(a)$).  One proceeds like in the proof Proposition~\ref{FCOEuroDiscTime} to prove by (three) backward induction(s) that if $|H_k|\le   h_k$, $k=0\!:\!n$, then 
\[
\E\, \Phi(X)\le \E\, \Phi(Y).
\]
To be more precise, let us introduce by analogy with this proposition the sequence $(\Psi_k)_{0\le k\le n}$ of  functions recursively  defined by 
\[
\Psi_n =\Phi,\; \Psi_k(x_{0:k})= \left(Q_{k+1}\Psi_{k+1}(x_{0:k}, x_k+.)\right)\!(h_k),\, x_{0:k}\!\in \R^{k+1}, \,k=0,\ldots,n-1.
\]
First note that the functions $\Psi_k$  satisfy a linear dynamic programing principle
\[
\Psi_k(Y_{0:k})= \E\big(\Psi_{k+1}(Y_{0:k+1})\,|\, {\cal F}^Z_k\big),\; k=0,\ldots,n-1
\]
so that by the chaining rule for conditional expectations, we have
$$
\Phi_k(Y_{0:k})= \E\big(\Phi(Y_{0:n})\,|\, {\cal F}^Z_k\big),\; k=0,\ldots,n.
$$
 Furthermore,  owing to the properties of the operator $Q_{k+1}$,  we already proved that for any convex function $G:\R^{k+2}\to \R$ such that $|G(x)|\le C(1+|x|^r)$, the function 
$$
(x_{0:k},u)\mapsto (Q_{k+1}G(x_{0:k},x_k+.))(u)= \E \,G(x_{0:k},x_k+uZ_{k+1})
$$
 is convex and even as a function of $u$ for every fixed $x_{0:k}$. As a consequence, it also satisfies the  maximum principle established in 
 Lemma~\ref{Lem:Jensen}$(c)$ since the random variables $Z_k$ have symmetric distributions. 
 
Now, let us  introduce the martingale induced by $\Phi(X_{0,n})$, namely 
\[
M_k = \E\big( \Phi(X_{0:n})\,|\, {\cal F}^Z_k)\big),\; k\!\in \{0,\ldots,n\}.
\]
We will show by a backward induction that $M_k\le \Psi_k(X_{0:k})$ for every $k\!\in \{0,\ldots,n\}$. If $k=n$, this is trivial. Assume now
that $M_{k+1}\le \Psi_{k+1}(X_{0:k+1})$ for a $k\!\in \{0,\ldots,n-1\}$. Then we get the following string of inequalities 
\begin{eqnarray}
\nonumber M_k = \E(M_{k+1}\,|\,{\cal F}^Z_k)&\le& \E(\Psi_{k+1}(X_{0:k+1})\,|\,{\cal F}^Z_k)\\
\nonumber&=& \E(\Psi_{k+1}(X_{0:k}, X_k+H_kZ_{k+1})\,|\,{\cal F}^Z_k)\\
\nonumber &=& \Big(\E(\Psi_{k+1}(x_{0:k}, x_k+uZ_{k+1})\,|\,{\cal F}^Z_k)\Big)_{| x_{0;k}= X_{0:k}, u=H_k}\\
\nonumber&=& 
\Big(Q_{k+1}\Psi_{k+1}(x_{0:k}, x_k+.)(H_k)\Big)_{| x_{0;k}= X_{0:k}}\\
\label{eq:H-le-h}&\le& \Big(Q_{k+1}\Psi_{k+1)}(x_{0:k}, x_k+.)(h_k)\Big)_{| x_{0;k}= X_{0:k}}= \Psi_k( X_{0:k})
\end{eqnarray}
where we used in the fourth line that $Z_{k+1}$ is independent of ${\cal F}^Z_k$ and in the penultimate line  the assumption $|H_k|\le h_k$ and the  maximum  principle. Finally,  at $k=0$, we get $\E \, \Phi(X_{0,n})= M_0 \le \Phi_0(x_0)= \E \Phi(Y_{0:n})$ which is  the announced conclusion. 

%
%
%
%
%

\smallskip 
\noindent {\sc Step~2} (Approximation-Regularization). We temporarily assume that the function $h$ (has a modification which) is bounded by a real constant so that $\P(d\omega)$-$a.s.$ $\|H (\omega)\|_{\sup} \vee \|h\|_{\sup} \le K$. We first need a technical lemma inspired by  Lemma~2.4 in~\cite{KASH} (p.132, $2^{nd}$ edition) about approximation of progressively measurable processes by {\em simple} processes, with in mind
 the preservation of the domination property  requested in our framework. 

\begin{Lem} \label{lem:techapprox}$(a)$ For every $\varepsilon \!\in (0,T)$ and every $g\!\in L^2([0,T], dt)$ we define 
\[
\Delta_{\varepsilon}g(t) \equiv t\longmapsto \frac{1}{\varepsilon}\int_{(t-\varepsilon)_+}^t g(s)ds\!\in {\cal C}([0,T],\R).
\]
The operator $\Delta_{\varepsilon}:L^2_{_T}\to {\cal C}([0,T], \R)$ is non-negative.   In particular,  if $g,\, \gamma\!\in L^2_{_T}$ with $|g|\le \gamma$ $\lambda_1$-$a.e.$, then $|\Delta_{\varepsilon} g| \le \Delta_{\varepsilon} \gamma$ and $\|\Delta_{\varepsilon} g\|_{\sup} \le |g|_{L^{\infty}_{_T}}$. 

\medskip
\noindent $(b)$ If $g\!\in {\cal C}([0,T],\R)$, define for every integer $m\ge 1$, the stepwise constant c\`agl\`ad (for left continuous right limited) approximation $\tilde g_n$ of $g$ by  
\[
\widetilde g^m(t)=g(0) \mbox{\bf 1}_{\{0\}}(t)+\sum_{k=1}^m g\big(t^m_{k-1}\big) \mbox{\bf 1}_{(t^m_{k-1},t^m_k]}.
\]
Then $\displaystyle \widetilde g^m \stackrel{\|\,.\,\|_{\sup}}{\longrightarrow} g$ as $m\to +\infty$. Furthermore, if $g$, $\gamma\!\in {\cal C}([0,T], \R)$ and $|g|\le \gamma$, then for every $m\ge 1$, $|\widetilde g^m| \le \widetilde \gamma^m$.
\end{Lem} 

The details  of the proof are left to the reader. 

\smallskip By the Lebesgue  fundamental  Theorem of Calculus we know that 
\[
\big|\Delta_{\frac 1n} H- H\big|_{L^2_{_T}} \longrightarrow 0\quad \P\mbox{-}a.s.
\]
Since $| \Delta_{\frac 1n} H -H|_{L^2_{_T}}\le 2K$, 
the Lebesgue dominated convergence theorem implies that 
\begin{equation}\label{eq:DeltaH1n}
\E \int_0^T | \Delta_{\frac 1n} H_t-H_t|^2dt \longrightarrow 0 \; \mbox{ as }\; n\to +\infty.
\end{equation}
By construction, $\Delta_{\frac 1n}H$ is an $({\cal F}_t)_t$-adapted pathwise continuous process satisfying the domination property $|\Delta_{\frac 1n}H| \le \Delta_{\frac 1n}h$ so that, in turn, using this time claim  $(b)$ of the above lemma, for every $n,\,m\ge 1$, 
\[
|\widetilde{\Delta_{\frac 1n}H_t}^m| \le \widetilde{\Delta_{\frac 1n}h_t}^m.
\]
On the other hand, for every $n\ge 1$, the $a.s.$ uniform continuity of $\Delta_{\frac 1n}H $ over $[0,T]$ implies 
$$
 \int_0^T \big|\widetilde{\Delta_{\frac 1n}H_t}^m- \Delta_{\frac 1n}H_t\big|^2dt\le \sup_{t\in [0,T]}|\widetilde{\Delta_{\frac 1n}H_t}^m -\Delta_{\frac 1n}H_t|^2\to 0\; \mbox{ as }\;m\to+\infty\;\P\mbox{-}a.s.
$$
One concludes again by the Lebesgue dominated convergence theorem that, for every $n\ge1$,
\[
\E \int_0^T \Big| \widetilde{\Delta_{\frac 1n}H_t}^m- \Delta_{\frac 1n}H_t \Big|^2dt\longrightarrow 0\;\mbox{ as }\; m\to +\infty.
\]
One shows likewise for the function $h$  it self that
\[
\big| \Delta_{\frac 1n} h - h\big|_{L^2_{_T}}\rightarrow 0 \; \mbox{ as }\; n\to +\infty
\]
and, for every $n\ge 1$,
\[
\big|\widetilde{\Delta_{\frac 1n}h}^m- \Delta_{\frac 1n}h\big|_{L^2_{_T}}\to 0\; \mbox{ as }\;m\to+\infty.
\]
Consequently there exists an increasing subsequence $m(n)\uparrow +\infty$ such that 
\[
\E \int_0^T \Big| \widetilde{\Delta_{\frac 1n}H_t}^{m(n)}- \Delta_{\frac 1n}H_t \Big|^2dt+ \int_0^T \Big| \widetilde{\Delta_{\frac 1n}h_t}^{m(n)}- \Delta_{\frac 1n}h_t \Big|^2dt\longrightarrow 0\;\mbox{ as }\; n\to +\infty
\]
which in turn implies, combined with~\eqref{eq:DeltaH1n} (and its deterministic counterpart for $h$),  
\[
\E \int_0^T \Big| \widetilde{\Delta_{\frac 1n}H_t}^{m(n)}- H_t \Big|^2dt+ \int_0^T \Big| \widetilde{\Delta_{\frac 1n}h_t}^{m(n)}- h_t \Big|^2dt\longrightarrow 0\;\mbox{ as }\; n\to +\infty.
\]

At this stage, we set  for every integer $n\ge 1$,
\begin{equation}~\label{eq:Hnhn}
H_t^{(n)} = \widetilde{\Delta_{\frac 1n} H_t}^{m(n)}\quad\mbox{ and}\quad h^{(n)}_t = \widetilde{\Delta_{\frac 1n} h_t}^{m(n)}
\end{equation}
which satisfy
\begin{equation}\label{eq:approxHh}
\E |H-H^{(n)} |^2_{L^2_{_T}} +|h-h^{(n)}|_{L^2_{_T}} \longrightarrow 0\;\mbox{ as }\; n\to +\infty.
\end{equation}
It should be noted that  these processes $H^{(n)}$, $H$ and these functions  $h^{(n)}$,  $h$ are all bounded by $2K$.
%

\medskip
We consider now the continuous modifications of the four (square integrable) Brownian martingales associated to the integrands  $H^{(n)}$, $H$, $h^{(n)}$ and $h$ (the last two being of Wiener type in fact). It is clear by Doob's Inequality that
\[
\sup_{t\in [0,T]}\Big|\int_0^t H^{(n)}_sdW_s - \int_0^t H_sdW_s \Big| + \sup_{t\in [0,T]}\Big|\int_0^t h^{(n)}_sdW_s - \int_0^t h_sdW_s \Big| \stackrel{L^2(\P)}{\longrightarrow}  0\;   \mbox{ as }\; n\to +\infty.\]

In particular $\displaystyle \Big(\int_0^{.} H^{(n)}_{s}dW_s\Big)_{t\in [0,T]}$ functionally weakly converges  to $ \displaystyle \Big(\int_0^. H_{s}dW_s\Big)_{t\in [0,T]}$ for the $\|\,.\,\|_{\rm sup}$-norm topology. We also have, owing to the $B.D.G.$  Inequality, that for every $p\!\in (0,+\infty)$, 
\begin{equation}\label{eq:BDG}
\E \sup_{t\in [0,T]}\Big|\int_0^t H^{(n)}_sdW_s\Big|^p \le  c^p_p\, \E \|H^{(n)}\|^p_{L^2_{_T}}\le c_p K^p
\end{equation}
where $c_p$ is the universal constant involved in the $B.D.G.$ inequality. The same holds true for the three other integrals related to $h^{(n)}$, $H$,  and $h$. 

Let $n\ge 1$. Set $H^{n}_k= H^{(n)}_{t^{m(n)}_k},\; h^{n}_k= h^{(n)}_{t^{m(n)}_k},\;  k=0,\ldots,m(n)$ and $Z^{n}_k = W_{t^{m(n)}_k}-W_{t^{m(n)}_{k-1}}
$, $k=1,\ldots,n(m)$. One easily checks that $\displaystyle \int_{0}^{t^{m(n)}_k} H^{(n)}_sdW_s= \sum_{\ell=1}^kH^n_{\ell}Z^n_{\ell}$, $k=0,\ldots,m(n)$ so that 
\[
I_{m(n)}\left(\int_0^. H^{(n)}_sdW_s\right) = i_{m(n)} \left(\Big(\sum_{\ell=1}^kH^n_{\ell}Z^n_{\ell} \Big)_{k=0:m(n)}\right).
\]
Let $F_{m(n)}$ be defined by~(\ref{eq:Fn}) from the  convex functional $F$ (with $(r,\|\,.\,\|_{\sup})$-polynomial growth). It  is  clearly  convex. One derives from Step~1 with applied with horizon $m(n)$ and discrete time random sequences  $(Z^n_k)_{k=1:m(n)}$, $(H^n_k)_{k=0:m(n)}$, $(h_k)_{k=0:m(n)}$  that
\begin{eqnarray*}
\E\, F\!\circ\! I_{m(n)}\Big(\int_0^.H^{(n)}_{s}dW_s\Big)&=& \E\, F_{m(n)} \left(\Big(\sum_{\ell=1}^kH^n_{\ell}Z^n_{\ell} \Big)_{k=0:m(n)}\right)\\
&\le & \E\, F_{m(n)} \left(\Big(\sum_{\ell=1}^kh^n_{\ell}Z^n_{\ell} \Big)_{k=0:m(n)}\right) = \E\, F\!\circ \!I_{m(n)}\Big(\int_0^.h^{(n)}_{s}dW_s\Big).
\end{eqnarray*}
 Combining the above functional weak convergence, Lemma~\ref{Interpol} and the uniform integrability   derived form~\eqref{eq:BDG} (with any $p>r$)  yields the expected inequality by letting $n$ go to infinity.
 
 \smallskip 
\noindent {\sc Step~3.} (Second approximation) Let $K\!\in \N$ and $\chi_K:\R\to \R$ the thresholding function defined by $\chi_K(u)= (u\wedge K)\vee (-K)$. It follows from the $B.D.G.$ Inequality that for every $p\!\in (0,+\infty)$
\begin{eqnarray}
\nonumber \E \sup_{t\in [0,T]}\Big|\int_0^t H_sdW_s-\int_0^t \chi_K(H_s)dW_s\Big|^p &\le & c^p_p \,\E |H-\chi_K(H )|^p_{L^2_{_T}}\\
\label{eq:BGG21}&=&c^p_p \, \E |\big(|H|-K\big)_+|^p_{L^2_{_T}} \\
\label{eq:BDG22}&\le&c^p_p \, |\big(|h|-K\big)_+|^p_{L^2_{_T}}
\end{eqnarray}
where $u_+=\max(u,0)$, $u\!\in \R_+$. The same bound obviously holds  when replacing $H$ by $h$. This shows that the convergence holds in every $L^p(\P)$ space, $p\!\in (0,+\infty)$ as $K\to +\infty$. Hence, one may let $K$ go to infinity in the inequality
\begin{equation}\label{eq:ineqchi-K}
\E F\left(\int_0^. \chi_K(H_s) dW_s \right)\le \E F\left(\int_0^. \chi_K(h_s) dW_s \right) = \E F\left(\int_0^. h_s\wedge K dW_s \right)  
\end{equation}
which yield the expected inequality.

\smallskip
\noindent $(b)$  We consider the same steps as for the upper-bound established in~$(a)$ with the same notations.

\smallskip 
\noindent {\sc Step~1}: First, in a discrete time setting, we assume that $0\le h_k\le H_k\!\in L^r(\P)$ and we aim at showing that by backward induction that $M_k \ge \Psi_k(X_{0:k})$ where $M_k= \E \big(\Phi(X_{0,n})\,|\, {\cal F}^Z_k\big)$.

If $k=n$, the inequality hold as an equality since $\Psi_n=\Phi$. Now assume $M_{k+1}\ge \Psi_{k+1}(X_{0:k+1)}$. Then,  like in $(a)$,
\begin{eqnarray*}
M_k&=&  \E \big(M_{k+1}\,|\, {\cal F}^Z_k\big) \\\
&\ge&   \E \big(\Phi(X_{0,k+1})\,|\, {\cal F}^Z_k\big)=  \E \big(\Phi(X_{0,k}, X_k+H_k Z_{k+1})\,|\, {\cal F}^Z_k\big)=\Big( Q_k\Psi_{k+1}(x_{0:k},x_k+. \,)(H_k)\Big)_{|x_{0:k}=X_{0:k}}\\
&\ge&\Big( Q_k\Psi_{k+1}(x_{0:k},x_k+ . \,)(h_k)\Big)_{|x_{0:k}=X_{0:k}}= \Psi_k(X_{0:k}).
\end{eqnarray*}  

\smallskip
\noindent {\sc Step~2.} This step is devoted to approximation in a bounded setting where $0\le h_t \le H_t\le K$. It  follows the lines of its counterpart in claim~$(a)$ taking advantage of the global boundedness by $K$. 

\medskip
\noindent {\sc Step~3.} This last step is devoted to the approximation procedure in the general setting. It differs from the above one  since there is no longer a deterministic upper-bound provided by the function $h\!\in L^2_{_T}$. Then, the key is to show that  the process $ \int_0^.\chi_K(H_s)dW_s$ converges for the sup norm  over $[0,T]$ in $L^{r'}(\P)$ toward the process $ \int_0^.H_sdW_s$. In fact, it follows from~\eqref{eq:BGG21} applied with $p=r'$ that 
\[
  \E \sup_{t\in [0,T]}\Big|\int_0^t H_sdW_s-\int_0^t \chi_K(H_s)dW_s\Big|^{r'} \le  c^p_p \, \E |\big(|H|-K\big)_+|^{r'}_{L^2_{_T}}.
 \]
As $|H|_{L^2_{_T}}\!\in L^{r'}(\P)$, one concludes by the Lebesgue dominated convergence theorem by letting $K\to+\infty$.~$\cqfd$
%

\bigskip
\noindent{\bf Remarks.}  $\bullet$ Step~1 can be extended to non-symmetric, centered independent random variables $(Z_k)_{1\le k\le n}$ if the sequences $(H_k)_{0\le  k\le n-1}$ and $(h_k)_{0\le k\le n-1}$ under consideration satisfy $0\le H_k\le h_k$, $k=0,\ldots,n-1$.

\smallskip
\noindent $\bullet$ When $H$ has left continuous paths, the proof can be significantly simplified by considering the simpler approximating sequence $H^{(n)}_t = \widetilde H_t^n$ which clearly converges toward $H$ $d\P\otimes dt$-$a.e.$ (and in the appropriate $L^p(dP\otimes dt)$-spaces as well).

\subsection{L\'evy-It\^o martingales}
\begin{Pro}\label{NonMarkovItoLevy} Let $Z=(Z_t)_{t\in [0,T]}$ be an integrable centered  L\'evy process with L\'evy measure $\nu$ satisfying $\nu(|x|^{p}{\bf 1}_{\{|x|\ge 1\}})<+\infty$ for a real exponent $p>1$.  Let $F :I\!\!D([0,T],\R)\to\R$ be a convex Skorokhod-continuous  functional with $(p,\|\,.\,\|_{\sup})$-polynomial growth.
Let $(H_t)_{t\in [0,T]}$ be an $({\cal F}_t)$-predictable process and let $h=(h_t)_{t\in [0,T]}\!\in\|h\|_{L_{_T}^{p\vee 2}}<+\infty$. 

\smallskip
\noindent $(a)$ If $ 0\le H_t \le   h_t$  $dt$-$a.e.$, $\P$-$a.s.$ then 
\[
\E \,F\left(\int_0^.H_sdZ_s\right)\le \E\, F\left(\int_0^.h_sdZ_s\right).
\]
If furthermore $Z$ is symmetric, the result holds as soon as $|H_t| \le   h_t$ $dt$-$a.e.$, $\P$-$a.s.$.

\smallskip
\noindent $(b)$ If  $H_t \ge h_t \ge 0$ $dt$-$a.e.$, $\P$-$a.s.$  and $|H|_{L^{p\vee2}_{_T}}\!\in L^{p}(\P)$, then 
\[
\E\, F\left(\int_0^.H_sdZ_s\right)\ge \E\,F \left(\int_0^. h_sdZ_s\right).
\]

\noindent $(c)$  If the L\'evy process $Z$ has no Brownian component, the above claims  claims $(a)$ and $(b)$ remain true if we only assume  $h\!\in L^p_{_T}$ and $|H|_{L^{p}_{_T}}\!\in L^{p}(\P)$ respectively.
\end{Pro}

\noindent {\bf Proof.} $(a)$ This proof follows the   approach introduced for the  is an extension of the Brownian-It\^o case up to the technicalities induced by L\'evy processes.

\smallskip
\noindent  {\sc Step~1} (Discrete time). This step does not differ from that developed for Brownian-It\^o martingales, except that in the the L\'evy setting we rely on claim~$(a)$ of Lemma~\ref{Lem:Jensen} since the marginal distribution of the increment of a L\'evy process has no reason to be symmetric. 
%

\smallskip
\noindent  {\sc Step~2} (Approximation-Regularization). Temporarily assume that $h$ is bounded. We consider the  approximation procedure of $H$ by stepwise constant c\`agl\`ad $({\cal F}_t)_t$-adapted (hence predictable) processes $H^{(n)}$ already defined by~\eqref{eq:Hnhn} in the proof of the previous proposition.     Then, we first consider the L\'evy-Khintchine decomposition of the L\'evy martingale $Z$ 
\[
\forall\, t\!\in [0,T], \qquad  Z_t = a\,W_t + \widetilde Z^{\eta}_t + Z^{\eta}_t, \quad a\ge 0,
\]
where $ \widetilde Z^{\eta}$ is a martingale with jumps of size at most $\eta$ and L\'evy measure $\nu(\,.\cap\{|z|\le \eta\})$ and $Z^{\eta}$ is a compensated Poisson process with (finite) L\'evy measure $\nu(. \cap \{|z|>\eta\})$. Let $n$ be a positive integer. We will perform  a ``cascade"  procedure to make $p$ decrease thanks to st the $B.D.G.$ Inequality. This --~classical~--~method is more detailed in the proof of Proposition~\ref{NonMarkovIto} in Appendix~\ref{app:B} (higher moments of L\'evy driven diffusions).

\smallskip We first assume  assume that $p\!\in (1,2]$.  Combining  Minkowski's and   $B.D.G.$'s  Inequalities yields 
\begin{eqnarray*}
\left\|\sup_{t\in [0,T]} \Big|  \int_0^t H_sdZ_s-\int_0^t H^{(n)}_{s}dZ_s\Big|\right\|_p&\le& c_p \,a
\Big\| |H-H^{(n)}|_{L^2_{_T}}  \Big\|_{p}
\\
&&+ c_p  \Big\| \sum_{0<s\le T} (H_s-H^{(n)}_{s})^2(\Delta Z_s)^2\mbox{\bf 1}_{\{|\Delta Z_s|>\eta\}}\Big\|^{\frac 12}_{\frac p2}\\
&& + c_p  \Big\| \sum_{0<s\le T} (H_s-H^{(n)}_{s})^2(\Delta Z_s)^2\mbox{\bf 1}_{\{|\Delta Z_s|\le \eta\}}\Big\|^{\frac 12}_{1}
\end{eqnarray*}
where we used in the   last line the monotony of $L^p(\P)$-norm  $\frac p2\le 1$.

 Using now  the compensation formula  and again that $\frac p2\!\in (0,1]$, it follows   
\begin{eqnarray*}
 \E\, \Big| \sum_{0<s\le T} (H_s-H^{(n)}_{s})^2(\Delta Z_s)^2\mbox{\bf 1}_{\{|\Delta Z_s|>\eta\}}\Big|^{\frac p2} &\le &  \E \sum_{0<s\le T} |H_s-H^{(n)}_s|^p |\Delta Z_s|^p \mbox{\bf 1}_{\{ |\Delta Z_s|>\eta\}}\\
 &= &\E |H -H^{(n)}  |^p_{L^p_T}  \,  \nu (|z|^p\mbox{\bf 1}_{\{|z|> \eta\}})\\
 &\le& T^{1-\frac p2}\E |H -H^{(n)}  |^p_{L^2_{_T}}  \,  \nu (|z|^p\mbox{\bf 1}_{\{|z|> \eta\}})\\
 &\le &  T^{1-\frac p2}\left(\E |H -H^{(n)}  |^2_{L^2_{_T}}  \right)^{\frac p2}  \nu (|z|^p\mbox{\bf 1}_{\{|z|> \eta\}}).
\end{eqnarray*}
On the other hand,
\[
 \E\, \Big| \sum_{0<s\le T} (H_s-H^{(n)}_{s})^2(\Delta Z_s)^2\mbox{\bf 1}_{\{|\Delta Z_s|\le \eta\}}\Big|= \E  |H -H^{(n)}  |^2_{L^2_{_T}}  \nu (z^2\wedge  \eta).
\]

We derive from~\eqref{eq:approxHh} that the above  three terms go to $0$ as $n$ goes to infinity so that 
\[
\sup_{t\in [0,T]} \Big|  \int_0^t H^{(n)}_sdZ_s-\int_0^t H_{s}dZ_s\Big|\stackrel{L^p(\P)}{\longrightarrow}0.
\]
Then, Lemma~\ref{Lem:6.3} applied to the subsequence $(m(n))_{n\ge 1}$  implies that the stepwise constant process $\displaystyle \left(\int_0^{\underline t_{m(n)}} H^{(n)}_{s}dZ_s\right)_{t\in [0,T]}$ satisfies
\[
{\rm dist}_{Sk} \left( \int_0^{\underline{.}_{m(n)}} H^{(n)}_sdZ_s, \int_0^. H_{s}dZ_s\right)\stackrel{\P}{\longrightarrow}0
\]
which in turn implies the functional $Sk$-weak convergence. Furthermore, the above $L^p$-convergence implies that  the sequence $\displaystyle\left( \sup_{t\in [0,T]}\Big|\int_0^t H^{(n)}_{s}dZ_s\Big|\right)_{n\ge 1}$  is uniformly $L^p$-integrable which is also clearly true for $\displaystyle\left( \sup_{t\in [0,T]}\Big|\int_0^{\underline{t}_{m(n)}} H^{(n)}_{s}dZ_s\Big|\right)_{n\ge 1}\!\!$. Following the same lines and still using Lemma~\ref{Lem:6.3}, we get 
\[
{\rm dist}_{Sk} \left( \int_0^{\underline{.}_{m(n)}} h^{(n)}_sdZ_s, \int_0^. h_{s}dZ_s\right)\stackrel{\P\mbox{-}a.s.}{\longrightarrow}0\;\mbox{and}\;\displaystyle\left( \sup_{t\in [0,T]}\Big|\int_0^t h^{(n)}_{s}dZ_s\Big|\right)_{n\ge 1}\; \mbox{ is uniformly $L^p$-integrable}.
\]

Since $0\le H_t\le h(t)$   $dt$-$a.e.$ $\P$-$a.s.$ (or $0\le |H_t|\le h_t$ if $Z$ is symmetric), for every fixed integer $n\ge 1$, we have, owing to Step~1 and  following the lines of Step~3 of the proof of Proposition~\ref{NonMarkovIto},
\[
\E \left(F\Big(\int_0^{\underline{t}_{m(n)}} H^{(n)}_sdZ_s\Big)_{t\in[0,T]}\right)\le \E\left(F\Big(\int_0^{\underline{t}_{m(n)}} h^{(n)}_sdZ_s\Big)_{t\in[0,T]}\right).
\]
Letting $n\to +\infty$ yields the announced result since $F$ is $Sk$-continuous with $(p,\|\,.\,\|_{\sup}$)-polynomial growth (owing to the above uniform $L^p$-integrability results).

\smallskip 
\noindent Assume now  $p>2$. First note that since $h$ is bounded one can extend~\eqref{eq:approxHh} as follows: there exists a sequence $m(n)\uparrow+\infty$ such that the processes $H^{(n)}$ and the functions $h^{(n)}$ defined by~\eqref{eq:Hnhn} satisfy
\begin{equation}\label{eq:approxHhp}
\E |H-H^{(n)} |^{p}_{L^p_{_T}} +|h-h^{(n)} |_{L^p_{_T}} \longrightarrow 0\;\mbox{ as }\; n\to +\infty.
\end{equation}
To this end, we introduce the dyadic logarithm m of $p$ $i.e.$ the integer $\ell_p$ such that where $2^{\ell_p}< p\le 2^{\ell_p+1}$.
%
Thus, if  $p\!\in (2,4]$ $i.e.$ $\ell_p=1$,
\begin{equation}\label{eq:HZBDG}
\hskip -0,75 cm \left\|\sup_{t\in [0,T]}\! \Big|  \int_0^t H_sdZ_s-\int_0^t H^{(n)}_{s}dZ_s\Big|\right\|_p\!\le\! c_p\, \left(\kappa 
\Big\| |H-H^{(n)}|_{L^2_{_T}}  \Big\|_{p} \!+ \!  \Big\| \sum_{0<s\le T} (H_s-H^{(n)}_{s})^2(\Delta Z_s)^2\Big\|^{\frac 12}_{\frac p2} \right).
\end{equation}
Now, Minkowski's Inequality applied with $\|.\|_{\frac p2}$ yields
\begin{eqnarray*}
\Big\| \sum_{0<s\le T} (H_s-H^{(n)}_{s})^2(\Delta Z_s)^2\Big\|_{\frac p2} &\le &\Big\| \sum_{0<s\le T} (H_s-H^{(n)}_{s})^2 (\Delta Z_s)^2-\nu(z^2)\int_0^T(H^{(n)}_s-H_s)^2ds \Big\|_{\frac p2} \\
&&+ \nu(z^2) \big\||H^{(n)}-H|_{L^2_{_T}}\big\|^{2}_p.
\end{eqnarray*}
In turn, the $B.D.G.$  Inequality applied to the martingale 
$$
M^{(1)}_t=  \sum_{0<s\le t} (H_s-H^{(n)}_{s})^2 (\Delta Z_s)^2-\nu(z^2)\!\!\int_0^t\!\!(H^{(n)}_s-H_s)^2ds,\; t\in [0,T],
$$ 
yields  
\begin{eqnarray*}
\Big\| \sum_{0<s\le T} (H_s-H^{(n)}_{s})^2 (\Delta Z_s)^2-\nu(z^2)\int_0^T(H^{(n)}_s-H_s)^2ds\Big)\Big\|_{\frac p2}& \le& c_{\frac p2} \Big\| \sum_{0<s\le T} (H_s-H^{(n)}_{s})^4 (\Delta Z_s)^4 \Big\|^{\frac12}_{\frac p4}\\
&\le& c_{\frac p2} \Big(\E  \sum_{0<s\le T} (H_s-H^{(n)}_{s})^p |\Delta Z_s|^p  \Big)^{\frac{2}{p}}\\
&= & c_{\frac p2} \Big (\nu(|z|^p)\E  \int_0^T |H_s-H^{(n)}_{s}|^p ds\Big)^{\frac{2}{p}}\\
&=&  c_{\frac p2} \Big ( \nu(|z|^p)\Big)^{\frac{2}{p}}\big\| |H-H^{(n)}|_{L^p_{_T}}\big\|_p^2
\end{eqnarray*}
where we successively used  that $\frac p4\le 1$ in the second line and  the compensation formula in the third  line. 
 Finally, we note that, as $p\ge2$, 
 $$ 
 \big\||H^{(n)} -H |_{L^2_{_T}} \big\|_p
\le T^{\frac 12-\frac 1p}   \big\||H^{(n)}-H|_{L^p_{_T}} \big\|_p\le T^{\frac 12-\frac 1p}   \big\||H^{(n)} -H |_{L^p_{_T}} \big\|_{p}\to 0 \mbox{ as } \; n\to +\infty
$$ 
owing to~\eqref{eq:approxHhp}. This shows that both  terms in the right hand side of~\eqref{eq:HZBDG} converge to $0$ as $n\to +\infty$,  so that 
\[
\left\|\sup_{t\in [0,T]} \Big|  \int_0^t H^{(n)}_sdZ_s-\int_0^t H_{s}dZ_s\Big|\right\|_p  \longrightarrow 0\;\mbox{ as }n\to +\infty.
\]
We show likewise
\[
\left\|\sup_{t\in [0,T]} \Big|  \int_0^t h^{(n)}_sdZ_s-\int_0^t h_{s}dZ_s\Big|\right\|_p     \longrightarrow 0\;\mbox{ as }n\to +\infty.
\]
These two convergences imply the $L^p(\P)$-uniform integrability of both sequences $\displaystyle\left( \sup_{t\in [0,T]}\Big|\int_0^t H^{(n)}_{s}dZ_s\Big|\right)_{n\ge 1}$ and $ \displaystyle\left( \sup_{t\in [0,T]}\Big|\int_0^t h^{(n)}_{s}dZ_s\Big|\right)_{n\ge 1}$. At this stage, one concludes like in the case $p\!\in(1,2]$.

In the general case,  one proceeds by a classical ``cascade" argument  based on repeated applications of the $B.D.G.$  Inequality involving  the martingales (see the proof of Proposition~\ref{pro:highmoment} in Appendix~\ref{app:B} for a more detailed implementation this cascade procedure in a similar situation)
\[
M^{(k)}_t= \sum_{0\le s\le t} (H^{(n)}_s-H_s)^{2^k} (\Delta Z_s)^{2^k}- \nu(|z|^{2^k})\!\int_0^t  (H^{(n)}_s-H_s)^{2^k}ds,\; t\ge 0,\; k=1,\ldots, \ell_p.
\]
We show  by switching from $p$ to $p/2,\, p/2^2,\ldots, p/2^k,\ldots$ until we get $p/2^{\ell_p}\!\in (1,2]$ when  $k=\ell_p$, that
\begin{eqnarray*}\label{eq:HZBDGfinal}
 \left\|\sup_{t\in [0,T]}\! \Big|  \int_0^t H_sdZ_s-\int_0^t H^{(n)}_{s}dZ_s\Big|\right\|_p&\le&  c_p\,\kappa 
\Big\| |H-H^{(n)}|_{L^2_{_T}}  \Big\|_{p} \\
&&+ \kappa_{p,\nu}\sum_{\ell=1}^{\ell_p}  \Big\| |H^{(n)}-H|_{L^{2^{\ell}}_T}\Big\|^2_{p}+\Big\| |H^{(n)}-H|_{L^{p}_T}\Big\|^2_{p}.\end{eqnarray*}
One shows likewise the   counterpart related to $h$ and $h^{(n)}$.



\smallskip 
\noindent {\sc Step~3} (Second approximation).  Now we have to get rid of the boundedness of $h$.  Like in the Brownian It\^o case, we approximate $h$ by $h\wedge K$ and $H$ by $\chi_K(H)$ where the thresholding function $\chi_K$ have been introduced in Step~3 of the proof of Theorem~\ref{FCOEuroJump} (to take into account at the same time the symmetric and the standard settings for the L\'evy process $Z$). Let $p\!\in (1,+\infty)$.

\begin{eqnarray*}
\left\|\sup_{t\in [0,T]} \Big|  \int_0^t H_sdZ_s-\int_0^t \chi_K(H_{s})dZ_s\Big|\right\|_p&\le& c_p \left(\kappa
\Big\| |H - \chi_K(H  )|_{L^2_{_T}}  \Big\|_{p} +   \Big\| \sum_{0<s\le T} (H_s- \chi_K(H_{s}))^2(\Delta Z_s)^2\Big\|_{\frac p2} \right)\\
&=&  c_p   \left(\kappa
\Big\| |(|H |-  K)_+|_{L^2_{_T}}  \Big\|_{p} +   \Big\| \sum_{0<s\le T} (|H_s|-  K)_+^2(\Delta Z_s)^2\Big\|_{\frac p2} \right)\\
&\le &c_p   \left(\kappa
   |(h -  K)_+|_{L^2_{_T}}   +   \Big\| \sum_{0<s\le T} (h_s-  K)_+^2(\Delta Z_s)^2\Big\|_{\frac p2} \right).
\end{eqnarray*}

We derive again by  this cascade argument that $ \Big\| \sum_{0<s\le T} (h_s-  K)_+^2(\Delta Z_s)^2\Big\|_{\frac p2} $ can be upper-bounded by linear combinations of quantities of the form
\[
 |(h -  K)_+|_{L^{2^k}_{_T}} \nu(z^{2^k}), \; 0\le k\le \ell_p
\]
and 
\[
\E \sum_{0<s\le T} (h_s-  K)_+^p|\Delta Z_s|^p = |(h -  K)_+|^p_{L^{p}_{_T}}\nu(|z|^p).
\]
Consequently, if $h\!\in L^p_{_T}$, all these quantities go to zero as $K\to +\infty$n owing to the Lebesgue dominated convergence theorem. In turn this implies that
\[
\left\|\sup_{t\in [0,T]} \Big|  \int_0^t H_sdZ_s-\int_0^t \chi_K(H_{s})dZ_s\Big|\right\|_p\longrightarrow 0\; \mbox{ as } K \to +\infty.
\]
The same holds with $h$ and $h\wedge K$. So it is possible to let $K$ go to infinity in the inequality
\[
\E F\Big(\int_0^.\chi_K(H_s)dZs\Big)\le \E F\Big(\int_0^.\chi_K(h_s)dZs\Big)
\]
to get the expected result.

%
\smallskip
\noindent $(b)$ is proved  adapting the  lines of the proof Proposition~\ref{NonMarkovIto}$(b)$ as we did for $(a)$. The main point is to get rid of the boundedness of $h$ $i.e.$ to obtain the conclusion of the above Step~3 without ``domination property" of $H$ by $h$. The additional assumption $|H|_{L^{p\vee2}_{_T}}\!\in L^p$ clearly yields to the expected  conclusion. 

\smallskip
\noindent $(c)$ This follows from a careful reading of the proof, having in mind that terms of the form $\Big\| |H-H^{(n)}|_{L^2_{_T}}  \Big\|_{p}
$ vanish when $\kappa=0$.~$\cqfd$

\subsubsection{Dol\'eans (Brownian) martingales}

The dame methods applied to  Dol\'eans exponential yields similar result holds   with direct  applications to the robustness of Black-Scholes formula for option pricing. First we recall that the Dol\'eans exponential of a  continuous local martingale $(M_t)_{t\in [0,T}$ is continuous local martingale defined by 
\[
{\cal E}\big(M\big)_t = e^{M_t -\frac 12 \langle M\rangle_t},\; t\!\in [0,T].
\] 
It is  a martingale on $[0,T]$ if an sonly if  $\E\, e^{M_t -\frac 12 \langle M\rangle_t}=1$. A practical  criterion, due to Novikov, says that, s  a martingale on $[0,T]$ as soon as  $\E \, e^{\frac 12\langle M\rangle_t }<+\infty$.

\begin{Pro}\label{pro:Dolexp} Let $(H_t)_{t\in [0,T]}$ and   $h=(h_t)_{t\in [0,T]}$  be like in Proposition~\ref{NonMarkovIto}.    Let $F : {\cal C}([0,T],\R_+)\to\R$ be a convex    functional with $(r,\|\,.\,\|_{\sup})$-polynomial growth ($r\ge 1$).

\smallskip
\noindent $(a)$ If ($\,|H_t| \le   h_t$ $dt$-$a.e.$) $\P$-$a.s.$, then 
\[
\E \,F\left({\cal E}\Big(\int_0^.H_sdW_s\Big)\right)\le \E\, F\left({\cal E}\Big(\int_0^.h_sdW_s\Big)\right).
\]

\noindent $(b)$ If  ($H_t \ge h_t \ge 0$ $dt$-$a.e.$)   $\P$-$a.s.$ and there exists $\varepsilon>0$ such that 
$$
\E \Big(e^{\frac{r^2+\varepsilon}{2}|H|_{L^2_{_T}}}\Big)< +\infty,
$$
then 
\[
\E\, F\left(\!{\cal E}\Big(\int_0^.H_sdW_s\Big)\!\right)\ge \E\,F \left(\!{\cal E}\Big(\int_0^.h_sdW_s\Big)\!\right).
\]
\end{Pro}

\noindent {\bf Proof.}  $(a)$ {\sc Step~1}: For a fixed integer $n\ge 1$, we consider the sequence of random variables $(\Xi^n_k)_{k=0:n}$ recursively defined in a forward way by 
\[
 \Xi^n_0=1\quad\mbox{and} \quad  \Xi^n_{k}= \Xi^n_{k-1}\exp{\Big(H_{t^n_{k-1}}\Delta W_{t^n_k} -\frac{T}{2n}H^2_{t^n_{k-1}} \Big)}, \quad k=1,\ldots,n,
\]
(where $\Delta W_{t^n_k} =  W_{t^n_k} - W_{t^n_{k-1}} $) and the  sequences $(\xi^{n,k}_{\ell})_{\ell= k:n}$defined, still  in a  recursive forward way, by  
\[
\xi^{n,k}_{k} =1,\; \xi^{n,k}_{\ell}= \xi^{n,k}_{\ell-1} \exp{\Big(h_{t^n_{\ell-1}}\Delta W_{t^n_{\ell}} -\frac{T}{2n}h^2_{t^n_{\ell-1}} \Big)},  \; \ell=k+1,\ldots,n.
\]
We denote  by $\widetilde Q^{(n)}$ the operator  defined on Borel functions $f:\R_+\to \R$ with polynomial growth by 
\[
\forall\, x,\, h\!\in \R_+,\quad \widetilde Q^{(n)}(f)(x,h) = \E\, f\Big(x\exp{\big(hW_{\frac Tn}-\frac{T}{2n}h^2\big)}\Big).
\]
It is clear that 
$\Big(\exp{\big(hW_{\frac Tn}-\frac{T}{2n}h^2\big)}\Big)_{h\ge 0}$ is increasing for the convex order ($i.e.$ a peacock as already mentioned on the introduction  since $\displaystyle \exp{\big(hW_{\frac Tn}-\frac{T}{2n}h^2\big)}\stackrel{d}{\sim} \exp{\big(W_{h^2\frac Tn}-\frac{1}{2}\frac Tn h^2\big)}$ and $(e^{W_u-\frac u2})_{u\ge 0}$ is a martingale.
Hence, as soon as $f$ is convex,   
\begin{equation}\label{eq:tildeQn}
h\mapsto \widetilde Q^{(n)}(f)(x,h) \mbox{ satisfies the maximum principle  $i.e.$  is even and)non-decreasing on }\; \R_+.
\end{equation}
In turn, it implies that the function $(x,h)\mapsto \widetilde Q^{(n)}(f)(x,h) $ is convex on $ \R\times \R_+$ since for every $x,x'\!\in \R_+$, $h,h'\!\in \R$, $\lambda\!\in [0,1]$,
\begin{eqnarray*}
&&\E\, f\Big(\lambda x\exp{\big(\lambda hW_{\frac Tn}-\frac{T}{2n}(\lambda h)^2\big)}+(1-\lambda) x'\exp{\big((1-\lambda) h'W_{\frac Tn}-\frac{T}{2n}((1-\lambda) h')^2\big)} \Big) \\
&&\hskip 2 cm \le \lambda \E\,f\Big(x\exp{\big(\lambda hW_{\frac Tn}-\frac{T}{2n}(\lambda h)^2\big)}\Big)+(1-\lambda)\E\, f\Big(x'\exp{\big((1-\lambda)h'W_{\frac Tn}-\frac{T}{2n}((1-\lambda)h')^2\big)}\Big)\\
&&\hskip 2 cm \le \lambda \E\,f\Big(x\exp{\big( |h|W_{\frac Tn}-\frac{T}{2n}h^2\big)}\Big)+(1-\lambda)\E\, f\Big(x'\exp{\big(|h'|W_{\frac Tn}-\frac{T}{2n}(h')^2\big)}\Big)\\
&& \hskip 2 cm =\lambda \E\,f\Big(x\exp{\big( hW_{\frac Tn}-\frac{T}{2n}h^2\big)}\Big)+(1-\lambda)\E\, f\Big(x'\exp{\big(h'W_{\frac Tn}-\frac{T}{2n}(h')^2\big)}\Big)
\end{eqnarray*}
where we used the convexity of $f$ in the first inequality and~\eqref{eq:tildeQn} in the second  one. 
From now on, we consider the discrete time filtration  ${\cal G}^n_k={\cal F}^W_{t^n_k}$ and set $\E_{k}=\E(\,.\, |{\cal G}^n_k)$.  

We temporarily assume that for every $k=0,\ldots,n$, $|H_{t^n_k}|\le h_{t^n_k}$ $\P$-$a.s.$.
Let $F:{\cal C}([0,T],\R)\to \R$ be a (Borel) functional with $(r,\|\,.\,\|_{\sup})$-polynomial  growth  and let $F_n =F\circ i_n$. We will show by induction that, for every $k\!\in \{1,\ldots,n\}$, 
\begin{equation}\label{eq:bof}
\E_{k-1} F_n(\Xi^n_{0:k-1}, \Xi^n_{k}\,\xi^{n,k}_{k:n})\le \E_{k-1} F_n(\Xi^n_{0:k-2}, \Xi^n_{k-1}\xi^{n,k-1}_{k-1:n})
 \end{equation}
with the obvious convention  $\Xi^n_{0:-1}=\emptyset$. Starting from the identity
\begin{eqnarray*}
 F_n(\Xi^n_{0:k-1},\Xi^n_k  \,\xi^{n,k}_{k:n}) & =  &  F_n\Big(\Xi^n_{0:k-1},\Xi^n_{k-1}\exp{\big(H_{t^n_{k-1}}\Delta W_{t^n_{k}} -\frac{T}{2n}H^2_{t^n_{k-1}} \big)}  \xi^{n,k}_{k:n}\Big), \\
 \mbox{we derive } \hskip 2,75 cm &&\\
\E_{k-1} F_n(\Xi^n_{0:k-1},\Xi^n_k  \,\xi^{n,k}_{k:n}) & = & \left(\E \Big(F(x_{0:k-1}, x_{k-1}\exp{\big(\eta\Delta W_{t^n_{k}} -\frac{T}{2n}\eta^2 \big)}   \xi^{n,k}_{k:n})\Big)\right)_{|x_{0:k-1}= \Xi^n_{0:k-1}, \eta = H_{t^n_{k-1}}}
\end{eqnarray*}
since  $(\Xi^n_{0:k-1}, H_{t^n_{k-1}})$ is ${\cal G}^n_{k-1}$-measurable and $(\Delta W_{t^n_k}, \xi^{n,k}_{k:n})$ is independent of ${\cal G}^n_{k-1}$. Now set, for every $x_{0:k-1}\!\in \R_+^k$, $\widetilde x_{k}\!\in \R_+$,
\[
G_{n,k}(x_{0:k-1}, \widetilde x_{k})= \E F_n\big(x_{0:k-1}, \widetilde x_{k} \,\xi^{n,k}_{k:n}\big)
\]
so that 
\[
\widetilde Q^{(n)}(G_{n,k}(x_{0:k-1},.))(x_{k-1},\eta) = \E\, F_n \Big(x_{0:k-1},x_{k-1}\exp{\big(\eta\Delta W_{t^n_{k}} -\frac{T}{2n}\eta^2 \big)}   \xi^{n,k}_{k:n}\Big).
\]
The function $F_n$ being convex on $\R_+^{n+1}$, it is clear that $G_{n,k}$ is convex on $\R_+^{k+1}$ as well. It is in particular convex in the variable $\widetilde x_{k}$ which in turn implies by~\eqref{eq:tildeQn} that $\eta\mapsto\widetilde Q^{(n)} (G_{n,k}(x_{0:k-1},.))(x_{k-1},\eta)$ satisfies the maximum principle $i.e.$  is even and convex. As a consequence, $|H_{t^n_{k-1}}|\le h_{t^n_{k-1}}$ implies 
\begin{eqnarray*}
\E_{k-1} F_n(\Xi^n_{0:k-1},\Xi^n_k  \,\xi^{n,k}_{k:n}) &=& \left[ \widetilde Q^{(n)}\big(G_{n,k}(x_{0:k-1}, \,.\,)\big)(x_{k-1},\eta) \right ]_{|x_{0:k-1}= \Xi^n_{0:k-1}, \eta = H_{t^n_{k-1}}} \\
&=& \left[\widetilde Q^{(n)}\big(G_{n,k}(x_{0:k-1}, \,.\,)\big)(x_{k-1},\eta)\right]_{|x_{0:k-1}= \Xi^n_{0:k-1}, \eta = |H_{t^n_{k-1}}|} \\
&\le & \left[\widetilde Q^{(n)}\big(G_{n,k}(x_{0:k-1}, \,.\,)\big)(x_{k-1},\eta)\right]_{|x_{0:k-1}= \Xi^n_{0:k-1}, \eta = h_{t^n_{k-1}}} \\
&=&  \E_{k-1} \Big(F_n\big(\Xi^n_{0:k-1}, \Xi^n_{k-1}\exp{\big( h_{t^n_{k-1}}\Delta W_{t^n_{k}} -\frac{T}{2n}h_{t^n_{k-1}}^2 \big)}  \xi^{n,k}_{k:n}\big)\Big)\\
&=&   \E_{k-1} \Big(F_n\big(\Xi^n_{0:k-2}, \Xi^n_{k-1} \xi^{n,k-1}_{k-1:n}\big)\Big)
\end{eqnarray*}
where we  used once again  that $\xi^{n,k}_{k:n}$ is independent of ${\cal G}^n_{k-1}$ in the penultimate line.

One derives by taking expectation of the resulting  inequality that the sequence $\E_{k-1} F_n(\Xi^n_{0:k-1},\Xi^n_k  \,\xi^{n,k}_{k:n})$, $k=1:n$,  is non-increasing. Finally, by comparing the terms for $k=n$ and $k=0$, we get
\[
\E F\big(X^{n,n}\big) = \E\, F_n(\Xi^n_{0:n})  \le  \E\, F_n(\xi^{n,0}_{0:n})= \E F\big(X^{n,0} \big).
\]


\smallskip
\noindent {\sc Step 2} (Approximation-Regularization). We closely  follow the approach developed in  Steps~2 and~3 of Proposition~\ref{NonMarkovIto}. First, we temporarily assume that $h$ is bounded by a real constant $K$ and we introduce the stepwise constant c\`agl\`ad  processes $(H^{(n)})_{t\in [0,T]}$ and $(h^{(n)}_t)_{t\in[0,T]}$ defined by~\eqref{eq:Hnhn}  (and  satisfying~\eqref{eq:approxHh}), namely
\[
\Big\| |H^{(n)}  -H |_{L^2_{_T}} \Big\|_{_2} +\big|h^{(n)}-h\big|_{L^2_{_T}}\longrightarrow 0\;\mbox{ as }\; n \to +\infty.
\]
In particular 
\[
\sup_{t\in[0,T]}\Big| \int_0^t(H^{(n)}_s)^2ds -\int_0^t H^2_sds\Big|\le  2K \big |H^{(n)}-H\big |_{L^1_{_T}}\le  2K\sqrt{T}  \big |H^{(n)}-H\big |_{L^2_{_T}}.
\]
As a consequence
\begin{eqnarray*}
\sup_{t\in[0,T]} \left |   \int_0^tH^{(n)}_sdW_s-\frac 12 \int_0^t (H_s^{(n)})^2ds -\left( \int_0^tH_sdW_s-\frac 12 \int_0^t H_s^2ds\right)\right|&\le&\sup_{t\in[0,T]} \left |   \int_0^t(H^{(n)}_s-H_s)dW_s \right|\\\
&&+  K\sqrt{T}  \big |H^{(n)}-H\big |_{L^2_{_T}}.
\end{eqnarray*}
Set for notational convenience 
\[
X^{(n)}_t ={\cal E}\Big(\int_0^.H^{(n)}_sdW_s\Big)_t\quad\mbox{ and }\quad X_t ={\cal E}\Big(\int_0^.H_sdW_s\Big)_t,\;t\!\in [0,T].
\]
which are both  true martingales owing to Novikov's criterion. The above inequality combined with Doob's Inequality implies that 
\[
\sup_{t\in[0,T]} \Big|\log X^{(n)}_t -\log X_t\Big| \stackrel{L^2}{\longrightarrow} 0\; \mbox{ as } \;n\to +\infty.
\]
As a consequence,  $\displaystyle X^{(n)}\stackrel{{\cal L}(\|\,.\,\|_{\sup})}{\longrightarrow} X$ since the exponential function is continuous. Denoting by $x^{(n)}$ and $x$ the counterpart of these processes for the functions $h^{(n)}$ and $h$, we   get likewise  $\displaystyle x^{(n)}\stackrel{{\cal L}(\|\,.\,\|_{\sup})}{\longrightarrow} x$. Owing once again to Lemma~\ref{Interpol}, the continuity of the exponential again, and the chain rule for weak convergence, we finally get
\[
e^{I_{m(n)}(\log X^{(n)})} \stackrel{{\cal L}(\|\,.\,\|_{\sup})}{\longrightarrow} e^{\log X}= X\;\mbox{ and }\;e^{I_{m(n)}(\log x^{(n)})} \stackrel{{\cal L}(\|\,.\,\|_{\sup})}{\longrightarrow} e^{\log x}= x \;\mbox{ as } \;n\to +\infty.
\]

Applying Step~1 with $X^{(n)}$ and $x^{(n)}$
\[
\forall\, n\!\in \N,\quad \E\, F(X^{(n)}) \le \E\, F(x^{(n)}).
\]
To let $n$ go to infinity in this inequality, we again need  a uniform integrability argument namely  that  $\|X^{(n)}\|_{\sup} $ and $\|x^{(n)}\|_{\sup} $ are both $L^p$-bounded for a $p>r$ since the functional $F$ has at most a $(r,\|\,.\,\|_{\sup})$-polynomial growth. So, let $p>r\vee1$. It follows from Doob's Inequality applied to the non-negative sub-martingale$(X^{(n)})^p$ that 
\begin{eqnarray*}
\E\Big(\sup_{t\in[0,T]} (X^{(n)}_t)^p\Big)& \le& \Big(\frac{p}{p-1}\Big)^p\E \big(X^{(n)}_{_T})^p\\
&\le & \Big(\frac{p}{p-1}\Big)^p \E\left( {\cal E}\big(p\int_0^. H^{(n)}_s dW_s\big)_{_T}\right)  e^{\frac{p(p-1)}{2} \int_0^T (H^{(n)}_s)^2ds}\\
&\le & \Big(\frac{p}{p-1}\Big)^p e^{\frac{p(p-1)}{2} K^2T}
\end{eqnarray*}
where we used that  $ \Big({\cal E}\Big(p\int_0^.H^{(n)}_sdW_s\Big)_t \Big)_{t\ge 0}$ is a true martingale (owing to Novikov' criterion).
%
%
The case of $F(x^{(n)})$ follows likewise.

\smallskip
\noindent {\sc Step~3}: The extension to $h\!\in L^2_{_T}$ is similar that performed in  the former propositions: first note that 
\[
{\cal E}\Big(\int_0^.\chi_K(H_s)dW_s\Big)\stackrel{{\cal L}(\|\,.\,\|_{\sup})}{\longrightarrow} {\cal E}\Big(\int_0^.H_sdW_s\Big)\quad\mbox{as }\; K\to +\infty.
\]
  The uniform integrality  of $\displaystyle \sup_{t\in[0,T]}  {\cal E}\Big(\int_0^.\chi_K(H_s)dW_s\Big)_t$ as $K$ grows to infinity follows form its   $L^p(\P)$- boundedness   for a $p\!\in (1,+\infty)$ which in turn is a consequence of Doob's inequality:     
\begin{eqnarray*}
\E \sup_{t\in[0,T]}\left ( {\cal E}\Big(\int_0^.\chi_K(H_s)dW_s\Big)_t\right)^p&\le& \Big(\frac{p}{p-1}\Big)^p \E\, {\cal E}\Big(\int_0^.\chi_K(H_s)dW_s\Big)^p_T\\
&\le &  \Big(\frac{p}{p-1}\Big)^p e^{\frac{p(p-1)}{2} \int_0^T \chi_K^2(h_s)ds} \, \E \,{\cal E}\Big(p\int_0^.\chi_k(H_s)dW_s\Big)_T  \\
&= &    \Big(\frac{p}{p-1}\Big)^p   e^{\frac{p(p-1)}{2} \int_0^T \chi_K^2(h_s)ds} \\
&\le &   \Big(\frac{p}{p-1} \Big)^p  e^{\frac{p(p-1)}{2}|h|_{L^2_{_T}}}<+\infty
\end{eqnarray*}
which yields $L^p$ -boundedness with respect to the threshold $K$.

\smallskip
\noindent $(b)$  The discrete time part can be established by adapting  item~$(a)$ in the spirit of  Proposition~\ref{NonMarkovIto}$(b)$. The approximation step follows like above as well, except for the final uniform integrability argument which needs  specific care.  It suffices to show that for an $r'>r$, $\sup_{t\in[0,T]} {\cal E}\Big(\int_0^.\chi_K(H_s)dW_s\Big)_t$ is $L^{r'}$-bounded as $K\to +\infty$.
 
 \smallskip 
 \noindent $\rhd$ If $r\!\in(0,1)$, one may choose $r'\!\in (r,1)$. Then,  one checks that
 \[
\E\left[ \sup_{t\in[0,T]} {\cal E}\Big(\int_0^.\chi_K(H_s)dW_s\Big)^{r'}_t\right]\ le \E\left[ \sup_{t\in[0,T]} {\cal E}\Big(\int_0^.\chi_K(H_s)dW_s\Big)_t\right].
 \]
Now if $p>1$,  Doob's Inequality implies 
\[
\E\left[ \sup_{t\in[0,T]} {\cal E}\Big(\int_0^.\chi_K(H_s)dW_s\Big)^p_t\right] \le \Big(\frac{p}{p-1}\Big)^p \E\left[ \sup_{t\in[0,T]} {\cal E}\Big(\int_0^.\chi_K(H_s)dW_s\Big)^p_{_T}\right] \le   \Big(\frac{p}{p-1}\Big)^p.
\]
which yields the announced result.

\smallskip
\noindent $\rhd$ If $r\ge 1$, let $p>r$. Combining  successively   Doob's Inequality and H\"older's Inequality, for every $p>r\vee1$ and every H\"older conjugate exponents $\lambda,\mu= \frac{\lambda}{\lambda-1}>1$,  leads to 
\begin{eqnarray*}
\E\, \sup_{t\in[0,T]}\left ( {\cal E}\Big(\int_0^.\chi_K(H_s)dW_s\Big)_t\right)^p&\le&\Big(\frac{p}{p-1}\Big)^p \E\, {\cal E}\Big(\int_0^.\chi_K(H_s)dW_s\Big)^p_T\\
&\le&    \Big(\frac{p}{p-1}\Big)^p  \left[\underbrace{\E\, {\cal E}\Big(\lambda p\int_0^.\chi_K(H_s)dW_s\Big)_T}_{=1}\right]^{\frac{1}{\lambda}} \left[ \E\, e^{\frac{\lambda p^2}{2} \int_0^T \chi_K(H_s)^2ds}\right]^{\frac{\lambda-1}{\lambda}}\\
&\le &    \Big(\frac{p}{p-1}\Big)^p  \left[ \E\, e^{\frac{\lambda p^2}{2} | H|_{L^2_{_T}}^2}\right]^{\frac{\lambda- 1}{\lambda}}.
\end{eqnarray*}
Consequently, for $\lambda$ close enough to $1$ and $p$ close enough to $r$, we have $\lambda p^2\le r^2+\varepsilon$ which ensures the $L^p(\P)$-boundedness as $K\uparrow +\infty$.~$\cqfd$ 

%
%

\subsubsection{A counter-example}
The counter-example below shows that Theorem~\ref{NonMarkovIto} is no longer  true if we relax the assumption that the dominating process $(h_t)_{t\in [0,T]}$ is deterministic.

\medskip Let $X=X^{\sigma}= (X^{\sigma}_{0:2})$ be a two period process satisfying  
\begin{eqnarray*}
X_0= 0,\quad X_1=  \sigma Z_1&\mbox{ and } &X_2 = X_1 + \sqrt{2v(Z_1)}Z_2
\end{eqnarray*}
where $Z_{1:2}\stackrel{\cal L}{\sim} {\cal N}(0;I_2)$, $\sigma\ge 0$,  and $\varphi: \R\to \R_+$ is a bounded non-increasing function.

\smallskip
Let $f(x) =e^x$ and let $\varphi: \R_+\to \R$ be the function defined by 
\[
\varphi(\sigma) := \E f(X_2)=\E \big(e^{\sigma Z_1+v(Z_1)}\big).
\]
Differentiating $\varphi$ yields 
\[
\varphi'(\sigma) = \E\big(e^{\sigma Z_1+v(Z_1)}Z_1\big)
\]
so that
\[
\varphi'(0) = \E \big(e^{v(Z_1)}Z_1\big)< \E \,e^{v(Z_1)} \E \, Z_1=0
\]
by a standard one-dimensional co-monotony argument:~both functions $z\!\mapsto \!e^{v(z)},\,z\!\mapsto \!z$ are  non-decrea\-sing which implies $\varphi'(0)\!\le\! 0$ but  none of them  are $\P_{Z_1}\!$-$a.s.$ constant, hence   equality cannot hold. As a consequence, $\varphi $ is (strictly) decreasing on a right neighbourhood $[0, \sigma_0]$, $\sigma_0>0$, of $0$.

To include this into  a Brownian stochastic integral framework, one proceeds as follows: let $W$ be a standard Brownian motion and  $\sigma, \widetilde \sigma \!\in (0,\sigma_0]$, $\sigma < \widetilde \sigma$. 
\[
H_t = \sigma \mbox{\bf 1}_{[0,1]} (t) + \sqrt{2v(W_1)}\mbox{\bf 1}_{(1,2]}(t), \;  \widetilde H_t = \widetilde \sigma \mbox{\bf 1}_{[0,1]} (t)+ \sqrt{2v(W_1)}\mbox{\bf 1}_{(1,2]}(t).
\]

It is clear that $0\le H_t \le \widetilde H_t$, $t\!\in [0,2]$, whereas
\[
\E \Big( e^{\int_0^2H_sdW_s}\Big) > \E \Big( e^{\int_0^2\widetilde H_sdW_s}\Big).\]
This makes up a counter-example to the conclusion of Proposition~\ref{NonMarkovIto}.

\medskip It has to be noted that if the function $v$ is non-decreasing, then choosing $f(x)=e^{-x}$ leads to a similar result since 
\[
\psi(\sigma) := \E f(X_2)=\E \big(e^{-\sigma Z_1+v(Z_1)}\big)
\]
satisfies $\Psi'(\sigma) = -\E \big(e^{-\sigma Z_1 +v(Z_1)}\big)$. In particular one still has by a co-monotony argument that $\psi'(0)<0$ since $v$ is not constant. 

\subsubsection{A comparison theorem for Laplace transforms of Brownian stochastic integrals}

Applying our paradigm, we start by a discrete time result with  its own interest for applications. 

\begin{Pro}\label{ComparLaplW} Let $(Z_k)_{1\le k\le n}$ be a sequence of ${\cal N}(0;1)$-random variables. We set $S_0=0$ and $S_k=Z_1\dots+Z_k$, $k\!=\!1,\ldots,n$ (partial sums). We consider the two discrete time  stochastic integrals
\[
X_k = \sum_{\ell=1}^k f_\ell(S_{\ell-1})Z_\ell\quad \mbox{and}\quad Y_k= \sum_{\ell=1}^k g_\ell(S_{\ell-1})Z_\ell, \quad k=1,\dots,n,\quad X_0=Y_0=0
\]
where $f_k,\,g_k:\R\to \R_+$, $k=1,\ldots,n$ are non-negative Borel functions satisfying:

\bigskip
\centerline{either all  $f_k$, $k=1,\ldots,n$, are non-decreasing or all   $g_k$, $k=1,\ldots,n$, are non-decreasing.}

\bigskip
If, furthermore, $0\le f_k\le g_k$ for all $k=1, \ldots,n$, then
$$
\forall\, \lambda \ge 0,\quad \E\, e^{\lambda X_n}\le \E\,e^{\lambda Y_n}.
$$
\end{Pro}

\noindent {\bf Proof.} We start from the Cameron-Martin identity which reads on Borel function $\varphi:\R\to \R$

\[
\forall\sigma\!\in \R,\quad \E\, e^{\sigma Z +\varphi(Z)}= e^{\frac{\sigma^2}{2}}\E \,e^{\varphi(Z+\sigma)}\le +\infty.
\]

First, we define in a backward way   functions $\widetilde f_k$ and $\widetilde g_k$, $k=1,\ldots,n+1$ by $\widetilde f_{n+1}=\widetilde g_{n+1}\equiv 0$,
\begin{equation}\label{eq:expStchInt}
\widetilde f_k(x) = \frac{\lambda^2}{2}f_{k}^2(x) + \log \E\Big(e^{\widetilde f_{k+1}(x+\lambda f_k(x) +Z)  }\Big),\; k=0,\ldots,n,
\end{equation}
where $Z\sim {\cal N}(0;1)$. The functions $\widetilde g_k$ are defined from the $g_k$ the same way round. Then, relying on the chaining rule for conditional expectations, we check by a backward  induction  that  
\[
\E \,e^{\lambda X_n} = \E \,e^{\lambda X_{k}+ \widetilde f_{k+1}(S_k)}, \; k=1,\ldots,n. 
\]
In particular, when $k=0$, we get 
$$
\displaystyle \E \,e^{\lambda X_n}  =  e^{\widetilde f_1(0)}.
$$
It follows from~(\ref{eq:expStchInt}) and  a second backward induction that, if the  functions $f_k$  are non-decreasing for every $k=1,\ldots,n$, so are the functions $\widetilde f_k$. The same holds for $\widetilde g_k$ with respect to the functions $g_k$. Assume $e.g.$ that all the functions $\widetilde f_k$ are non-decreasing. Then,  a (third) backward induction shows: $\widetilde f_k\le \widetilde g_k$ since $\widetilde f_n\le \widetilde g_n$,  and,  if $\widetilde f_{k+1}\le \widetilde g_{k+1}$, then  for every $x\!\in \R$,
\[
\widetilde f_{k+1}(x+\lambda f_k(x)+Z) \le \widetilde f_{k+1}(x+\lambda g_k(x) +Z)\le \widetilde g_{k+1}(x+\lambda g_k(x) +Z).
\]
Plugging this inequality in~(\ref{eq:expStchInt}) combined with $f^2_k\le g^2_k$, one concludes that  $\widetilde f_k \le \widetilde g_k$. A similar reasoning can be carried out if the functions $\widetilde g_k$ are non-decreasing.~$\cqfd$

\bigskip
By the standard weak approximation method detailed in the former results, we derive the following  continuous time version of this result involving (non-decreasing) {\em completely monotone}   functions defined below.
\begin{Dfn} A {\em non-decreasing} function $\varphi:\R\to\R$ is completely  monotone  if it is the Laplace transform of a non-negative Borel measure $\mu$ supported by the non-negative real line, namely
\[
\forall\, x\!\in \R,\quad \varphi(x) = \int_{\R_+} e^{\lambda x}\mu(d\lambda).
\]
\end{Dfn}

\begin{Thm} Let $f,g: [0,T]\times \R\to \R_+$ two bounded Borel functions such that 
\begin{equation}\label{eq:Condfg}
\left\{\begin{array}{ll}
(i)& f, \,g \mbox{ are } dt\otimes dx\mbox{-$a.e.$ continuous},\\
(ii)& 0\le f\le g,\\
(iii)& \Big(\forall\, t\!\in [0,T],\; f(t,.) \mbox{ is non-decreasing}\Big) \mbox{ or }\Big(\forall\, t\!\in [0,T],\; g(t,.) \mbox{ is non-decreasing}\Big).
\end{array}\right.
\end{equation}
Then,
\[
\forall\, \lambda \ge 0,\quad \E\, e^{\lambda\int_0^T f(t,W_t)dW_t}\le \E\, e^{\lambda\int_0^T g(t,W_t)dW_t}
\]
so that, for every non-decreasing {\em completely monotone } function $\varphi:\R\to \R_+$ 
\[
 \E\, \varphi\left(\int_0^T f(t,W_t)dW_t\right)\le \E\, \varphi \left(\int_0^T g(t,W_t)dW_t\right).
\]
\end{Thm}

\noindent {\bf Remarks.} $\bullet$ The finiteness of these integrals follows from Novikov's criterion.
 
\noindent $\bullet$ One derives from~(\ref{eq:Condfg}) the seemingly more general result
\begin{equation}\label{eq:Condfgbis}
\left\{\begin{array}{ll}
(i)& f, g \mbox{ are } dt\otimes dx\mbox{-$a.e.$ continuous},\\
(ii)& \exists\, h:[0,T]\times \R\to \R_+\mbox{ such that }\left\{\begin{array}{ll}(a)&  0\le f\le h\le g \mbox{ and} \\
(b)& \forall\, t\!\in [0,T],\; h(t,.) \mbox{ is non-decreasing}.\end{array}\right.
\end{array}\right.
\end{equation}

\noindent{\bf Proof.} Assume $e.g.$ that $f(t,.)$ is non-decreasing for every $t\!\in [0,T]$. First note that by Fubini's Theorem  and It\^o's isometry
\[
\Big\|  \int_0^T f(s,W_s)dW_s -  \int_0^T f(\underline s_n,W_{\underline s_n})dW_s \Big\|_2^2 = \int_0^T \E \big( f(s,W_s)-f(\underline s_n,W_{\underline s_n})\big)^2 ds.
\]
Now, if we denote $C_s=\{x\!\in \R\,|\, f \mbox{ is continuous at } (s,x)\}$ for every $t\!\in [0,T]$, it follows from  Assumption~(\ref{eq:Condfg})$(i)$  that $\lambda(^cC_s)=0$  $ds$-$a.e.$ still by Fubini's Theorem.  As $\P_{X_s}$ is equivalent to the Lebesgue measure, one derives that $\P_{s}(C_s)= 1$ $ds$-$a.e.$. As a consequence, $\E \big( f(s,W_s)-f(\underline s_n,W_{\underline s_n})\big)^2\to 0$ $ds$-$a.e.$ as $n\to +\infty$ since $(\underline s_n, W_{\underline s_n}) \to (s,W_s)$. One concludes by the dominated Lebesgue theorem that $\Big\|  \int_0^T f(s,W_s)dW_s -  \int_0^T f(\underline s_n,W_{\underline s_n})dW_s \Big\|_2\to 0$  since $f$ is bounded. 


\smallskip Now, define for every $k=1,\ldots,n$,  
$$
X_k = \int_0^{t^n_k}f(\underline s_n, W_{\underline s_n})dW_s =\sum_{\ell=1}^k \sqrt{\frac Tn} f(t^n_{\ell-1}, W_{t^n_{\ell-1}})U^n_{\ell}
$$
where $U^n_{\ell}= \sqrt{\frac nT} (W_{t^n_\ell}-W_{t^n_{\ell-1}})$, $\ell= 1,\ldots,n$. We define likewise $(Y_k)_{k=0:n}$ with respect to the function $g$. It is clear that both $(X_k)$ and $(Y_k)$ satisfy the assumptions of the above Proposition~\ref{ComparLaplW} so that 
\[
\forall\, \lambda \ge 0, \qquad \E \,e^{\lambda  \int_0^{T}f(\underline s_n, W_{\underline s_n})dW_s} \le  \E \,e^{\lambda \int_0^{T}g(\underline s_n, W_{\underline s_n})dW_s}.
\]
One concludes by combining the above quadratic (hence weak) convergence and  the uniform integrability argument which follows from
\[
\forall\, \lambda >0, \qquad \sup_n \E \,e^{\lambda  \int_0^{T}f(\underline s_n, W_{\underline s_n})dW_s} \le e ^{\frac{\lambda^2}{2}\|f\|_{\rm sup} T}<+\infty.\quad \cqfd
\] 

\section{Convex order for the {\em r\'eduite} and  applications to path-dependent American options}\label{sec:Amer}
In this section, we aim at applying the methodology developed in the former sections to Optimal Stopping Theory, $i.e.$, as far as financial applications are concerned, to Bermuda and American style options. For general background on Optimal Stopping theory, we refer to~\cite{NEV} (Chapter~VI) and \cite{DCD} (Chapter 5.1) in discrete time and, among others,  to~\cite{NEK, KASH, SHI} in continuous time. For a discussion (and results) on comparison methods for American option prices, which usually includes an analytic component involving variational inequalities,  we refer to~\cite{BRusch1} and the references therein.

\subsection{Bermuda options}
We start from  the  discrete time dynamics introduced in the  ``European" framework. Let $(Z_k)_{1\le k\le n}$ be a sequence of independent $\R^d$-valued random vectors satisfying $Z_k\!\in L^r(\Omega, {\cal A}, \P)$, $r\ge 1$ and $\E\,Z_k=0,\, k=1,\ldots,n$. Let $(X_k)_{0\le k\le n}$ and $(Y_k)_{0\le k\le n}$ be the two sequences of random vectors defined by~\eqref{modele} $i.e.$ 
\[
X^x_{k+1}= X^x_k + \sigma_k(X^x_k) Z_{k+1}, \quad Y^x_{k+1}= Y^x_k + \theta_k(Y^x_k) Z_{k+1}, \; 0\le k\le n-1,\; X^x_0=Y^x_0=x
\]
where $\sigma_k$, $\theta_k$, $k=0,\ldots,n$ are functions from $\R$ to $\R$, all  with linear growth. This implies by a straightforward  induction that the random variables $X^x_k$ and $Y^x_k$ all lie in $L^r$ since, $e.g.$, $\sigma_k(X^x_k)$ are adapted to ${\cal F}^Z_k$ hence independent of $Z_{k+1}$, $k=0,\ldots,n-1$.  

Let ${\cal F}=({\cal F}_k)_{0\le k\le n}$ and ${\cal G}=({\cal G}_k)_{0\le k\le n}$ two filtrations on $(\Omega, {\cal A}, \P)$ such that $X^x$ is ${\cal F}$-adapted and $Y^x$ is ${\cal G}$-adapted. Let $F_k: \R^{k+1}\to \R_+$, $k=0,\ldots,n$ be a sequence of non-negative   functions with {\em $r$-polynomial growth} ($i.e.$ $0\le F_k(x_{0:k})\le C(1+|x_{0:k}|^r)$, $k=0,\ldots,n$), $r\ge 1$. Then  the processes $\big(F_k(X^x_{0:k})\big)_{0\le k\le n }$ and  $(F_k(Y^x_{0:k}))_{0\le k\le n}$  are called  {\em payoff} or {\em obstacle} processes (${\cal F}$-adapted and  ${\cal G}$-adapted respectively).

\smallskip 
We define the  ${\cal F}$- and $
{\cal G}$-``{\em r\'eduites}" associated to these payoff processes by  
\[
u_0(x) =  \sup\big\{ \E\,F_{\tau}(X^x_{0:\tau}), \, \tau \;{\cal F}\mbox{-stopping time}\big\} 
\quad
\mbox{and}
\quad
v_0(x)=\sup\big\{\E \,F_{\tau}(Y^x_{0:\tau}), \, \tau \;{\cal G}\mbox{-stopping time}\big\}
\]
respectively. These quantities are closely related to the optimal stopping problems attached to these dynamics since they represent the supremum of possible gains among ``honest" stopping strategies  ($i.e.$ non-anticipative with respect to the filtration) in a game where one wins $F_k(X^x_{0;k})$ ($F_k(Y^x_{0:k})$ respectively) when leaving the game at time $k$. Owing to the dynamic programing formula (see the proof of Proposition below) and the Markov property shared by both dynamics $X^x$ and $Y^x$, it is clear that we may assume without loss of generality that ${\cal F}= {\cal F}^X$ (natural filtration of $X^x$) and ${\cal G}= {\cal F}^{Y}$  (natural filtration of $Y^x$)  or even  ${\cal F}= {\cal G}={\cal F}^Z$  without changing the value of the {\em r\'eduites}. 

\medskip
The proposition below is the counterpart of Proposition~\ref{FCOEuroDiscTime} in discrete time for ``European" options. 

\begin{Pro}\label{Snelldiscret} Let $F_k: \R^{k+1}\to \R_+$, $k=0,\ldots,n$, be a sequence of non-negative functions with $r$-polynomial growth ($r\ge 1$). Assume that all these functions $F_k$ are {\em convex}, $k=0,\ldots,n$.


\medskip
\noindent $(a)$ {\em Partitioning function}: If, for every $k\!\in \{0,\ldots,n-1\}$,  there exists a convex function $\kappa_k$  such that   $0\le \sigma_k\le \kappa_k \le \theta_k$,  then, for every $x\!\in \R$, 
\[
u_0(x)\le v_0(x).
\]

\noindent $(b)$ {\em Dominating function}:  If the random variable $Z_k$ have symmetric distributions, the functions $\theta_k$, $k=1,\ldots,n$, are convex and  $|\sigma_k|\le \theta_k$, $k=1,\ldots,n$,  then the above inequality remains holds true.
\end {Pro}

\noindent {\bf Remark.} An equivalent formulation of claim $(a)$ is: assume that both $(\sigma_k)_{0\le k\le n}$ and $(\theta_k)_{0\le k\le n}$ are non-negative convex functions with $r$-linear growth, then for every sequence $(\kappa_k)_{0\le k\le n}$ of  functions such that $\sigma_k\le \kappa_k\le \theta_k$, $k=0,\ldots,n-1$, 
\[
u_0(x)\le  c_{\kappa}(x) \le v_0(x)
\]
where $c_{\kappa}(x)$ is the {\em {\em r\'eduite}} of $(F_k(K^x_{0:k}))_{0\le k\le n}$ where $(K^x_k)_{0\le k\le n}$ satisfies the discrete time dynamics
\[
K^x_{k+1}= K^x_k + \kappa_k(K^x_k) Z_{k+1}, \,k=0,\ldots,n-1,\;  K^x_0= x.
\]
This follows from $(a)$ applied successively to the pair $(\sigma_k, \kappa_k)_{0\le k\le n}$ and $(\kappa_k, \theta_k)_{0\le k\le n}$.

\bigskip
\noindent {\sc Proof.} $(a)$ It is clear that this claim is equivalent to proving the expected inequality either if all the functions $(\sigma_k)_{0\le k\le n}$ or all the functions $(\theta_k)_{0\le k\le n}$ are convex. 

We introduce $U^x=(U^x_k)_{0\le k\le n}$ and $V^x=(V^x_k)_{0\le k\le n}$  the $(\P,{\cal F})$-Snell  envelopes of $\big(F_k(X^x_{0:k})\big)_{0\le k\le n}$ and $\big(F_k(Y^x_{0:k})\big)_{0\le k\le n}$ respectively $i.e.$
\[
U^x_k= \P\mbox{-}{\rm supess} \Big\{\E\big(F_{\tau}(X^x_{0:\tau})\,|\, {\cal F}_k\big), \, \tau\; {\cal F}\mbox{-stopping time}, \tau \ge k\Big\}
\]
and
\[
V^x_k= \P\mbox{-}{\rm supess} \Big\{\E\big(F_{\tau}(Y^x_{0:\tau})\,|\, {\cal F}_k\big), \, \tau\; {\cal G}\mbox{-stopping time}, \tau \ge k\Big\}
\]
The connection between {\em {\em r\'eduite}} and Snell envelope is a classical fact  from Optimal Stopping  Theory for which we refer $e.g.$ to~$e.g.$~\cite{NEV}, Chapter~VI), namely
\[
u_0(x) = \E\, U^x_0
\]f
(idem for $v_0$, $V^x_0$ for $Y^x$). It is also classical background on Optimal stopping theory (see again $e.g.$~\cite{NEV}, Chapter~VI)  that  the $(\P,{\cal F})$-Snell envelope $U^x$ satisfies  the following Backward Dynamic Programming principle  
\[
U^x_n=F_n(X^x_{0:n}),\; U^x_k = \max\big(F_k(X^x_{0:k}),\E(U_{k+1}\,|\,{\cal F}_k)\big),\; k=0,\ldots,n-1.
\]
Then, we derive from the dynamics satisfied by the $X^x_k$ and the independence of the random vectors $Z_k$ that $(X^x_k)_{k=0:n}$ is a Markov chain. In turn a  a first backward induction shows that  $U^x_k = u_k(X^x_{0:k})$ $a.s.$, $k=0,\ldots,n$, where the  Borel functions $u_k:\R^{k+1}\to \R$, $k=0,\ldots,n$,  satisfy the backward induction
\begin{equation}\label{eq:BDPP-Amer}
u_n = F_n,\; u_k(x_{0:k})= \max\Big(F_k\big(x_{0:k}), Q_{k+1}u_{k+1}(x_{0:k}, x_k+.)\big)(\sigma_k(x_k))\big)\Big), \; k=0,\ldots, n-1.
\end{equation}
We define likewise the functions $v_k:\R^{k+1}\to \R$, $k=0,\ldots,n$, related to the $(\P,{\cal G})$-Snell  envelopes of $\big(F_k(Y^x_{0:k})\big)_{0\le k\le n}$.

\smallskip To emphasize the analogy with the proof of Proposition~\ref{FCOEuroDiscTime}  we will detail  the case  where all  the functions $\sigma_k= \kappa_k$ are convex, $k=0,\ldots,n$ and satisfy $0\le \sigma_k\le \theta k$.  Following the lines of the proof of this proposition, we show, still by induction, that the  functions $u_k:\R^{k+1}\to \R$ are convex by combining Lemma~\ref{Jen} and~\eqref{eq:BDPP-Amer}. The additional argument to ensure the propagation of convexity is to note that the function $(u,v)\mapsto \max(u,v)$ is convex and increasing in each of its variable $u$ and $v$.

On the other hand, as $0\le \sigma_k\le \theta_k$, $k=0,\ldots,n$ and $\sigma_k$ are all convex, we can show by a new backward induction that $u_k \le v_k$, $k=0,\ldots,n$. If $k=n$ this is obvious. If it holds true with $k+1 \le n$, then for every $x_{0:k}\!\in \R^{k+1}$,
\begin{eqnarray*}u_k(x_{0:k})&\le &\max\Big(\!F_k\big(x_{0:k}\big), \!\big(Q_{k+1}u_{k+1}(x_{0:k},x_k+.)\big)(\theta_k(x_k))\Big)\\
&\le& \max\Big(\!F_k\big((x_{0:k}\big), \!\big(Q_{k+1}v_{k+1}(x_{0:k},x_k+.)\big)(\theta_k(x_k))\!\Big)\!=\! v_k(x_{0:k})
\end{eqnarray*}
where we used successively that  $u\mapsto \big(Q_{k+1}u_{k+1}(x_{0:k},x_k+.)\big)(u)$ is non-decreasing on $\R_+$ since $u_{k+1}$ is convex and that $u_{k+1}\le v_{k+1}$. Finally, the inequality for $k=0$ reads
$$
u_0(x)=\E\, U^x_0 \le\E \,V^x_0 =  v_0(x)
$$ 
which yields the announced result. Other cases   follow the same way round following the lines of the proof of Proposition~\ref{eq:BDPP-Amer}. $\cqfd$

\subsection{Continuous time optimal stopping and American options}
\subsubsection{Brownian diffusions}
In this section,  we switch to the continuous time setting. We will investigate  the (functional) convex order properties of the  {\em {\em r\'eduite}}  (or the  Snell envelope) of payoff processes   obtained as adapted convex functionals of  Brownian martingale diffusion processes $i.e.$ of the form $(F(t, X^t))_{t\in [0,T]}$  where $X^t $ denotes the stopped process $(X_s)_{s\in [0,T]}$ at time $t\!\in [0,T]$ where  $X$ itself is a martingale Brownian diffusion of type $X^{(\sigma)}$  as defined in~\eqref{eq:martdiff}. This embodies most pricing problems for American options in local volatility models. 

In particular, the  results of this section can be seen as an extension to path-dependent ``payoff functionals" of El Karoui-Jeanblanc-Shreve's Theorem (see~\cite{NEKJESH}) which mainly deals with convex functions of the marginal of the processes at time $T$  (see also~\cite{HOB} devoted to parhwise-dependent lookback options).  The proposition below is also  very close to former results by Bergenthum and R\"uschendorf  by combining Theorems~3.2 and~3.6 from~\cite{BRusch1} with Theorem 4.1 from~\cite{BRusch2}. Here,  we focus on the partitioning function. 

\begin{Pro}\label{pro:SnellBrown} Let   $\sigma,\,\theta:[0,T]\times \R$  be two    Lipschitz continuous functions in $(t,x)$
 and let $W$ be a standard ${\cal F}= ({\cal F}_t)_{t\ge 0}$-Brownian motion    defined on a probability space $(\Omega,{\cal A}, \P)$ where  ${\cal F}$ satisfies the usual conditions. Let $(X^{(\sigma),x}_t)_{t\in [0,T]}$  and $(X^{(\theta),x}_t)_{t\in [0,T]}$ be the martingale diffusions, unique strong solutions starting at $x\!\in \R$ to~\eqref{eq:martdiff} (where $W^{(\sigma)}=W^{(\theta)}=W)$.  
 
 Assume  that  there exists  a partitioning function $\kappa :[0,T]\times\R \to \R$ such that $\kappa(t,.)$ is convex for every $t\!\in [0,T]$ with linear growth in $x$ uniformly in $t\!\in [0,T]$ and
$$
0\le \sigma(t,.) \le\kappa(t,.)\le  \theta(t,.), \; t\!\in [0,T].
$$ 
Let  $F:[0,T]\times {\cal C}([0,T],\R)\to \R_+$ be a  $\|\,.\,\|_{\rm \sup}$-continuous functional  with $(r,\|\,.\,\|_{\sup})$-polynomial growth ($r\ge 1$) in $\alpha \!\in {\cal C}([0,T],\R)$, uniformly in $t\!\in [0,T]$.
Moreover, assume that, for every $t\!\in[0,T]$, $F(t,.)$ is {\em convex} on $ {\cal C}([0,T],\R)$. Let $u_0(x)$ and $v_0(x)$ denote the ${\cal F}$-{\em r\'eduites} of $\big(F(t,(X^{(\sigma),x})^t)\big)_{t\in [0,T]}$ and $\big(F(t,(X^{(\theta),x})^t)\big)_{t\in [0,T]}$ respectively defined by 
\[
u_0(x) = \sup\Big\{\E \,F\big(\tau,(X^{(\sigma),x})^{\tau}\big), \; \tau \!\in {\cal T}^{\cal F}_{[0,T]} \Big\} \quad\mbox{ and }\quad v_0(x)= \sup\Big\{\E\,F\big(\tau,(X^{(\theta),x})^{\tau}\big), \; \tau\!\in {\cal T}^{\cal F}_{[0,T]}\Big\}
\]
where ${\cal T}^{\cal F}_{[0,T]}=\{\tau:\Omega\to [0,T], {\cal F}\mbox{-stopping time}\}$. Then
\[
u_0(x)\le v_0(x).
\]
\end{Pro}

\noindent{\bf Remark.} All the quantities involved in the above theorem do exist since the sup norm of $X^{(\sigma),x}$ and $X^{(\theta),x}$ have polynomial moments at any order. Moreover, the Lipschitz continuity assumption is  too stringent but we adopt it to   shorten the proof of  the``approximation" step  from discrete to continuous time dynamics.

\bigskip
\noindent {\bf Proof.}  {\sc Step~1} (Euler schemes) We consider  the Euler schemes   $\bar X^{(\sigma),n}$ and  $\bar X^{(\theta),n}$ (with step $\frac Tn$) of both diffusions (we drop the dependence on the starting value $x$). Both schemes are adapted to the filtration ${\cal F}^{(n)}:=({\cal F}_{t^n_k})_{0\le k\le n}$.

It follows from Proposition~\ref{Snelldiscret} that the $(\P,{\cal F}^{(n)})$-Snell envelopes $\bar U^{(n)}=(\bar U^{(n)}_{t^n_k})_{0\le k\le n}$,  $\bar K^{(n)}=(\bar K^{(n)}_{t^n_k})_{0\le k\le n}$ and $\bar V^{(n)}=(\bar V^{(n)}_{t^n_k})_{0\le k\le n}$ of  the ${\cal F}^{(n)}$-adapted payoff processes $F\Big(t^n_k,\big[ I_n\big(\bar X^{(\sigma),n}\big)\big]^{t^n_k}\Big)$, $k=0,\ldots,n$,   $F\Big(t^n_k,\big[ I_n\big(\bar X^{(\kappa),n}\big)\big]^{t^n_k}\Big)$, $k=0,\ldots,n$,  and $F\Big(t^n_k,\big[ I_n\big(\bar X^{(\theta),n}\big)\big]^{t^n_k}\Big)$, $k=0,\ldots,n$,  satisfy
\begin{equation}\label{eq:  ineqreduites-n}
\E \,\bar U^n_0 \le \E \bar K^n_0 \le  \E\, \bar V^n_0.
\end{equation}
Note that it is always possible to define the Euler scheme associated to the function $\kappa$ regardless of its convergence toward the related SDE.

\noindent  {\sc Step~2} (Convergence) First,  set for convenience $\bar Y^{(n)}_{t^n_k}= F\Big(t^n_k,\big[ I_n\big(\bar X^{(\sigma),n}\big)\big]^{t^n_k}\Big)$, $k=0,\ldots,n$, so that
\[
\bar U^{(n)}_{t^n_k}= \P\mbox{-} {\rm supess}\Big\{\E( \bar Y^{(n)}_{\tau}\,|\, {\cal F}_{t^n_k}),\; \tau\!\in {\cal T}^{(n)}_{t^n_k,T}\Big\},\;k=0,\ldots, n,
\]
where $ {\cal T}^{(n)}_{t^n_k,T}=\Big\{ \tau :\Omega\to \{t^n_k, \ldots,t^n_{\ell}, \ldots, t^n_n\}, \;{\cal F}^{(n)}\mbox{-stopping time}\Big\} $; we also know that  the $(\P, {\cal F})$-Snell envelope of the process $Y_t=F(t, X^t)$, $t\!\in [0,T]$, is defined by 
\[
U_t= \P\mbox{-}{\rm supess} \Big\{\E\big(Y_{\tau}\,|\, {\cal F}_{t}\big), \, \tau\!\in {\cal T}^{\cal F}_{t,T}   \Big\}, \;t\in [0,T],
\]
where  ${\cal T}^{\cal F}_{t,T}=  \Big\{\tau :\Omega\to[t,T], \; {\cal F}\mbox{-stopping time}\Big\}$. This Snell envelope is well-defined since $\|X\|_{\sup}$ lies in every $L^p(\P)$, $p\!\in (0, +\infty)$, which implies in turn that $\|Y\|_{\sup}$ lies in every $L^p(\P)$.   As the obstacle process $(F(t, X^t))_{t\in [0,T]}$ has continuous paths and is uniformly integrable, it is regular for optimal stopping and  $t\mapsto \E\, U_t$ is continuous  (see~\cite{NEK,LAMB}. Hence, the super-martingale  $(U_t)_{t\in [0,T]}$ has a (non-negative) c\`adl\`ag modification     whose compensator  is continuous (and non-decreasing). More generally, if a sequence os stopping times $\tau_n\uparrow \tau<+\infty$ and $U_{\tau}\!\in L^1$, then $\E U_{\tau_n} \to \E U_{\tau}$.
For technical purpose, we introduce an   intermediate  quantity defined by
\[
\widetilde U_{t^n_k} =\P\mbox{-} {\rm supess}\Big\{\E(Y_{\tau}\,|\, {\cal F}_{t^n_k}),\; \tau\!\in {\cal T}^{(n)}_{t^n_k,T} \Big\} \le U_{t^n_k},\;k=0,\ldots, n.
\]
 Our aim is to prove,  after having canonically extended  $\bar U^{(n)}$  into a c\`adl\`ag stepwise constant process by setting $\bar U^{(n)}_t=\bar U^{(n)}_{t^n_k}$, $t\in [t^n_k,t^n_{k+1})$, that $\bar U^{(n)}_t$ converges to $U_t$  in $L^p$ for every $t\!\in [0,T]$. 
  We start from the fact that 
  \begin{equation}\label{eq:decamp}
  |U_t-\bar U^{(n)}_{\underline t_n} |\le | U_t- U_{\underline t_n}|+ U_{\underline t_n} - \widetilde U^{(n)}_{\underline t_n}+| \widetilde U^{(n)}_{\underline t_n} - \bar U^{(n)}_{\underline t_n}|.
  \end{equation}
  
 Once again,  regularity for optimal stopping  of $U$ implies in particular that $(U_t)_{t\in [0,T]}$  is $L^1$-left continuous in $t$. In particular $\E|U_t-U_{t^n_k}| \to 0$ as $n\to +\infty$.

\smallskip  
As concerns the second term in the right hand side of~\eqref{eq:decamp}, we proceed as follows
\[
0\le U_{t^n_k} - \widetilde U^{(n)}_{t^n_k} \le  \P\mbox{-}{\rm supess} \Big\{\E \big(Y_{\tau}-Y_{\tau^{(n)}}\,|\, {\cal F}^{(n)}_{t^n_k}\big),\; \tau \!\in {\cal T}_{t^n_k,T}\Big\}
\]
where $\tau^{(n)}=\sum_{\ell= k }^n \frac {\ell T}{n} \mbox{\bf 1}_{\{\frac{(\ell -1)T}{n}<\tau\le \frac{\ell T}{n}\}} = \sum_{\ell= k }^n \bar t^n \mbox{\bf 1}_{\{t^n_{\ell-1}<\tau\le t^n_{\ell}\}}\!\in  {\cal T}^{(n)}_{t^n_k,T}\subset  {\cal T}_{t^n_k,T}$ so that  
\[
0\le U_{t^n_k} - \widetilde U^{(n)}_{t^n_k} \le \E\big(\sup_{t\ge t^n_k}|Y_{t}-Y_{\overline t^n}|\,|\, {\cal F}_{t^n_k} \big)\le \E\big(\sup_{t\in [0,T]}|Y_{t}-Y_{\overline t^n}|\,|\, {\cal F}_{t^n_k} \big).
\]
Doob's Inequality applied to the martingale $M_n = \E\big(\sup_{t\in [0,T]}|Y_{t}-Y_{\overline t^n}|\,|\, {\cal F}_{t^n_k} \big)$, $n\ge 1$,   implies that  for every $p\!\in (1,+\infty)$, 
\begin{eqnarray*}
\Big \|\max_{0\le k\le n} (U_{t^n_kn} - \widetilde U^{(n)}_{t^n_k}) \Big \|_p\le \frac{p}{p-1}\|M_n\|_{_p} =     \frac{p}{p-1} \Big\| \sup_{t\in [0,T]} |Y_t - Y_{\overline t_n}|\Big\|_p&\to& 0 \;\mbox{ as } n\to +  \infty
\end{eqnarray*}
since  $X^{\overline t_n}$ $a.s.$ converges towards $X^t$  for the sup-norm owing to the pathwise continuity of $X$. This in turn, implies that $F(\overline t^n, X^{\overline t^n})$  $a.s.$ converges toward $F( t, X^t)$ since $F$ is continuous. The $L^p$-convergence follows by uniform integrability, still since $\|Y\|_{\sup}$ has polynomial moments at any order. 

Now we investigate the second term in the right hand side of~\eqref{eq:decamp}. 
\begin{eqnarray*}
| \widetilde U^{(n)}_{t^n_k}- \bar U^{(n)}_{t^n_k}  |&\le &  {\rm supess}\Big\{\E\big(|Y_{\tau}  -\bar  Y^{(n)}_{\tau}   |\,|\, {\cal F}_{t^n_k}\big), \; \tau \!\in {\cal T}^{(n)}_{t^n_k,T}\Big\}\\
&\le& \E\Big( \max_{0\le k\le n} |\bar Y^{(n)}_{t^n_k}-Y_{t^n_k} | \,|\, {\cal F}_{t^n_k}  \Big).
\end{eqnarray*}

On the other hand, 
\begin{eqnarray}
\nonumber  \max_{0\le k\le n} |\bar Y^{(n)}_{t^n_k}-Y_{t^n_k}|&\le &  \max_{0\le k\le n} \big|F\big(t^n_k, (I_n(\bar X^{(\sigma),n}))^{t^n_k}\big)- F\big(t^n_k, (X^{(\sigma)})^{t^n_k}\big)\big|\\
\label{eq:Snell-major}&\le &  \sup_{t\in [0,T]}\big|F\big(t, (I_n(\bar X^{(\sigma),n}))^{t}\big)- F\big(t, (X^{(\sigma)})^{t}\big)\big|.
\end{eqnarray}
Now, note that the functional  $\displaystyle \alpha \mapsto \Big( t \mapsto  F\big (t, \alpha^t\big)\Big)
$ defined from $({\cal C}([0,T], \R), \|\,.\,\|_{\sup})$ into itself is continuous: if $(t_n,\alpha_n)\to (t,\alpha)$ for the product topology on  $[0,T]\times ({\cal C}([0,T], \R)$, then 
\[
\|\alpha_n^{t_n}-\alpha^t \|_{\sup} \le \|\alpha_n-\alpha\|_{\sup} + w(\alpha, |t-t_n|)
\]
so that $ (t_n, \alpha^{t_n})\to (t,\alpha^t)$. As a consequence,  the functional $F$ being continuous  on $[0,T]\times {\cal C}([0,T],\R)$, $F(t_n,\alpha^{t_n})\to F(t,\alpha^t)$ which in turn implies that 
$\sup_{t\in [0,T]}|F(t,\alpha_n^{t})-  F(t,\alpha^{t})|\to 0$. As  $I_n( \alpha) \to \alpha$ for the sup norm as $n\to +\infty$, we derive that if $\alpha_n \to \alpha$ for the sup norm then $\displaystyle \sup_{t\in [0,T]} |F(t,I_n(\alpha_n)^{t})-F(t,\alpha^{t})|\to 0$ as $n\to +\infty$. 

Then under the Lipschitz continuity assumption on $\sigma$, we know that the Euler scheme $\bar X^{(\sigma),x,n})\to X^{(\sigma),x}$ $\P$-$a.s.$ as $n\to +\infty$ $a.s.$ (see $e.g.$~\cite{BOLE}, Theorem B.14, p.276).  The $(r, \|\,.\,\|_{\sup})$-polynomial growth  assumption made  on $F$ and the the fact that $\displaystyle \sup_{n\ge 1}\E \|\bar X^{(\sigma),x,n}\|^{p}_{\sup}<+\infty$ for any $p>r$ implies the $L^1$-convergence to $0$ of the term in~\eqref{eq:Snell-major}. Finally, this shows that 
\[
\E\, \bar U^n_0\to u_0(x) \quad \mbox{ as }\quad n\to +\infty. 
\]
%
The conclusion   follows from~\eqref{eq:  ineqreduites-n} in Step~1 by letting $n\to +\infty$ in the resulting inequality $\E\, \bar U^n_0\le \E\,\bar V^n_0$.
$\cqfd$

\bigskip 
\noindent {\sc Applications to comparison theorems for American options in  local volatility models}. By specifying our diffusion dynamics as a local volatility model as defined by~(\ref{LocalVol}), we can extend the comparison result~(\ref{ComparLocVol}) to path-dependent American options provided the ``payoff" functionals $F(t,.)$ are convex with polynomial growth as specified in the above theorem. 

\subsubsection{The case of jump martingale diffusions}

In what follows the product space  $[0,T] \times \D([0,T],\R)$ is endowed with the  product topology $ |\,.\,| \otimes Sk$.  The  notation $X_t(\alpha)= \alpha(t)$, $\alpha\!\in \D([0,T],\R)$ still denotes the canonical process  on $\D([0,T],\R)$  and $\theta$   denotes the canonical random variable on $[0,T]$ ($i.e.$ $\theta(t)=t$, $t\!\in [0,T]$).  

\medskip
 Let $({\cal F}_t)_{t\in [0,T]}$ be a right continuous  filtration on a probability space $(\Omega,{\cal A}, \P)$ and let $Y$ be an $({\cal F}_t)_{t\in [0,T]}$-adapted  c\`adl\`ag process defined on this probability space. We introduce the so-called (${\cal H}$)-assumption (also known as filtration enlargement assumption) which reads as follows:
\[
({\cal H})\equiv  \forall\, H:\Omega\to \R,  \mbox{ bounded and } {\cal F}^Y_{_T}\mbox{-measurable}, \; \E \big(H\,|\,{\cal F}_t\big)= \E \big(H\,|\,{\cal F}^Y_t\big)\; \P\mbox{-}a.s.
\]
This filtration enlargement assumption is equivalent to the following more tractable condition: there exists $D\subset [0,T]$, everywhere dense in $[0,T]$, with $T\!\in D$, such that 
\[
\forall\; n\ge1,\; \forall\, t_1,\ldots,t_n\!\in D,\; \forall\, h\!\in {\cal C}_0(\R^n,\R), \; \E \big(h(Y_{t_1},\ldots, Y_{t_n})\,|\,{\cal F}_{t}\big) = \E \big(h(Y_{t_1},\ldots, Y_{t_n})\,|\,{\cal F}^Y_{t}\big) \; \P\mbox{-}a.s. 
\]
where ${\cal C}_0(\R^n,\R)= \{f\!\in {\cal C}(\R^n,\R)\, \mbox{such that}\, \lim_{|x|\to+\infty} f(x)=0\}$.  We still consider the jump diffusions of the form~(\ref{LevyEDS}) $i.e.$
\[
dX_t =\kappa(t,X_{t_-})dZ_t
\]
where $\kappa:[0,T]\times \R\to \R$ is a continuous function, Lipschitz continuous in $x$ uniformly in $t\!\in[0,T]$.

The aim of this section is to extend the result obtained for convex order for Brownian diffusions to such jump diffusions. We will rely on an    abstract convergence result  for  {\em r\'eduites} established in~\cite{LAPA} (Theorem~3.7 and the remark that follows) that we recall below. To this end, we need to recall two classical definitions on stochastic processes.

\begin{Dfn}\label{def:DetAldous} $(a)$ {\em Class $(D)$ processes:} A c\`adl\`ag process $(Y_t)_{t\in [0,T]}$ is of class $(D)$ if
\begin{equation}\label{ClasseD}
\left\{Y_{\tau},\; \tau \!\in {\cal T}_{[0,T]} \right\}\mbox{ is uniformly integrable.}
\end{equation}

\noindent $(b)$ {\em Aldous's tightness criterion  (see $e.g.$~\cite{JASH2},  Chapter VI,  Theorem 4.5, p.356):} A sequence of  ${\cal F}^n$-adapted  
c\`adl\`ag processes $Y^n=(Y^n_t)_{t\in [0,T]}$, $n\ge 1$,  defined on filtered stochastic spaces $(\Omega^n, {\cal A}^n, {\cal F}^n,\P^n)$, $n\ge 1$, satisfies Aldous's tightness  criterion with respect to the filtrations ${\cal F}^n$, $n\ge1$,   if 
\begin{equation}\label{Aldous}
\forall\, \eta >0,\; \lim_{\delta \to 0}\limsup_n \sup_{\tau_n\le \tau_n'\le (\tau_n +\delta)\wedge T} \P^n\big(|Y^n_{\tau_n}-Y^n_{\tau_n'}|\ge \eta\big) =0
\end{equation}
where $\tau_n$and $\tau_n'$ run over $[0,T]$-valued ${\cal F}^{Y^n}$-stopping times. 

Then, the sequence $(Y^n)_{n\ge 1}$ is tight for the Skorokhod topology.
\end{Dfn}
\begin{Thm}\label{ContRed} $(a)$ Let $(X^n)_{n\ge 1}$ be a sequence of adapted quasi-left   c\`adl\`ag processes~(\footnote{A c\`adl\`ag $({\cal F}_t)_{t\in[0,T]}$-adapted process $X=(X_t)_{t\in [0,T]}$ is {\em quasi-left continuous} with respect to the right continuous  filtration ${\cal F}=({\cal F}_t)_{t\in[0,T]}$ if for every   ${\cal F}$-stopping time $\tau$ having values in $[0,T]\cup\{+\infty\}$ and every increasing sequence of ${\cal F}$-stopping times $(\tau_k)_{k\ge 1}$  with limit $\tau$, $\lim_k X_{\tau_k} = X_{\tau}$ on  the event $\{\tau<+\infty\}$ (see $e.g.$~\cite{JASH2}, Chapter ~I.2.25, p.22).}) defined on a probability spaces $(\Omega^n,{\cal F}^n,\P^n)$ of class $(D)$ and satisfying the above Aldous tightness criterion~\eqref{Aldous}. For every $n\ge 1$, let
\[
u^n_0= \sup\left\{\E\, X^n_{\tau},\; \tau \; [0,T]\mbox{-}valued\, {\cal F}^n\mbox{-stopping time}\right\}
\]
denote the ${\cal F}^n$-{\em r\'eduite} of $X^n$. Let $(\tau_n^*)_{n\ge 1}$ be a sequence of $ \big({\cal F}^{X^n}, \P^n)$-optimal stopping times(\footnote{$i.e.$ satisfying $\E\, X^n_{\tau_n^*}=u^n_0$.}). Assume furthermore that $(X^n)_{n\ge 1}$ satisfies
\[
X^n \stackrel{{\cal L}}{\longrightarrow} \P, \; \P\mbox{ probability measure on } (I\!\!D([0,T],\R), {\cal D}_T)\; \mbox{such that } \;
\E_{\P}\sup_{t\in [0,T]}|X_t|<+\infty.
\]

If every limiting value $\Q$ of ${\cal L}(X^n,\tau^*_n)$ on  $I\!\!D([0,T],\R)\times [0,T]$   satisfies the (${\cal H}$) property, then the $({\cal F}^n, \P^n)$-{\em r\'eduites} $u^n_0$ of $X^n$ converge toward the $({\cal D}, \P)$-{\em r\'eduites} $u_0$ of $X$ $i.e.$
\[
\lim_n u^n_0= u_0.
\]
Moreover, if the optimal stopping problem related to $(X,\Q, {\cal D}^{\theta})$ has a unique solution in distribution, $i.e.$ $\theta\stackrel{d}{=}\mu^*_{\tau^*}$, not depending on $\Q$, then $\displaystyle \tau^*_n\stackrel{{\cal L}([0,T])}{\longrightarrow} \mu^*_{\tau^*}$.

\noindent$(b)$ The same result holds when considering  a sequence of {\em companion processes} $Y^n$ having values in a Polish metric space $(E,d_{_E})$  $i.e.$ we consider that the filtration of interest at finite range $n$ is now $({\cal F}^{(X^n,Y^n)}_t)_{t\in [0,T]}$. We assume  that $X^n$ is quasi-left continuous with respect to this enlarged filtration. We will only ask  the couple $(X^n,Y^n)$ to converge for the {\em product topology} $i.e.$ on $(\D([0,T], \R), Sk_{\R})\times (\D([0,T], E), Sk_E)$ since this product topology spans the same Borel $\sigma$-field as the regular  Skorokhod topology on $\D([0,T], \R\times E)$.
\end{Thm}

The main result of this section is the following:
\begin{Thm} 
Let $Z=(Z_t)_{t\in [0,T]}$ be a martingale L\'evy process with L\'evy measure $\nu$ satisfying $\nu(|z|^p
)<+\infty$ for  $p\!\in [ 2,+\infty)$, so that the process $Z$ is an $L^2$-martingale null at $0$. Let $X^{(\kappa_i,x)}$, $i=1,2$,  be the martingale jump diffusions driven by $Z$ starting at (the same) $x\!\in \R$. Let $F:[0,T]\times \D([0,T], \R)\to \R_+$ be a {\em convex} functional satisfying  the following local Lipschitz  assumption (w.r.t. to the sup norm) combined with a Skorokhod continuity assumption, namely
\begin{equation}\label{HypoLipsup}
\left\{\begin{array}{ll}
(i) & F:[0,T]\times \D([0,T], \R)\to \R_+  \mbox{ is $Sk$-continuous},\\
(ii) &  |F(t,\beta)- F(s,\alpha)| \le C\Big(|t-s|^{\rho'} + \|\alpha-\beta\|^{\rho}_{\sup}\big(1+\|\alpha\|_{\sup}^{r-\rho}+\|\beta\|_{\sup}^{r-\rho} \big)\Big),\; \rho,\, \rho'\!\in (0,1],\; r\!\in[1,p), 
\end{array}\right.
\end{equation}
 Let $U^{(\kappa_i)}$ denote  the Snell envelopes of the processes $\big(F(t,(X^{\kappa_i}))^t\big)_{t\in [0,T]}$, $i=1,2$ respectively.

If there exist $\kappa_i:[0,T]\times \R\to \R$, $i=1,2$, two continuous functions
with linear growth in $x$, uniformly in $t\!\in [0,T]$, and a {\em partitioning 
 function} $\kappa:[0,T]\times \R\to \R$,  convex   in $x$ 
 for every $t\!\in [0,T]$, such that
\[
\kappa_1\le \kappa \le \kappa_2.
\]
 Then
 \[
 U_0^{(\kappa_1)}\le U_0^{(\kappa_2)}.
 \]
\end{Thm}

\noindent {\bf Remarks.} $\bullet$ Note that, since $p\ge 2$,
\[
\nu(|z|^p)<+\infty \longleftrightarrow \nu(|z|^p\mbox{\bf 1}_{\{|z|\ge 1\}})<+\infty  \longleftrightarrow Z_t \!\in L^p  \longleftrightarrow\sup_{t\in [0,T]} |Z_t| \!\in L^p.
\]

\smallskip \noindent $\bullet$ One proves likewise that, for every $t\!\in [0,T]$, 
\[
 \E (U_t^{(\kappa_1)})\le \E(U_t^{(\kappa_2)}).
\]

\noindent $\bullet$ If the functions $\kappa(t,.)$, $t\!\in [0,T]$  are all convex (but possibly not the functions $\kappa_i(t,.)$) then the same proof shows by coupling $(\kappa_1,\kappa)$ and $(\kappa, \kappa_2)$ that 
\[
 \E (U_0^{(\kappa_1)})\le  \E(U_0^{(\kappa)}) \le \E(U_0^{(\kappa_2)}).
\]

\begin{Lem}\label{Lem:7.1} Let $X=(X_t)_{t\in [0,T]}$ be an $({\cal F}_t)_{t\in[0,T]}$-adapted c\`adla\`ag process defined on a probability space $(\Omega,{\cal A}, \P)$ where $({\cal F}_t)_{t\in[0,T]}$ is a c\`ad filtration. Let $G:[0,T]\times I\!\!D([0,T],\R)\to \R_+$ be a Skorokhod continuous functional such that $|G(\alpha)|\le C(1+\|\alpha\|^r_{\sup})$, $r\!\in (0,p)$. If $X$ is quasi-left continuous
and if $\|X\|_{\sup} \!\in L^p$, then  the ``obstacle process" $(G(t,X^t))_{t\in [0,T]}$ is {\em regular for optimal stopping}  $i.e.$ as soon as $\tau <+\infty$ $\P$-$a.s.$ $\E G(\tau_n,X^{\tau_n})\to \E\, G(\tau,X^{\tau})$.
\end{Lem}

\noindent {\bf Proof.} First one easily proves by coming back to the very definition of Skorokhod topology that   $\alpha_n\stackrel{Sk}{\longrightarrow} \alpha$ and $t_n\to t\!\in {\rm Cont}(\alpha)$ then 
$ \alpha_n^{t_n} \stackrel{Sk}{\longrightarrow}\alpha^t$.
Let  $(\tau_n)_{n\ge 1}$ be a sequence of ${\cal F}_t$-stopping times satisfying  $\tau_n\uparrow \tau<+\infty$ $\P$-$a.s.$, then 
$X_{\tau}=X_{\tau_-}$ $\P$-$a.s.$ $i.e.$ $\tau(\omega) \!\in {\rm Cont}(X(\omega))$ $\P(d\omega)$-$a.s.$. It follows that $(\tau_n X^{\tau_n})\to (\tau, X^{\tau})$ $\P$-$a.s.$.  The continuity assumption made on  $G$ implies that $G(\tau_n,X^{\tau_n})\stackrel{Sk}{\longrightarrow} G(\tau,X^{\tau})$. One concludes by a uniform integrability argument that $\E\, G(\tau_n,X^{\tau_n})\to \E\, G(\tau,X^{\tau})$ since  $\|X\|_{\sup}\!\in L^p$ implies that 
$\big(G(\tau_n,X^{\tau_n})\big)_{n\ge 1} $ is $L^{\frac pr}$-bounded. $\cqfd$

\bigskip
\noindent {\bf Proof.} {\sc Step~1} {\em Aldous tightness criterion.} We still consider the stepwise constant Euler scheme  $\bar X^n= (\bar X^n_t)_{t\in [0,T]}$  with step $\frac T n$defined by
\[
\bar X^n_{t^n_k}= \bar X^n_{t^n_{k-1}} +\kappa(t^n_{k-1}, \bar X^n_{t^n_{k-1}})(Z_{t^n_k}-Z_{t^n_{k-1}}), \; k=1,\ldots,n,\; \bar X^n_0=X_0
\] 
and  $i.e.$ $\bar X^n_t= \bar X^n_{\underline t_n}$. Let $\sigma_n$, $\tau_n\!\in {\cal T}^{{\cal F}^n}_{[0,T]}$, 
such that $\sigma_n\le \tau_n\le (\sigma_n+\delta)\wedge T$. In fact, following Lemma~\ref{Lem:6.3}, we may assume without loss of generality  that $\sigma_n$ and $\tau_n$ take values in $\{t^n_k,\, k=0,\ldots,n\}$. Then, owing to~\eqref{HypoLipsup},
\begin{eqnarray*}
\E\big|F(\tau_n, (\bar X^n)^{\tau_n})- F(\sigma_n, (\bar X^n)^{\sigma_n})\big|\le C\delta{^{\rho'}}+C\, \E\big(\| (\bar X^n)^{\tau_n}- (\bar X^n)^{\sigma_n}\|^{\rho}_{\sup} (1+2\|\bar X^n\|^{r-\rho}_{\sup})\big).
\end{eqnarray*}
H\"older Inequality  applied with the conjugate exponents $a=\frac{r}{\rho}$ and $b= \frac{r}{r-\rho}$ yields
\begin{eqnarray*}
\E\Big(\| (\bar X^n)^{\tau_n}- (\bar X^n)^{\sigma_n}\|^{\rho}_{\sup} \big(1+2\|\bar X^n\|^{r-\rho}_{\sup}\big)\Big)&\le& \Big\| \sup_{\sigma_n \le s\le (\sigma_n+\delta)\wedge T} |\bar X^n_s-\bar X^n_{\sigma_n}|      \Big\|^{\rho}_{r}\Big(1+2\big\|\sup_{t\in [0,T]}|\bar X^n_t| \big\|^{r-\rho}_{r}\Big).
\end{eqnarray*}

As $\nu(z^2)<+\infty$, we can decompose the L\'evy process $Z$ into $Z_t= a\,W_t+ \widetilde Z_t$, $a\ge 0$ 
where $W$ is a standard Brownian motion and $\widetilde Z$ is a pure jump square integrable martingale L\'evy process. 

\ms \ni $\bullet$ If $r\!\in [1,2]$: it follows from the $B.D.G.$  Inequality applied to the local martingale $(\bar X_{\sigma_n+\frac{iT}{n}}-\bar X^n_{\sigma_n})_{i\ge 0}$ that 
\begin{eqnarray*}
\Big\| \sup_{\sigma_n \le t^n_k\le (\sigma_n+\delta)\wedge T} |\bar X^n_{t^n_k}-\bar X^n_{\sigma_n}|\Big \|^r_{r}\!&\!\le \!&\!c_r a^r \Big\| \sum_{\sigma_n <t^n_k \le (\sigma_n+\delta)\wedge T} \hskip -0.75cm \kappa(t^n_{k-1},\bar X^n_{t^n_{k-1}})^2\big(W_{t^n_k} -  W_{t^n_{k-1}} \big)^2\Big\|_{L^{\frac r2}}^{\frac r2}\\\
\!&\!\!&\!+c_r \Big\| \sum_{\sigma_n <t^n_k\le (\sigma_n+\delta)\wedge T} \hskip -0.75cm \kappa(t^n_{k-1},\bar X^n_{t^n_{k-1}})^2\big( \widetilde Z_{t^n_k} -  \widetilde  Z_{t^n_{k-1}} \big)^2\Big\|_{L^{\frac r2}}^{\frac r2}\\
\!&\!\le\!&\! c_r a^r\Big\| \sum_{\sigma_n <t^n_k\le (\sigma_n+\delta)\wedge T} \hskip -0.75cm \kappa(t^n_{k-1},\bar X^n_{t^n_{k-1}})^2\big(W_{t^n_k} -  W_{t^n_{k-1}} \big)^2\Big\|_{L^{2}}^{\frac r2}\\
\!&\!\!&\! + c_r\E \left(\sum_{k} \mbox{\bf 1}_{\{\sigma_n <t^n_k\le (\sigma_n+\delta)\wedge T\}}|\kappa(t^n_{k-1}, \bar X^n_{t^n_{k-1}})|^r|Z_{t^n_k}-Z_{t^n_{k-1}}|^r\right).
\end{eqnarray*}
Now
\begin{eqnarray*}
\!\E \Big[ \!\sum_{\sigma_n <t^n_k\le (\sigma_n+\delta)\wedge T} \hskip -0.75cm \kappa(t^n_{k-1},\bar X^n_{t^n_{k-1}})^2\big(W_{t^n_k} \!-\!  W_{t^n_{k-1}} \big)^2\Big]\!&\!=\!&\! \frac Tn \E \Big[ \sum_{\sigma_n <t^n_k\le (\sigma_n+\delta)\wedge T} \hskip -0.75cm \kappa(t^n_{k-1},\bar X^n_{t^n_{k-1}})^2\Big]\\
\!&\!\le \!&\! \frac Tn \E |Z_{\frac Tn}|^r \E\Big[\! \max_{1\le k\le n}|\kappa(t^n_{k-1},\bar X^n_{t^n_{k-1}})|^2\!\!\times\!\mbox{\rm card}\{k\!:\!\sigma_n\! <\!t^n_k\!\le\! (\sigma_n +\delta)\!\wedge\! T\}\!\Big]\\
\!&\!\le \!&\! \frac Tn\E\Big[ \max_{1\le k\le n}|\kappa(t^n_{k-1},\bar X^n_{t^n_{k-1}})|^r\Big] \frac{\delta n}{T}\\
\!&\!=\!&\! \delta   \Big\| \max_{1\le k\le n}|\kappa(t^n_{k-1},\bar X^n_{t^n_{k-1}})|\Big\|_{2}^2
\end{eqnarray*}
 On the other hand,
\begin{eqnarray*} 
\Big\| \sum_{\sigma_n <t^n_k\le (\sigma_n+\delta)\wedge T} \hskip -0.75cm \kappa(t^n_{k-1},\bar X^n_{t^n_{k-1}})^2\big( \widetilde Z_{t^n_k} -  \widetilde  Z_{t^n_{k-1}} \big)^2\Big\|_{L^{\frac r2}}^{\frac r2}\!&\!\le \!&\!\E \sum_{k} \mbox{\bf 1}_{\{\sigma_n <t^n_k\le (\sigma_n+\delta)\wedge T\}}|\kappa(t^n_{k-1}, \bar X^n_{t^n_{k-1}})|^r|\widetilde Z_{t^n_k}-\widetilde Z_{t^n_{k-1}}|^r\\
\!&\!= \!&\! \E |\widetilde Z_{t^n_k}-\widetilde Z_{t^n_{k-1}}|^r \E  \Big[\sum_{k} \mbox{\bf 1}_{\{\sigma_n <t^n_k\le (\sigma_n+\delta)\wedge T\}}|\kappa(t^n_{k-1}, \bar X^n_{t^n_{k-1}})|^r\Big]\\
\!&\!=\!&\! \E |\widetilde Z_{\frac Tn}|^r \E\Big[ \!\max_{1\le k\le n}|\kappa(t^n_{k-1},\bar X^n_{t^n_{k-1}})|^r\!\!\times\!\mbox{\rm card}\{k\!:\!\sigma_n\! <\!t^n_k\!\le\! (\sigma_n +\delta)\!\wedge\! T\}\!\Big]\\
 \!&\!\le \!&\! \E |\widetilde Z_{\frac Tn}|^r \E\Big[ \max_{1\le k\le n}|\kappa(t^n_{k-1},\bar X^n_{t^n_{k-1}})|^r\Big] \frac{\delta n}{T}\\
 \!&\!\le\!&\! \delta \Big(\frac{n}{T} \E |\widetilde Z_{\frac Tn}|^r \Big) \E\Big[ \max_{0\le k\le n-1}|\kappa(t^n_{k-1},\bar X^n_{t^n_{k-1}})|^r\Big]\\
 \!&\!\le\!&\! C_{\kappa, \widetilde Z,T} \,\delta   \Big\| \max_{1\le k\le n}|\kappa(t^n_{k-1},\bar X^n_{t^n_{k-1}})|\Big\|_{r}^r
\end{eqnarray*}
where we used that $t\mapsto \frac 1t \E |\widetilde Z_t|^r$ remains bounded on the whole interval $(0,T]$.

Under the assumptions $\nu(z^2)<+\infty$ and $\kappa$ with linear growth (in $x$ uniformly in $t\!\in [0,T]$),  it follows form Proposition~\ref{pro:highmoment} in Appendix~\ref{app:B}   that  $\sup_{n\ge 1} \big\|\sup_{0\le k\le n}|\kappa(t^n_{k}, \bar X^n_{t^n_k})| \big\|_r <+\infty $ since $r\le 2$ (see the first remark below the statement of the theorem),
we get 
\[
\Big\| \sup_{\sigma_n \le t^n_k\le (\sigma_n+\delta)\wedge T} |\bar X^n_{t^n_k}-\bar X^n_{\sigma_n}|\Big \|_{r} \le C_{\rho,r,\kappa,Z,T} \big( \delta^{\frac 14}+ \delta^{\frac 1 r} \big)
\]
where the real constant $C_{\rho,r,\kappa,Z,T} $ does not depend on $n$, $\sigma_n$, $\tau_n$ and $\delta$. This implies in turn that
\[
\lim_{\delta\to 0}\limsup_n \sup_{\sigma_n <\tau_n\le (\sigma_n+\delta)\wedge T} \E \big| F(\tau_n,\bar X^{n, \tau_n})- F(\sigma_n,\bar X^{n,\sigma_n})\big|=0
\]
and the conclusion follows.
 
\ms \ni $\bullet$ If $r\!\in [2,4]$: One writes 
\begin{eqnarray*}
 \sum_{\sigma_n <t^n_k\le (\sigma_n+\delta)\wedge T} \hskip -0.75cm \kappa(t^n_{k-1},\bar X^n_{t^n_{k-1}})^2\big( Z_{t^n_k}-Z_{t^n_{k-1}} \big)^2&=&  \sum_{\sigma_n <t^n_k\le (\sigma_n+\delta)\wedge T} \hskip -0.75cm \kappa(t^n_{k-1},\bar X^n_{t^n_{k-1}})^2(\big( Z_{t^n_k}-Z_{t^n_{k-1}} \big)^2-\E |Z_{\frac Tn}|^2)\\
 && + \E |Z_{\frac Tn}|^2 \sum_{\sigma_n <t^n_k\le (\sigma_n+\delta)\wedge T} \hskip -0.75cm \kappa(t^n_{k-1},\bar X^n_{t^n_{k-1}})^2\big( Z_{t^n_k}-Z_{t^n_{k-1}} \big)^2.
\end{eqnarray*}
The second term of the sum in the right hand side of the above equality can be treated as above (it corresponds to $r=2$). As concerns the first one, note that  the $i.i.d.$ sequence  $\big(( Z_{t^n_k}-Z_{t^n_{k-1}})^2-\E |Z_{\frac Tn}|^2 \big)_{1\le k\le n} $ is centered and lies in $L^{\frac r2}(\P)$ with $\frac r2 \!\in [1,2]$. Hence it can be controlled like the former  case. 
Carrying on the process by a cascade induction as detailed $e.g.$ in the proof of Proposition~\ref{pro:highmoment} in Appendix~\ref{app:B}, one can  lower $r$ to $r/2,\dots, r/2^{\ell_r}\!\in (1,2]$  by induction, owing to $B.D.G.$  inequality.


%

\ms
\noindent {\sc Step~2}. It follows from Step~1 of Theorem~\ref{FCOEuroJump} (adapted to a $2$-dimensional  framework with $(\kappa,\mbox{\bf 1})$ as a drift) that 
\[
\Big(\bar X^n, I_n(Z)\Big)\stackrel{\L(Sk)}{\longrightarrow} (X,Z)\quad \mbox{ as }\quad n\to +\infty.
\]
 If we consider the discrete time Optimal Stopping problem(s) related to the Euler schemes $\bar X^{(n,\kappa_i)}$, $i=1,2$,  which turns out the be the same as in Step~1 of the proof of Proposition~\ref{pro:SnellBrown}, the existence of  optimal stopping times $\tau^{(i)}_n$, $i=1,2$,  taking values in $\{t^n_k,\, k=0,\ldots,n\}$ is straightforward  owing to the finite horizon of these problems (see~\cite{NEV}, Chapter~VI for more details).

\medskip
\noindent {\sc Step~3:} Let $\Omega_c = I\!\!D([0,T], \R)^2\times [0,T]$ be the canonical space  of the distribution of the sequence $(\bar X^n, I_n(Z),\tau^*_n)_{n\ge 1}$.
For every $(\alpha,u)\!\in I\!\!D([0,T], \R)^2\times [0,T]$, the canonical process is defined by  $\Xi_t(\alpha,u)= \alpha(t)=(\alpha^1(t),\alpha^2(t))\!\in \R^2$ and the canonical random times is given by $\theta(\alpha,u)=u$. Furthermore we will denote by $\Xi =(\Xi^1,\Xi^2)$ the two components of $\Xi$.

 Let 
\[
{\cal D}_t^{\theta}=\cap_{s>t}\sigma(\Xi_u, \{\theta\le u\}, \, 0\le u\le s\}\, \mbox{ if }\, t\!\in [0,T)\;\mbox{ and } {\cal D}_{T}^{\theta} = \sigma\big(\Xi_s, \{\theta\le s\}, \, 0\le s\le T\big)
\]
denote the canonical right-continuous filtration on $\Omega_c$. This canonical space $\Omega_c$ is equipped with  the product metric topology   $Sk^{\otimes 2}\otimes |.|$ where $|.|$ denotes 
the standard topology on $[0,T]$ induced by the absolute value.

\medskip In order to conclude to the convergence of the  {\em r\'eduites}, we need, following Theorem~\ref{ContRed} established in~\cite{LAPA}, to show that any limiting distribution $\Q= \lim_n \P_{((\bar X^n, I_n(Z)),\tau^*_n)}$  on the canonical space $\big(I\!\!D([0,T], \R^2)\times[0,T], Sk^{\otimes 2}\otimes |.|\big)$ satisfies  the $({\cal H})$-assumption, namely
\[
\E_{\Q} \big(H\,|\, {\cal D}_t^{\theta}\big) = \E_{\Q}\big (H\,|\, {\cal D}_t\big)\; \Q\mbox{-}a.s.
\]
for every random variable $H$ defined on $\Omega_c$.

\medskip
Let ${\rm Atom}_{\Q}(\theta)= \{s\!\in [0,T],\, \Q_{\theta}(\{s\})>0\}$ be the set, possibly empty, of $\Q$-atoms of $\theta$. Let $\Phi:I\!\!D([0,T], \R^2)\to \R$ and $\Psi:\D([0,T], \R)\to \R$ two  bounded functionals, $Sk^{\otimes 2}$- and $Sk$-continuous respectively and let $u\!\notin {\rm Atom}_{\Q}(\theta)$, $u\le s\le T$. Noting that $\Psi(I_n(Z)^s)\mbox{\bf 1}_{\{\tau^*_n\le u\}}$  is  ${\cal F}^n_s$-measurable, we get
\begin{eqnarray*}
\E_{\Q}\big(\Phi(\Xi)\Psi(\Xi^{2,s})\mbox{\bf 1}_{\{\theta\le u\}}\big)&=& \lim_n \E\Big(\Phi(\bar X^n,I_n(Z))  \Psi(I_n(Z)^s)\mbox{\bf 1}_{\{\tau^*_n\le u\}}\Big)\\
&= &  \lim_n \E\Big(\E\big[\Phi(\bar X^n,I_n(Z))|{\cal F}^Z_s\big]  \Psi(I_n(Z)^s)\mbox{\bf 1}_{\{\tau^*_n\le u\}} \Big). 
\end{eqnarray*}

Up to an extraction $(n')$,  we may assume that $\E\big[\Phi(\bar X^{n'},I_{n'}(Z))|{\cal F}^Z_s\big]$ weakly converges to $\E\big[\Phi(X,Z)|{\cal F}^Z_s\big]$ since $\Phi(\bar X^{n'},I_{n'}(Z))$ weakly converges  toward $\Phi(X,Z)$. Up to a second  extraction, still denoted $(n')$, we may assume that  $\Psi(I_n(Z)^s)$  $a.s.$ converges toward $\Psi(Z^s)$ for the Skorokhod topology since $\P(\Delta Z_s\neq 0)=0$ (the stopping operator at time $s$,  $\alpha\mapsto \alpha^s$,  is $Sk$-continuous at functions $\alpha$ which are continuous at $s$).

Consequently, going back on the canonical space $\Omega_c$, we obtain
\[
\Big(\E\big[\Phi(\bar X^n,I_n(Z))|{\cal F}^Z_s\big], \, \Psi(I_n(Z)^s),\,\mbox{\bf 1}_{\{\tau^*_n\le u\}}\Big)\stackrel{{\cal L}}{\longrightarrow} {\cal L}_{\Q}\left(\E\big[\Phi(\Xi)|{\cal F}^{\Xi^2}_s\big], \Psi(\Xi^{2,s}),\mbox{\bf 1}_{\{\theta \le u\}} \right).
\]
which ensures that
\[
\E_{\Q}\left( \Phi(\Xi) \psi(\Xi^{2,s})\mbox{\bf 1}_{\{\theta\le u\}}  \right) = \E_{\Q}\left( \E_{\Q}\big[\Phi(\Xi)|{\cal D}_{s_-}\big]\Psi(\Xi^{2,s}) \mbox{\bf 1}_{\{\theta\le u\}} \right).
\]

One concludes by  standard functional monotone approximation arguments  that the equality holds true for any bounded measurable functional $\Phi$, $\Psi$ and any $u\!\in [0,T]$. Then, by considering a sequence $s_n \downarrow s$, $s_n>s$, we derive that  
\[
\E_{\Q}\left( \Phi(\Xi) \,|\, {\cal D}^{\theta}_s\right) =  \E_{\Q} \Big(\Phi(\Xi)\,|\,{\cal D}_{s}\Big).
\]
This shows that the $({\cal H})$-assumption is fulfilled so that by Theorem~\ref{ContRed}, $U^n_0$ converges toward $U_0$.$\cqfd$

\small

\appendix
\section{Appendix: Euler scheme for Brownian martingale   diffusions}\label{App:A}

\begin{Pro}\label{pro:EulerWeakCv}
Let $(\bar X^n_t)_{t\in[0,T]}$ be the genuine Euler scheme of step $\frac Tn$ of the $SDE$ $dX_t=\sigma(t,X_t)dW_t, \, X_0=x$ defined as the solution to
 \[
 d\bar X^n_t =\sigma(\underline{t}_n, \bar X^n_{\underline{t}_n})dW_t,\; \bar X^n_0=x.
 \]
If $\sigma:[0,T]\times \R\to \R$ is continuous  and satisfies the linear growth assumption
\[
\forall\, t\!\in [0,T],\; \forall\, x\!\in \R,\quad |\sigma(t,x)|\le C_{\sigma}(1+|x|)
\]
Then the sequence $(\bar X^n)_{n\ge 1}$ is $C$-tight on ${\cal C}([0,T], \R)$ and any of its limiting distribution is a weak solution to the above $SDE$. In particular if a weak uniqueness assumption holds, then $\displaystyle \bar X^n \stackrel{(\|\,.\,\|_{\sup}) }{\longrightarrow} X$. 
\end{Pro}
Following $e.g.$~\cite{BOLE} (Lemma B.1.2, p.275, see also~\cite{KLPL, PAGtxt}), we first show that, owing to the linear growth assumption $|\sigma(t,x)|\le C_{\sigma}(1+|x|)$ made on $\sigma$, the non-decreasing function $\varphi_{p,n}(t)= \E \sup_{s\in [0,t]} |\bar X^n_s|^p $, $p\!\in [1,+\infty)$ is finite for every $t\!\in [0,T]$. Using Doob's Inequality and Gronwall's Lemma, it follows that there exists a real constant $C=C'_{p,\sigma}>0$ such that 
 \[
 \varphi_{p,n}(t) \le \varphi_p(t) :=C e^{C t}(1+|x|^p).
 \]
 Consequently, it follows from the the successive application of the $L^p$-B.D.G. and H\"older inequalities  for $p\!\in (2,+\infty)$ that, for every $s,t\!\in [0,T]$, $s\le t$,
 \begin{eqnarray*}
 \E |\bar X^n_{t}-\bar X^n_s|^p &\le & c_p^p \E \left(\int_s^t|\sigma(\underline{u}_n,\bar X^n_{\underline{u}_n})|^2du\right)^{\frac p2}\\
 &\le& c_p^p|t-s|^{\frac p2}\big(1+\varphi_p(T)\big).
  \end{eqnarray*}
  Kolmogorov's criterion (see~\cite{BIL}, Theorem 12.3, p.95) implies that the sequence $M_n= (W_t,\bar X^n_t)_{t\in [0,T]}$ is $C$-tight ($i.e.$  tight as $({\cal C}([0,T], \R), \|\,.\,\|_{\sup})$-valued random variables). From now on, we mainly rely on the results established in~\cite{JaMePa}. Let $n'$ be a subsequence such that $(\bar X^{n'},W)$ functionally weakly converges to a probability $\Q$ on  $({\cal C}([0,T], \R^2), \|\,.\,\|_{\sup})$, hence it  satisfies the  $U.T.$ (for {\em Uniform Tightness}) assumption (see Proposition~3.2 in~\cite{JaMePa}). 
  The function $\sigma$ being continuous on $[0,T]\times \R$, the sequence $(\sigma(\underline{t}_n, \bar X^n_{\underline{t}_n}))_{n\ge1}$ is $C$-tight on the Skorokhod space since $\big((\underline{t}_n, \bar X^n_{\underline{t}_n})_{t\in [0,T]}\big)_{n\ge 1}$ clearly is. One derives that, up to  a new extraction still denoted $(n')$, we may assume that $\big(\sigma(\underline{t}_{n'}, \bar X^{n'}_{\underline{t}_{n'}})_{t\in [0,T]}, \bar X^{n'}, W\big)_{n\ge1}$ functionally converges toward a probability $\P$ on $I\!\!D([0,T], \R^3)$.  By Theorem~2.6  from~\cite{JaMePa}  
for the functional convergence of stochastic integrals, we know that 
\[
\left( \sigma(\underline{t}_{n'}, \bar X^{n'}_{\underline{t}_{n'}}), (\bar X_t ^{n'},W _t), \int_0^t \sigma(\underline{s}_{n'}, \bar X^{n'}_{\underline{s}_{n'}})dW_s\right)_{t\in [0,T]}\stackrel{{\cal L}(Sk)}{\longrightarrow}  \Q\quad\mbox{ as }\quad n\to +\infty
\]
where $\Q$ is a probability distribution on $I\!\!D([0,T], \R^4)$ such that the canonical process $Y=(Y^i)_{i=1:4}$ satisfies  $Y\stackrel{{\cal L}}{\sim}\big(Y^1,(Y^2, B), \int_0^. Y^2_{s}dB_s)$ where $ B:Y^3$ is a standard $\Q$-Brownian motion with respect to the $\Q$-completed right continuous canonical  filtration $({\cal D}^4_t)_{t\in [0,T]}$ on $I\!\!D([0,T], \R^4)$.  Furthermore, we know that $Y ^1= \sigma(.,Y^2  )$ $\Q$-$a.s.$ since $\sup_{t\in [0,T]}|\sigma(\underline{t}_{n'}, \bar X^{n'}_{\underline{t}_{n'}})-\sigma(t, \bar X^{n'}_t)|$  converges to $0$ in probability. The former claim follows from the facts that   $\sup_{t\in [0,T]}|\bar X^n_t|$ is tight and $\sigma(t,\xi)$ is uniformly continuous on every compact 
set of $[0,T]\times \R$ with linear growth in $\xi$ uniformly in $t\!\in [0,T]$. On the other hand, we know that $\bar X^{n'}=x+\int_0^{.}  \sigma(\underline s_{n'}, \bar X^{n'}_{\underline{s}_{n'}})dW_s$ which in turn implies that $Y_{2,.} = x + \int_0^. \sigma(s, Y_{2,s})dW_s$. This shows the existence of a weak solution to the $SDE$ $X_t= x+\int_0^t \sigma(s,X_s)\,dW_s, \, t\!\in [0,T]$.

Under the weak uniqueness assumption, this distribution is unique hence is the only functional weak limiting  distribution for the tight sequence $(\bar X^n)_{n\ge 1}$, hence we get the convergence in distribution on ${\cal C}([0,T], \R)$. ~$\cqfd$

\medskip
\noindent {\bf Remark.} If the original $SDE$ has a unique strong solution, the same proof leads to the establish the convergence  in probability of the Euler scheme toward $X$. One just has to add the process $X$ itself to the sequence $\big((\sigma(\underline{t}_{n}, \bar X^{n}_{\underline{t}_{n}}))_{t\in [0,T]}, \bar X^n, W\big)_{n\ge1}$

\section{Appendix: Euler scheme for a L\'evy driven martingale diffusion}\label{app:B}
We consider the following $SDE$ driven by a martingale L\'evy process $Z$ with L\'evy measure $\nu$
\begin{equation}\label{eq:Levy}
  X _{t }=x +\int_{(0,t]} \kappa(s,  X_{s_-})dZ_s,  \; X_0 =x
\end{equation}
(where $\kappa$ is a Borel function on $[0,T]\times \R$) and its {\em genuine} Euler scheme defined by 
\begin{equation}\label{eq:EulerLevy}
\bar X^n_{t_{k+1}}= \bar X^n_{t_k}  +\  \kappa(t_k, \bar X_{t_k}) (Z_{t_{k+1}}-Z_{t_k}), \; k=1,\ldots,n,\;\bar X_0=X_0=x
\end{equation}
at discrete times $t^n_k$ and extended into a continuous time process  by setting  $\bar X^t= \bar X^n _{\underline t_n}$ so that
\begin{equation}\label{eq:EulerLevyb}
\bar X^n_t = x+\int_{(0,t]} \kappa(\underline{s}_{n-},\bar X^n _{\underline s_n-})dZ_s.
\end{equation}

\subsection{Convergence of the Euler scheme toward a solution to L\'evy driven $SDE$}
 \begin{Pro}\label{pro:WeakEulerLevy} $(a)$ Assume that $\nu(|z|^p)<+\infty$ for a $p\!\in(1,2]$, has no Brownian component and $\kappa(t,\xi)$ has linear growth in $\xi$ uniformly in $t\!\in [0,T]$. Then
 \[
 \sup_{n\ge 1}\big\|\sup_{t\in [0,T]} |\bar X^n_t|\big\|_p+\big\|\sup_{t\in [0,T]}| X_t|\big\|_p<+\infty.
 \]  
If moreover $\kappa$ is continuous, the $SDE$~\eqref{LevyEDS} has at least one weak  solution and, as soon as weak uniqueness holds for~\eqref{eq:Levy}, 
  \[
  \bar X^n_t\stackrel{{\cal L}(Sk)}{\longrightarrow} X.
  \]
 \noindent $(b)$ If $\nu(z^2)<+\infty$, the same  result remains true {\em mutatis mutandis} if $Z$ has a non-zero Brownian component. 
 \end{Pro}
 
\noindent  {\bf Remark.} In fact, as soon as~\eqref{eq:Levy} has a strong solution; one shows by the same argument as those developed below the stronger result
\[
\sup_{t\in [0,T]}|\bar X^n_t-X_t| \stackrel{\P}{\longrightarrow} 0\; \mbox{ as }\; n\to +\infty.
\]
We refer to~\cite{JAC} for a proof when $\kappa$ is homogeneous and $C^3$.

\bigskip
\noindent {\bf Proof.} $(a)$ We consider the L\'evy-Khintchine decomposition   of the L\'evy process  $Z=(Z_t)_{t\in [0,T]}$, namely   
$$
Z_t =  \widetilde Z_t+ Z^{1}, \; t\!\in [0,T],
$$ 
where $\widetilde Z^{\eta}$ is a pure jump  {\em square integrable} martingale with  jumps of size at most $1$ and  L\'evy measure $\nu(\,.\cap\{|z|\le 1\})$ (having a second moment by construction) and $Z^{1}$ is a compensated (hence martingale) Poisson process with (finite) L\'evy measure $\nu(\,.\cap\{|z| >1\})$.

It is clear from~\eqref{eq:EulerLevy}  that $\bar X^n_{t^n_k}\!\in L^p$ for every $k=0\ldots,n$. Then, as $\nu(|z|^p)<+\infty$,  it follows classically that 
$\sup_{u\in [t^n_k,t^n_{k+1}]}|Z_u| \stackrel{d}{\sim}\sup_{[0,\frac Tn]} |Z_u| \!\in L^p$ (see $e.g.$~\cite{SATO}). Combining these two results implies that $\varphi_{p,n}(t):= \big\|\sup_{s\in [0,t]} |\bar X^n_s|\big\|_p$ is finite for every $t\!\in [0,T]$.

It follows from Equation~\eqref{eq:EulerLevyb} satisfied by $\bar X$ that, for every $t\!\in [0,T]$,
\[
\sup_{s\in [0,t]}|\bar X^n_s|\le |x|+\sup_{s\in [0,t]}\left|\int_{(0,s]}\kappa(\underline{u}_{n-},\bar X^n_{\underline{u}_{n-} })dZ_u\right |.
\]
Consequently
\[
\varphi_{p,n}(t)\le |x| + \left\|\sup_{s\in [0,t]}\left|\int_{(0,s]}\kappa(\underline{u}_{n-},\bar X^n_{\underline{u}_{n-} })dZ_u\right| \right\|_p
\]
The $L^p$-B.D.G Inequality implies  (since $p>1$)  
\[
 \left\|\sup_{s\in [0,t]}\left|\int_{(0,s]}\kappa(\underline{u}_n,\bar X_{\underline{u}_{n} -})dZ_u\right|\right\|_p\le c_p
  \left\| \sum_{0<s\le t}\kappa(\underline{s}_n,\bar X_{\underline{s}_{n-}})^2(\Delta Z_s)^2 \right\|_{\frac p2}^{\frac 12}.
\]
As $p\!\in (1,2]$, $\frac p2\le 1$ which implies
\begin{eqnarray*}
\left\| \sum_{0<s\le t}\kappa(\underline{s}_n,\bar X_{\underline{s}_{n-}})^2(\Delta Z_s)^2 \right\|_{\frac p2}^{\frac 12}&\le &\left(\E \sum_{0<s\le t}|\kappa(\underline{s}_n,\bar X_{\underline{s}_{n-}})|^p|\Delta Z_s|^p\right)^{\frac 1p}=\left(\nu(|z|^p)\, \E \!\int_0^t| \kappa(\underline{s}_n,\bar X_{\underline{s}_{n-}})|^p ds\right)^{\frac 1p}\\
&\le& C_{\kappa,p}^p \nu(|z|^p)^{\frac 1p}\left(\int_0^t (1+\varphi(s)^p)ds\right)^{\frac 1p}
\end{eqnarray*}
where $C_{\kappa,p}$ is a real constant satisfying $|\kappa(s,\xi)|\le C(1+|\xi|^p)^{\frac 1p}$, $(s,\xi)\!\in [0,T]\times \R$.  

Finally, there exists a real constant $C'= C'_{\kappa,p, \nu}$ such that the function  $\varphi_{p,n}$ satisfies 
\[
\varphi_{p,n}(t) ^p\le C'\left(|x|^p + t+ \int_0^t  \varphi(s)^p ds\right)
\]
which in turn implies by Gronwall's Lemma
\[
\forall\, t\!\in [0,T], \quad \varphi_{p,n}(t)^p \le e^{C't}C'(1+|x|^p)
\]
or, equivalently,
\[
\forall\, t\!\in [0,T], \quad \varphi_{p,n}(t) \le \varphi(t) = e^{C''t}C''(1+|x|)\quad \mbox{ where} \quad C''= C'/p.
\]

To establish the Skorokhod tightness of the sequence $(\bar X^n)_{n\ge 1}$, we rely on the Aldous tightness criterion (see Definition~\ref{def:DetAldous}$(b)$ or~\cite{JASH2},  Theorem 4.5, p.356).  Let $\rho\!\in (0,1]$. 
Let $\sigma$ and $\tau $ be two $[0,T]$-valued ${\cal F}^Z$-stopping stopping times such that $\sigma\le \tau \le (\sigma+\delta)\wedge T$. 
 \begin{eqnarray*}
\E |\bar X^n_{\tau}-\bar X^n_{\sigma}|^{\rho}=\E\, \Big|\sum_{\sigma <u\le\tau } \kappa(\underline{u}_n, \bar X^n_{\underline{u}_{n-}})\Delta Z_u\Big|^{\rho}&\le& \E \Big(\sum_{\sigma <u\le\tau } |\kappa(\underline{u}_n, \bar X^n_{\underline{u}_{n-}})|^{\rho}|\Delta Z_u|^{\rho}\Big)\\
&=& \nu(|z|^{\rho})\E \int_{\sigma}^{(\sigma+\delta)\wedge T} |\kappa(\underline{u}_n, \bar X^n_{\underline{u}_{n-}})|^{\rho}\\
&\le&  \delta\,\nu(|z|^{\rho})\, \E \sup_{t\in [0,T]} |\kappa(t,\bar X^n_t)|^{\rho}\\
&\le& \delta\,\nu(|z|^{\rho})\,C_{\kappa}(1+\varphi_p(T))^{\frac{\rho}{p}}
   \end{eqnarray*}
where we used that $\rho\le 1\le p$ and $\nu(|z|^{\rho})\le \nu(|z|^{2}\wedge 1)+\nu(|z|^{p})<+\infty$.
%
%
%
%
%
Then 
\[
\sup\left\{\E |\bar X_{\tau}-\bar X_{\sigma}|^{\rho} + \E | Z_{\tau}- Z_{\sigma}|^{\rho},\;\sigma\le \tau\le (\sigma+\delta)\wedge T, \,{\cal F}^Z\mbox{-stopping times} \right\}\le  \nu(|z|^{\rho})(1+C_{\kappa}(1+\varphi_p(T))^{\frac{\rho}{p}} \big)\delta
\] 
which goes to $0$ as $\delta\to 0$. This  implies that the sequence $M_n = (\bar X^n, Z)$, $n\ge 1$, is  $Sk$-tight. 
%
Moreover, following Proposition~3.2  from~\cite{JaMePa},  the sequence $(M_n)_{n\ge 1}$ satisfies the $U.T.$  condition  it is $Sk$-tight and
 \begin{eqnarray*}
\E \sup_{t\in [0,T] }\big(|\Delta \bar X^n|\vee |\Delta Z_t| \big)&\le &\Big[ \E \Big(\sum_{0<t\le T} |\Delta \bar X_t^n|^p+ |\Delta Z_t |^p   \Big)\Big]^{\frac 1p}\\
&\le &\Big[ \nu(|z|^p)\E \int_0^T \big(1+|\kappa(\underline{t}_n, \bar X^n_{\underline{t}_n})|^p\big) dt  \Big]^{\frac 1p}\\
&\le& \big(\nu(|z|^p) \big)^{\frac 1p}\big(T+C_{\kappa,p}^p(1+\varphi_p(T))\big)^{\frac 1p}<+\infty.
\end{eqnarray*}

\noindent On the other hand, the sequence $\left((\kappa\big({\underline{t}_n}, \bar X^n_{\underline{t}_n}))_{t\in [0,T]}, M_n\right)_{n\ge 1}$ is $Sk$-tight owing to the following lemma.

\begin{Lem} Let ${\cal V}^+_{[0,T]}$ be the set of functions $\mu :[0,T]\to [0,T]$ such that $\mu(0)=0$ and $\mu(T)=T$ endowed with the sup norm. Assume $\kappa:[0,T]\times \R\to \R$ is continuous. Then  the mapping $\Psi:{\cal V}^+_{[0,T]}\times I\!\! D([0,T], \R^{d}) \to I\!\! D([0,T], \R^{1+d})$ defined by $\Psi(\mu, \alpha) = \big(\kappa(\mu(.), \alpha^1(.)), \alpha\big)$ is  continuous ($\alpha=(\alpha^1,\ldots,\alpha^d)$) for the product topology.
\end{Lem}

\noindent {\bf Proof (of the lemma).} Let $(\lambda_n)_{n\ge 1 }$ be a sequence of homeomorphisms of $[0,T]$ such that  $\lambda_n \to Id_{[0,T]}$ and $\alpha_n\circ \lambda_n \to \alpha$ uniformly and let $\mu_n \to \mu$ in ${\cal V}^+_{[0,T]}$ where $Id_{[0,T]}$ denotes the identity of $[0,T]$ to $[0,T]$. Then the closure of $(\alpha_n\circ \lambda_n(t))_{n\ge 1, t\in [0,T]}$ is a compact set $K$ of $\R^d$ hence the function $\kappa$ is uniformly continuous on $[0,T] \times K$. On the other hand
\[
\|\mu_n\circ \lambda_n -Id_{[0,T]}\|_{\sup} \le \|\mu_n -Id_{[0,T]}\|_{\sup} + \| \lambda_n-Id_{[0,T]}\|_{\sup} \quad \mbox{ as }\quad n\to +\infty
\] 
and $\displaystyle 
\|\alpha_n \circ \lambda_n -\alpha \|_{\sup}\to 0 \quad \mbox{ as }\quad n\to +\infty$. 
The conclusion follows.~$\cqfd$

\bigskip
Up to an extraction, we may assume that the triplet $\big(\big(\kappa(\underline{t}_{n'}, \bar X_{\underline{t}_{n'}}^{n'} \big)_{t\in [0,T]}, M_{n'}\big)_{n\ge 1}$ weakly converges for the Skorokhod topology toward a probability $\P$ on the canonical Skorokhod  space $(I\!\! D([0,T], \R^3),( {\cal D}_t)_{t\in [0,T]})$. 

By Theorem~2.6  from~\cite{JaMePa}   for the functional convergence of stochastic integrals, we know that 
\[
\Big(\kappa(\underline{t}_{n'}, \bar X^{n'}_{\underline{t}_{n'}}), (\bar X ^{n'}_t,Z_t ), \int_0^t \kappa(\underline{s}_{n'-}, \bar X^{n'}_{\underline{s}_{n'-}})dZ_s\Big)_{t\in [0,T]}\stackrel{{\cal L}(Sk)}{\longrightarrow}  \Q
\]
probability distribution on $I\!\!D([0,T], \R^4)$ such that the canonical process $Y=(Y^i)_{i=1:4}$ satisfies  $Y\stackrel{{\cal L}}{\sim}\big(Y^1,(Y^2, Z), \int_0^. Y^2_{s}dZ_s)$ where $Y^3$ is a L\'evy process  with respect to the $\Q$-completed right continuous canonical  filtration $({\cal D}^{\Q}_t)_{t\in [0,T]}$ on $I\!\!D([0,T], \R^4)$ having the distribution of $Z$ ($i.e.$  $\Q_{Y^3}= {\cal L}(Z)$).  Furthermore, we know that $Y^1= \kappa(.,Y^2{.})$ $\Q$-$a.s.$ since the mapping $(\mu, (\alpha^i)_{i=1:4})\mapsto \alpha^1-\kappa(\mu, \alpha^2)$ is continuous from ${\cal V}^+_{[0,T]}\times I\!\!D([0,T], \R^4)$ to $I\!\!D([0,T], \R)$ (and $\underline t_n$ converges uniformly to $Id_{[0,T]}$). 

On the other hand we know that $\bar X^{n'}_t=x+\int_0^{t}  \kappa(\underline s_{n'-}, \bar X^{n'}_{\underline{s}_{n'-}})dZ_s, \, t\!\in [0,T]$ which in turn implies that $(Y_t^2   = x + \int_0^t \kappa(s, Y^2_{s_-})dZ_s,\, t\!\in [0,T])$ $\Q$-$a.s.$. This shows the existence of a weak solution to the $SDE$ $X_t= x+\int_0^t \kappa(s,X_{s_-})dZ_s,\, t\!\in [0,T] $.

Under the weak uniqueness assumption,   the distribution $\Q_{Y^2}$ of $Y^2$ is unique equal, say, to~$\P_{X}$.

\medskip
\noindent $(b)$  We assume that the L\'evy measure has a finite second moment $\nu(z^2)<+\infty$ on the whole real line. Then one can decompose $Z$ as 
\[
Z_t = a\,W_t + \widetilde Z_t,\; t\!\in [0,T], \quad (a\ge 0)
\]
where $\kappa\ge 0$ and $\widetilde Z$ is a pure jump martingale L\'evy process with L\'evy measure $\nu$. Then one shows like in the Brownian case that $\varphi(t)= \E \sup_{s\in [0,t]}|\bar X^n_s|^2$ is finite over $[0,T]$ using that all $\bar X_{t_k}$ are square integrable and $\E \sup_{s\in [t_k,t_{l+1})}|Z_s-Z_{t_k}|^2 = \E \sup_{s\in  [0, \frac Tn]} |Z_s|^2  <+\infty$. Then, using Doob's Inequality, we show that 
\[
\varphi(t)\le 4C^2_{\kappa}(a^2+\nu(z^2)\big) \left(t+\int_0^t \varphi(s)ds\right)
\]
 where $C_{\kappa} $ is a real constant satisfying  $\kappa(t,\xi)\le C_{\kappa}(1+|\xi|^2)^{\frac 12}$, $\xi\!\in \R$. 

 To establish the Skorokhod tightness of the sequence, we rely on the Aldous tightness criterion (see Definition~\ref{def:DetAldous}$(b)$ or\cite{JASH2},  Theorem 4.5, p.356).  Let $\sigma$ and $\tau $ be two $[0,T]$-valued ${\cal F}^Z$-stopping stopping times such that $\sigma\le \tau \le (\sigma+\delta)\wedge T$. Then applying Doob's Inequality, this time   to the martingale $\Big(\int_{\sigma}^{\sigma+s}\kappa(\underline{u}_n, \bar X_{\underline{u}_{n-}})dZ_u\Big)_{s\ge 0}$, we get
\vskip -0.5cm 
 \begin{eqnarray*}
\E |\bar X_{\tau}-\bar X_{\sigma}|^2&\le& 4a^2 \E\left(\int_{\sigma}^{\tau} |\kappa(\underline{u}_n, \bar X^n_{\underline{u}_{n-}})|^2 du\right)+4\,\E\left(\sum_{\sigma<u\le\tau} |\kappa(\underline{u}_n, \bar X^n_{\underline{u}_{n-}})|^2 |\Delta Z_u|^2\right)\\
&=& 4\big(a^2+\nu(z^2)\big) \E\left(\int_{\sigma}^{\tau} |\kappa(\underline{u}_n, \bar X^n_{\underline{u}_{n-}})|^2 du\right)\\
&\le & 4\big(a^2+\nu(z^2)\big)  \E\left(\int_{\sigma}^{(\sigma + \delta)\wedge T} |\kappa(\underline{u}_n, \bar X^n_{\underline{u}_{n-}})|^2 du\right)\\
&\le & 4 (a^2+\nu(z^2)\big) \delta \,C_{\kappa}^2(1+\varphi(T)).
   \end{eqnarray*}

Then $\displaystyle \E |\bar X_{\tau}-\bar X_{\sigma}|^2 + \E |\bar Z_{\tau}-\bar Z_{\sigma}|^2\le 4 (a^2+\nu(z^2)\big)\nu(z^2)(1+\varphi(T) \big)\delta$ which clearly implies the $Sk$-tightness of the sequence $M_n = (\bar X^n, Z)$, $n\ge 1$. 

The sequence satisfies the $U.T.$  condition from~\cite{JaMePa} since $(M_n)_{n\ge 1}$ is $Sk$-tight and (see Proposition~3.2  from~\cite{JaMePa})
 \begin{eqnarray*}
\E \sup_{t\in [0,T] }\big(|\Delta \bar X^n|\vee |\Delta Z_t| \big)&\le &\Big( \E \sum_{0<t\le T} |\Delta \bar X_t^n|^2+ |\Delta Z_t |^2   \Big)^{\frac 12}\\
&\le &\Big( \nu(z^2)\E \int_0^T (1+|\kappa(\underline{t}_n, \bar X^n_{\underline{t}_n})|^2 dt  \Big)^{\frac 12}\\
&\le& \nu(z^2) (1+C_{\kappa}(1+\varphi(T)))<+\infty.
\end{eqnarray*}
From this point, the proof is quite similar to that of claim $(a)$.~$\cqfd$

\subsection{Higher moments}
Let $Z_t = a W_t +\widetilde Z_t$, $t\!\in [0,T]$, be the decomposition of the L\'evy process $Z$ where $W$ is a standard B.M. and $\widetilde Z$ is an independent pure jump L\'evy process.
\begin{Pro}\label{pro:highmoment}If $\nu(|z|^p)<+\infty$ for $p\!\in [2,+\infty)$, then
\[
\sup_{n\ge 1}\Big\| \sup_{t\in [0,T]}|\bar X^n_t|\Big\|_p <+\infty.\]
\end{Pro}

\noindent {\bf Proof.} If $p\!\in (1,2]$, the claim follows from the above Proposition~\ref{pro:WeakEulerLevy}. Assume from now on $p\!\in [2,+\infty)$. Let $\varphi_{p,n}(t) = \E \big(\sup_{t\in [0,T]}|\bar X^n_t|^p\big)$. Let $\ell_p$ be the unique integer defined by the inequality $2^{\ell_p}< p\le 2^{\ell_p+1}$. It is straightforward using the same arguments as above that $\varphi_{p,n}(T)<+\infty$ since  $\sup_{t\in [0,T]}|Z_t|^p \!\in L^1$ (see~\cite{SATO}, Theorem 25.18, p.~166) and $X_{t_k}\!\in L^p$ by induction using~\eqref{eq:EulerLevy}.  For convenience, we denote $\kappa_{s_-} = \kappa(\underline{s}_n, \bar X^n_{\underline{s}_n-})$.

Now,  combining the (integral and regular) Minkowski and the B.D.G. Inequalities implies 
\begin{eqnarray}
\nonumber \varphi_{p,n}(t)^{\frac 1p} &\le &|x| +c_p \Big \|a^2\int_0^t \kappa_{s_-}^2ds+ \sum_{0<s\le t} \kappa_{s_-}^2 (\Delta Z_s)^2\Big \|^{\frac 12}_{\frac p2}\\
\label{eq:varyhipn}&\le&  |x| +c_p \Big(a\Big\|\int_0^t \kappa_{s_-}^2ds \Big\|^{\frac 12}_{\frac p2}+ \Big\|\sum_{0<s\le t} \kappa_{s_-}^2(\Delta Z_s)^2\Big\|_{\frac p2}^{\frac 12}\Big)
\end{eqnarray}
where we used in the second inequality that $\sqrt{u+v}\le \sqrt{u}+\sqrt{v}$, $u,\, v\ge 0$. First note that by two successive applications of H\" older Inequality to $dt$ and $d\P$,  we obtain
\begin{equation}\label{eq:BDGW}
\left \|\int_0^t \kappa_{s_-}^2ds \right\|^{\frac 12}_{\frac p2}  \le   T^{\frac 12-\frac 1p}\left(\int_0^t \E\, |\kappa_{s_-}|^p ds\right)^{\frac 1p}
\end{equation}

\smallskip Now using that for every  $\ell\!\in \{1,\ldots,\ell_p\}$,  $\displaystyle \sum_{0<s\le t} |\kappa_{s_-}|^{2^{\ell}}|\Delta Z_s|^{2^{\ell}}-\int_0^t  |\kappa_{s_-}|^{2^{\ell}}ds\, \nu(|z|^{2^{\ell}})$, $ t\!\in [0,T]$, is a true martingale, we have by combining the Minkowski inequality, the  B.D.G. Inequality applied with  $\frac {p}{2^{\ell}}>1$ and the elementary inequality $(u+v)^{r} \le u^r+v^r$, $u$, $v\ge 0$, $r\!\in (0,1]$, yield 
\begin{eqnarray*}
\left\| \sum_{0<s\le t }  |\kappa_{s_-}|^{2^{\ell}}(\Delta Z_s)^{2^{\ell}}\right\|^{\frac{1}{2^{\ell}}}_{\frac{p}{2^{\ell}}}&\le& \left\|\sum_{0<s\le t }  |\kappa_{s_-}|^{2^{\ell}}(\Delta Z_s)^{2^{\ell}}-\int_0^t   |\kappa_{s_-}|^{2^{\ell}}ds\, \nu(|z|^{2^{\ell}} )\right\|^{\frac{1}{2^{\ell}}}_{\frac{p}{2^{\ell}}}+\left\| \int_0^t   |\kappa_{s_-}|^{2^{\ell}}ds\right\|^{\frac{1}{2^{\ell}}}_{\frac{p}{2^{\ell}}} \nu(|z|^{2^{\ell}})^{\frac{1}{2^{\ell}}}\\
&\le& c_{\frac{p}{2^{\ell}}}^{\frac{1}{2^{\ell}}}\left\| \sum_{0<s\le t }  |\kappa_{s_-}|^{2^{\ell+1}}(\Delta Z_s)^{2^{\ell+1}}\right\|^{\frac{1}{2^{\ell+1}}}_{\frac{p}{2^{\ell+1}}} +\left\| \int_0^t   |\kappa_{s_-}|^{2^{\ell}}ds\right\|^{\frac{1}{2^{\ell}}}_{\frac{p}{2^{\ell}}} \nu(|z|^{2^{\ell}})^{\frac{1}{2^{\ell}}}.
\end{eqnarray*}

Then two applications of H\"older Inequality applied to $dt$ and  $d\P$ successively imply 
\[
\left\| \int_0^t   |\kappa_{s_-}|^{2^{\ell}}ds\right\|^{\frac{1}{2^{\ell}}}_{\frac{p}{2^{\ell}}} \le T^{\frac{1}{2^{\ell}}-\frac 1p}\left(\int_0^t \E\, |\kappa_{s_-}|^p ds\right)^{\frac 1p}.
\]
Summing up these inequalities in cascade finally yields a positive real constant $K^{(0)}=K_{p, \nu,a,T}$ such that
\[
 \left\|\sum_{0<s\le t} |\kappa_{s_-}|^2(\Delta Z_s)^2\right\|_{\frac p2}^{\frac 12} \le  K^{(0)} \left(\left(\int_0^t \E\, |\kappa_{s_-}|^p ds\right)^{\frac 1p} +  \left\|\sum_{0<s\le t}| \kappa_{s_-}|^{2^{\ell_p+1}}(\Delta Z_s)^{2^{\ell_p+1}}\right\|_{\frac{p}{2^{\ell_p+1}}}^{\frac{1}{2^{\ell_p+1}}} \right).
\]
Now, as $\frac{p}{2^{\ell_p+1}}\le 1$, one gets by the compensation formula 
\[
\left \|\sum_{0<s\le t} |\kappa_{s_-}|^{2^{\ell_p+1}}|\Delta Z_s|^{2^{\ell_p+1}}\right \|_{\frac{p}{2^{\ell_p+1}}}^{\frac{1}{2^{\ell_p+1}}} \le \left(\E\sum_{0<s\le t} |\kappa_{s_-}|^{p}(\Delta Z_s)^{p}\right)^{\frac 1p}= \left(\int_0^t \E |\kappa_{s_-}|^p ds\right)^{\frac 1p}\nu(|z|^p)^{\frac 1p}.
\]
Hence, there exists a real constant $K^{1)}= K^{(1)}_{p, \nu,a,T}>0$
\begin{equation}\label{eq:BDGLevy}
 \left\|\sum_{0<s\le t} |\kappa_{s_-}|^2(\Delta Z_s)^2\right\|_{\frac p2}^{\frac 12}\le  K^{(1)}_{p, \nu,a,T}\left(\int_0^t \E |\kappa_{s_-}|^p ds\right)^{\frac 1p}.
\end{equation}
Finally, plugging~\eqref{eq:BDGW} and~\eqref{eq:BDGLevy} in~\eqref{eq:varyhipn}, there exist positive  real constants $K^{(\ell)} = K^{(\ell)}_{p, \nu,a,T}$, $\ell=2,3$, such that  
\begin{eqnarray*}
\varphi_{p,n}(t)^{\frac 1p} &\le &K^{(2)}_{p, \nu,a,T}\left(|x|+  \left(\int_0^t \E |\kappa_{s_-}|^p ds\right)^{\frac 1p}\right)
\\
&\le & K^{(3)}_{p, \nu,a,T}\left(|x|+1+  \left(\int_0^t  \varphi_{p,n} (s) ds\right)^{\frac 1p}\right)
\end{eqnarray*}
(where we used in the second inequality that $\kappa$ has linear growth) so that 
\[
\varphi_{p,n}(t) \le 2^{p-1}(K'^{(3)}_{p, \nu,a,T})^p\Big(\big(|x|+1\big)^p + \int_0^t  \varphi_{p,n} (s) ds\Big).
\]
Gronwall's lemma completes the proof since it implies that
\[
\varphi_{p,n}(t) \le e^{ 2^{p-1}(K'^{(3)}_{p, \nu,a,T})^p\,t}  2^{p-1}(K'^{(3)}_{p, \nu,a,T})^p(|x|+1)^p.\;\cqfd
\]

\end{document}